\documentclass[11pt,letterpaper]{amsart}
\usepackage{amsmath,amssymb,amsfonts}
\usepackage{mathtools}
\usepackage{bbm}
\usepackage{enumerate}
\usepackage{amscd, latexsym, tikz-cd}
\usepackage[all]{xy}
\usepackage{pb-diagram}
\usepackage{mathrsfs}
\usepackage{pictexwd,dcpic}
\usepackage[colorlinks,final,backref=page,hyperindex]{hyperref}
\usepackage{makecell}

\usepackage{extarrows}

\usepackage{pb-diagram}
\theoremstyle{plain}

\newcounter{theorem}[section]
\numberwithin{theorem}{section}
\numberwithin{equation}{subsection}

\newtheorem{Thm}[theorem]{Theorem}
\newtheorem{Conj}[theorem]{Conjecture}
\newtheorem{Cor}[theorem]{Corollary}
\newtheorem{Prop}[theorem]{Proposition}
\newtheorem{Lem}[theorem]{Lemma}
\newtheorem{Que}[theorem]{Question}
\newtheorem{Exm}[theorem]{Example}

\theoremstyle{remark}
\newtheorem{Def}[theorem]{Definition}
\newtheorem{Rmk}[theorem]{Remark}

\newcommand{\Hom}{\operatorname{Hom}}
\newcommand{\Ext}{\operatorname{Ext}}

\newcommand{\GL}{\operatorname{GL}}

\newcommand{\SL}{\operatorname{SL}}
\newcommand{\SO}{\operatorname{SO}}
\newcommand{\rO}{\operatorname{O}}
\newcommand{\Sp}{\operatorname{Sp}}

\newcommand{\U}{\operatorname{U}}

\newcommand{\Ind}{\operatorname{Ind}}
\newcommand{\ind}{\operatorname{ind}}
\newcommand{\Jac}{\operatorname{Jac}}

\newcommand{\Irr}{\operatorname{Irr}}

\newcommand{\csupp}{\operatorname{csupp}} 
 
\newcommand{\Bes}{\operatorname{Bes}} 

\newcommand{\LS}{\operatorname{LS}} 
\newcommand{\MVW}{\operatorname{MVW}} 

\newcommand{\id}{\operatorname{id}}

\newcommand{\C}{\mathbb{C}}

\newcommand{\D}{\mathbb{D}}

\newcommand{\Z}{\mathbb{Z}}
\newcommand{\R}{\mathbb{R}}

\newcommand{\bm}{\begin{multline*}}
\newcommand{\tu}{\end  {multline*}}

\newcommand{\disc}{{\rm disc}}

\newcommand{\unit}{{\rm unit}}

\def\calA{\mathcal{A}}

\def\calC{\mathcal{C}}

\def\calH{\mathcal{H}}

\def\calJ{\mathcal{J}}

\def\calL{\mathcal{L}}
\def\calM{\mathcal{M}}

\def\calO{\mathcal{O}}

\def\calR{\mathcal{R}}
\def\calS{\mathcal{S}}

\def\ra{\rightarrow}
\def\lra{\longrightarrow}

\def\LP{\Pi^{\rm L}} 
\def\AP{\Pi^{\rm A}} 
\def\gp{{\rm gp}} 
\def\np{{\rm np}} 
\def\rel{{\rm rel}} 
\def\Der{{\rm Der}} 

\title[Non-tempered GGP conjecture for classical groups: the Bessel case]
{Non-tempered Gan--Gross--Prasad conjecture for classical groups: the Bessel case}
\author{Rui Chen, Jing Quan Chong, and Chuijia Wang} 

\address{School of Mathematical Sciences, Zhejiang University, 866 Yuhangtang Road, Hangzhou 310058, China }
\email{rchenmat@zju.edu.cn}

\address{Department of Mathematics, Faculty of Science, National University of Singapore, Level 4, Block S17, 10 Lower Kent Ridge Road, Singapore, 119076}
\email{e0360564@u.nus.edu}

\address{Department of Mathematics, The University of Hong Kong, Pokfulam, Hong Kong, China}
\email{wangchuijia@u.nus.edu}

\begin{document}

\begin{abstract}
In their recent paper, Gan, Gross and Prasad proposed a precise conjecture to describe the restrictions of representations of classical groups that lie in an L-packet of Arthur type. On top of local root numbers, the main notion considered is the relevance of two A-parameters. In this paper, we establish that the relevance of two A-parameters is sufficient for the existence of a Bessel period for a pair of representations in the associated L-packets. Moreover, such a pair is unique. We also prove the necessity of relevance for co-tempered A-parameters, and give several counterexamples in the general case.
\end{abstract}

\maketitle



\section{Introduction}

\subsection{Background: tempered branching laws}

A central problem in the representation theory of reductive groups is the branching laws, which aims to understand how an irreducible representation of a reductive group decomposes when restricted to a suitable subgroup. The Gan--Gross--Prasad problem (``GGP'' for short) is one of the most well-known instances of this type of problem. Starting from the pioneering work of Gross and Prasad on special orthogonal groups \cite{MR1186476, MR1295124} and its extension by Gan, Gross, and Prasad \cite{MR3202556}, this problem was placed into a unified local and global framework connecting branching laws with arithmetic invariants of Langlands parameters.

\vskip 5pt

Let us recall the local Bessel setting of the GGP problem, which is the primary focus of this paper. Let $E$ be a non-Archimedean local field of characteristic zero, $\iota$ an involution of $E$, and $F$ the subfield of $E$ fixed by \(\iota\). In particular, $E$ is either $F$ itself or a quadratic extension of $F$. Consider a pair of non-degenerate quadratic (when $E=F$) or Hermitian (when $E \neq F$) spaces $W \subset V$ over $E$ such that the orthogonal complement $W^\perp$ is odd-dimensional and split, and let $G\left(V\right)$ and $G\left(W\right)$ denote the corresponding orthogonal or unitary groups. Using this embedding $W \subset V$, we can define the Bessel subgroup $\Bes\left(V,W\right) \subset G\left(V\right) \times G\left(W\right)$ and a generic character $\nu$ of it as in \cite[Sect. 12]{MR3202556}. Given irreducible representations $\pi$ of $G\left(V\right)$ and $\sigma$ of $G\left(W\right)$, one is interested in understanding the Bessel multiplicity
\[
    m\left(\pi,\sigma\right) \coloneqq \dim \Hom_{\Bes\left(V,W\right)}\left(\pi \boxtimes \sigma, \nu\right),
\]
which is known to be at most $1$ \cite{MR2529937, MR2680495, MR2874638}. An analogous Fourier--Jacobi model is defined when $\dim W^\perp$ is even, for symplectic--metaplectic and unitary groups, by replacing \(\nu\) with a Weil representation \cite[Sect. 12]{MR3202556}.

\vskip 5pt

When the representations belong to \emph{tempered} (or more generally, \emph{generic}) local L-packets, Gan, Gross, and Prasad predicted a remarkably elegant answer in \cite{MR3202556}. Namely, for any pair of tempered L-parameters $\left(\phi_V, \phi_W\right)$, as $\left(\pi,\sigma\right)$ runs over all irreducible representations of relevant pure inner forms in the associated Vogan L-packet $\LP_{\phi_V} \times \LP_{\phi_W}$, the sum of multiplicities is always equal to $1$, and the unique distinguished pair is indexed by a character of the component group defined using symplectic local root numbers (recalled in Section \ref{sub:local_gan_gross_prasad_conjecture_tempered_case}). This tempered local conjecture is now a theorem in both the Bessel and Fourier--Jacobi settings over any local field. Over non-Archimedean local fields, this was established by the landmark works of Waldspurger \cite{MR2684300, MR3155344, MR3155345, MR3202558, MR3202559} and M{\oe}glin--Waldspurger \cite{MR3155346} in the orthogonal case, Beuzart-Plessis \cite{MR3251763, MR3676153} in the Hermitian case, Gan--Ichino \cite{MR3573972} in the skew-Hermitian case, and Atobe \cite{MR3788848} in the symplectic-metaplectic case (see also \cite{MR4146145, MR4568696, MR5115665, luo2020local, chenluo2022local} for Archimedean cases).

\subsection{The non-tempered conjecture and main result}\label{sub:introNTGGPandMainResult}

With the tempered case settled, it is natural to ask what happens for other representations, especially those arising as local components of automorphic representations in the discrete automorphic spectrum. Because the generalized Ramanujan conjecture fails for classical groups, such local components need not be tempered. Instead, Arthur's endoscopic classification organizes them via local \emph{Arthur parameters} (A-parameters), which are admissible homomorphisms
\[
    \psi : WD_F \times \SL_2(\C) \longrightarrow {}^L G
\]
with bounded image on the Weil group $W_F$ and algebraic restriction to the two (``Deligne'' and ``Arthur'') $\SL_2\left(\C\right)$ factors. Each local A-parameter $\psi$ gives rise to two closely related finite sets of representations:
\begin{itemize}
    \item Pulling back the A-parameter $\psi$ along the homomorphism 
    \[
        WD_F \lra WD_F\times \SL_2\left(\C\right), \quad w\mapsto \left(w,\left(\begin{array}{cc}
                |w|^{\frac{1}{2}} & \\
                 & |w|^{-\frac{1}{2}}
            \end{array}\right)\right)
    \]
    produces an L-parameter $\phi_\psi$ of \emph{Arthur type}, and hence a Vogan L-packet $\LP_{\phi_\psi}$. When the Arthur $\SL_2\left(\C\right)$ acts non-trivially, $\phi_\psi$ is non-tempered.

    \vskip 5pt

    \item (Twisted) endoscopic character relations attach to the A-parameter $\psi$ a local A-packet $\AP_\psi$ consisting of unitary representations and containing $\LP_{\phi_\psi}$ as a subset. These local A-packets $\AP_\psi$ are often strictly larger than $\LP_{\phi_\psi},$ and need not be disjoint for distinct A-parameters.
\end{itemize}
\vskip 5pt
In the non-tempered setting, the sum of multiplicities over the relevant subset of a Vogan L-packet $\left(\LP_{\phi_{\psi_V}} \times \LP_{\phi_{\psi_W}}\right)^\rel$ is no longer always equal to $1$: it vanishes for most pairs of A-parameters.
In \cite{MR4190046}, Gan, Gross, and Prasad proposed a precise combinatorial condition on a pair of A-parameters that governs the total multiplicity over a Vogan L-packet $\left(\LP_{\phi_{\psi_V}} \times \LP_{\phi_{\psi_W}}\right)^\rel$, calling it \emph{relevance} (see \cite[Sect. 3]{MR4190046}, and also Definition \ref{D:relevance-definition}). Roughly speaking, two A-parameters form a relevant pair if their non-tempered Jordan blocks can be matched in pairs whose (Arthur) $\SL_2\left(\C\right)$-types differ by $1$.
This notion is a refinement of the nilpotent-orbit constraints in the unitary $L^2$-setting studied by Clozel and Venkatesh \cite{Clo04, Ven05}.

\vskip 5pt

As a first test ground, one can consider the pair $\left(\GL_{n+1}\left(F\right), \GL_n\left(F\right)\right)$, which can be obtained formally by replacing $E$ with the split quadratic algebra $F \times F$. This is also known as the ``\emph{Rankin--Selberg model}''. In this setting, both A-packets and L-packets are singletons, and Gan, Gross, and Prasad conjectured that for a pair of irreducible Arthur type representations $\left(\pi,\sigma\right)$ of $\GL_{n+1}\left(F\right)\times \GL_n\left(F\right)$, the Hom space 
\[
    \Hom_{\GL_n\left(F\right)}\left(\pi,\sigma\right)
\]
is non-zero if and only if the corresponding A-parameters form a relevant pair \cite[Conj. 5.1]{MR4190046}. Over non-Archimedean local fields, special cases of this conjecture for $\GL_n$ had already been shown in their original paper. Later, Gurevich \cite{Gur22} proved the ``period implies relevance'' direction of this conjecture using the machinery from quantum groups. Subsequently, this conjecture is completely resolved by Chan \cite{MR4373242}. A key ingredient in Chan's proof is the fully-faithfulness of the parabolic induction functor $u_\rho\left(a,b\right)\times-$ under suitable cuspidal support conditions. See also \cite{chenchenNTGGP4ArchiGL} for the Archimedean case.

\vskip 5pt

However, much less is known for classical groups. Several difficulties arise in this setting, namely, L-packets and A-packets are no longer singletons, different A-packets can have complicated intersections, the theory of derivatives is considerably more intricate, and there is no direct passage to quantum affine algebras as in \cite{Gur22, MR4373242}. Focusing on the Bessel setting, the non-tempered local conjecture of Gan, Gross, and Prasad is as follows (see Section \ref{sec:bessel_periods_conjectures_and_some_inequalities} for precise conventions regarding full versus special orthogonal groups and relevant pure inner forms).

\begin{Conj}[{\cite[Conj. 6.1 \& 7.1]{MR4190046}}]\label{Conj:Non-temperedGGP-intro}
Let $W \subset V$ be non-degenerate quadratic or Hermitian spaces over $E$ as above, and let $\psi_V$, $\psi_W$ be local A-parameters for the corresponding orthogonal or unitary groups. Denote by $\left(\LP_{\phi_{\psi_V}} \times \LP_{\phi_{\psi_W}}\right)^\rel$ and $\left(\AP_{\phi_{\psi_V}} \times \AP_{\phi_{\psi_W}}\right)^\rel$ the subsets of the Vogan L-packet and A-packet consisting of representations of relevant pure inner forms respectively.
\begin{enumerate}
       \item (``Relevance implies period'') If $\left(\psi_V,\psi_W^\vee\right)$ is relevant, then the strong multiplicity one property
\[
    \sum_{\left(\pi,\sigma\right)\in \left(\LP_{\phi_{\psi_V}}\times \LP_{\phi_{\psi_W}}\right)^\rel} m\left(\pi,\sigma\right)=1
\]
        still holds, and the unique member $\left(\pi,\sigma\right)$ such that $m\left(\pi,\sigma\right)\neq 0$ is indexed by the distinguished character given by the same recipe as in \cite[Sect. 6]{MR3202556}.

        \vskip 5pt

        \item (``Period implies relevance'') If there exists $\left(\pi,\sigma\right)\in\left(\AP_{\psi_V}\times \AP_{\psi_W}\right)^\rel$ such that $m\left(\pi,\sigma\right)\neq 0$, then there exist local A-parameters $\psi_{V,\pi}$ and $\psi_{W,\sigma}$, such that $\pi\in\AP_{\psi_{V,\pi}}$, $\sigma\in\AP_{\psi_{W,\sigma}}$, and $\left(\psi_{V,\pi},\psi_{W,\sigma}^\vee\right)$ is relevant. Moreover, we may take $\psi_{V,\pi}=\psi_V$ and $\psi_{W,\sigma}=\psi_W$ if either of the following holds:
        \begin{itemize}
            \item[a.] $\left(\pi,\sigma\right)\in\left(\LP_{\phi_{\psi_V}}\times \LP_{\phi_{\psi_W}}\right)^\rel$; or

            \vskip 5pt

            \item[b.] $\psi_V$ and $\psi_W$ are co-tempered, i.e. the Deligne $\SL_2\left(\C\right)$ acts trivially.
        \end{itemize}
\end{enumerate}
\end{Conj}

\vskip 5pt

When $\psi_V$ and $\psi_W$ are tempered, $\left(\psi_V, \psi_W^\vee\right)$ is automatically relevant and Conjecture \ref{Conj:Non-temperedGGP-intro} reduces to the tempered case. Note also that because $\LP_{\phi_\psi}$ is contained in $\AP_\psi$, (2.b) immediately implies (2.a) when $\psi_V$ and $\psi_W$ are co-tempered. By contrast, for general A-packets, Gan, Gross, and Prasad already showed in \cite[Sect. 7]{MR4190046} that intersections of A-packets force the more intricate formulation. The main result of this paper is as follows.

\begin{Thm}[Theorem \ref{T:mainresult}]\label{Thm:main-intro}
In the non-Archimedean Bessel setting (orthogonal and unitary groups), the ``relevance implies period'' direction holds in full generality, and the ``period implies relevance'' direction holds for co-tempered A-parameters; that is, parts (1) and (2.b) of Conjecture \ref{Conj:Non-temperedGGP-intro} hold.
\end{Thm}

\vskip 5pt

\begin{Rmk}
Here, we should mention that our proof of this theorem relies on several widely expected results on A-packets, especially for non-quasi-split pure inner forms. More precisely, we assume that Conjecture \ref{Conj:APacketsNProperties}, Lemma \ref{L:AtobeJacLLCLem.4.1}, Theorem \ref{T:AtobeJacLLCMainThm4.2} and Proposition \ref{P:unitaryInd} hold for all orthogonal and unitary groups under consideration. While these results have already been established in many cases, for certain pure inner forms they are not yet available in the literature.   
\end{Rmk}

\vskip 5pt

While part (2.b) of the conjecture holds, one can see already in low rank cases that part (2.a) of the conjecture fails. Specifically, consider the embedding of the compact two-dimensional unitary group $\U(V_2^-)$ into the split three-dimensional unitary group $\U(V_3^+)$. Because $\U(V_2^-)$ is compact, its trivial representation $\mathbbm{1}_2$ is tempered, whereas the trivial representation $\mathbbm{1}_3$ of $\U(V_3^+)$ lies in the Arthur type L-packet associated to the A-parameter $\mathbbm{1} \boxtimes S_3$. Trivially
\[
    \Hom_{\U(V_2^-)}(\mathbbm{1}_3, \mathbbm{1}_2) \neq 0,
\]
yet the A-parameter of $\mathbbm{1}_3$ and the tempered A-parameter of $\mathbbm{1}_2$ are not relevant. Fortunately, this example does not violate the statement for general case in part (2). Indeed, it's easy to see that $\mathbbm{1}_2$ also belongs to the A-packet attached to the A-parameter $\mathbbm{1} \boxtimes S_2$, and the pair $\left(\mathbbm{1} \boxtimes S_3,\mathbbm{1} \boxtimes S_2\right)$ is relevant. Nevertheless, we still believe that part (2.a) could be true after some modifications. Indeed, besides the co-tempered case we also establish part (2.a) in some special cases (see Theorem \ref{T:ReductionToGoodParity}(2) and Corollary \ref{C:PeriodImpliesRelevanceEvenParameter} for more details), and in Section \ref{sec:counterexample} we propose several ways to repair part (2.a), by restricting to split forms, or relaxing the relevance condition to Qadri's notion of strong Ext-relevance \cite{MR4846726, cq26}.

\subsection{Motivation and broader connections}

Beyond its intrinsic representation-theoretic interest, the local non-tempered GGP conjecture is closely tied to global automorphic periods, algebraic cycles on Shimura varieties, and relative Langlands duality.

\vskip 5pt

$\bullet$ \textit{Global periods and automorphic descent.}
Over a number field, the global non-tempered GGP conjecture predicts that a (suitably regularized) Bessel or Fourier--Jacobi period of non-tempered cuspidal or discrete automorphic representation is non-vanishing if and only if the global A-parameters are relevant, local distinction holds at every place, and a certain L-value is non-zero at $s=0$ \cite[Conj. 9.1]{MR4190046}. For discrete global A-parameters $\Psi, \Psi'$ of orthogonal groups, with $\Psi$ symplectic and $\Psi'$ orthogonal, this L-value takes the form
\[
    L\left(\Psi,\Psi',s\right) = \frac{L\left(\phi_\Psi \otimes \phi_{\Psi'},s+\frac{1}{2}\right)}{L\left({\rm Sym}^2\phi_\Psi\oplus\wedge^2\phi_{\Psi'},s+1\right)},
\]
which is holomorphic at $s=0$ whenever $\left(\Psi, \Psi'\right)$ is relevant \cite[Thm. 9.7]{MR4190046}. This global non-tempered conjecture remains largely open (see Haan--Kwon \cite{HK26} for recent progress in certain low-rank cases). 

\vskip 5pt

One of the main motivations for these non-tempered conjectures is the theory of \emph{automorphic descent}, where non-tempered Bessel and Fourier--Jacobi periods inevitably arise. Indeed, Ginzburg--Rallis--Soudry \cite{GRS11} construct generic cuspidal representations of classical groups from Fourier coefficients of residual Eisenstein series, and the representations involved are typically non-tempered even when the target parameter is tempered. Jiang--Zhang's twisted automorphic descent \cite{JZ20} also uses auxiliary (non-tempered) Bessel periods to detect non-generic members of global L-packets. As explained in \cite[Sect. 12]{MR4190046}, Theorem \ref{Thm:main-intro} provides the local branching input needed to predict the parameters and representations that can be reached by such descent constructions.

\vskip 5pt

$\bullet$ \textit{Cycles on Shimura varieties.}
Pairing algebraic cycles arising from subgroup embeddings with automorphic cohomology classes provides a geometric incarnation of period integrals. Because non-tempered A-packets contribute cohomology outside the middle degree, understanding how cycle classes pair against such cohomology requires non-tempered branching laws. For instance, Rapoport--Smithling--Zhang observed that describing the Hecke module generated by ordinary unitary diagonal cycles calls for non-tempered GGP \cite[Rmk. 6.16]{RSZ20}, Sweeting used explicit non-tempered $\SO_4 \subset \SO_5$ periods to determine when sub-Shimura correspondences act non-trivially on cohomology attached to Yoshida lifts \cite[Thm. A \& C]{Swe22}, and Xue recently applied non-tempered GGP periods to construct algebraic Tate classes in arbitrary codimension on $\U(n,1)$ Shimura varieties \cite{Xue26}.

\vskip 5pt

$\bullet$ \textit{Symplectic geometry and relative Langlands duality.}
Locally, Gan, Gross, and Prasad (following Zhiwei Yun) showed that relevance admits a geometric interpretation in terms of correlators and moment maps for reductive dual pairs \cite[Sect. 4]{MR4190046}. In the relative Langlands framework of Ben-Zvi--Sakellaridis--Venkatesh \cite{bzsv}, a local A-parameter of $G$ determines a subvariety of the dual Hamiltonian $G^\vee$-variety as the moment-map preimage of a Slodowy slice, and relevance translates into the existence of fixed points on this subvariety under the extended action induced by this A-parameter \cite[Sect. 14.3]{bzsv}. As highlighted by Gan and Wang \cite{gw25}, one may then view the non-tempered GGP conjecture as the dual problem to Adams' conjecture on theta correspondence.

\subsection{Strategy of proof and organization}

Let us briefly outline the proof of Theorem \ref{Thm:main-intro}. For both parts (1) and (2) of Conjecture \ref{Conj:Non-temperedGGP-intro}, our first step is to reduce to the case where $\psi_V$ and $\psi_W$ are of \emph{good parity} (i.e. every irreducible summand is (conjugate) self-dual of the same type as the ambient parameter). Since general A-packets can be constructed through unitary parabolic induction of Speh representations and good-parity A-packets, we combine Mackey theory with Kret--Lapid's Jacquet module formula for ladder representations \cite{MR2996769} to separate good parity part and bad parity part. In fact, we work in a slightly more general setting, and as a bonus, we obtain a partial result on ``period implies relevance'' direction (Corollary \ref{C:PeriodImpliesRelevanceEvenParameter}). Once reduced to good parity case, our proofs of parts (1) and (2.b) proceed along different routes:
\begin{itemize}
    \item \emph{Proof of Part (1):} By using M{\oe}glin--Waldspurger's result on the multiplicities for standard modules (Proposition \ref{mwprop13}) together with the tempered local GGP conjecture, we can deduce the multiplicity at most one property for the Vogan L-packet and the epsilon dichotomy. To prove the existence of distinguished representation under the assumption of relevance, we induct on a certain non-negative integer which measures of non-temperedness of $\left(\psi_V, \psi_W\right)$, using inequalities provided by Mackey theory to reduce the non-temperedness of the parameters. An essential bottleneck arises when one parameter is tempered and the other has Arthur $\SL_2(\C)$-blocks of size at most $2$. 
    We overcome this by explicitly decomposing the relevant unitary and non-unitary parabolic inductions (building on \cite{MR4055180, MR4412277, clw26}) and comparing distinguished characters to rule out unwanted constituents.

    \vskip 5pt

    \item \emph{Proof of Part (2.b):} Co-tempered A-packets are obtained from tempered L-packets via the Aubert--Zelevinsky involution. Therefore, just like tempered L-packets, two co-tempered A-packets are either identical or disjoint, and their (semi-simplified) derivatives can be tracked explicitly by using Atobe's results on Jacquet modules \cite{MR4055180}. We prove by induction on $\dim V + \dim W$. The case that both $\pi$ and $\sigma$ are supercuspidal has been dealt in the original paper of Gan, Gross and Prasad \cite[Thm. 7.7]{MR4190046}. When $\pi$ or $\sigma$ is not supercuspidal, we extract a derivative of minimal exponent from $\left(\pi,\sigma\right)$ and construct a pair of representations on smaller groups that is ``more tempered'' and still admits a non-zero Bessel period. Applying the induction hypothesis yields almost all the information that we need to show the relevance of A-parameters, and the remaining information is obtained through some delicate case-by-case analysis on other derivatives of $\pi$ and $\sigma$.
\end{itemize}

\vskip 5pt

The paper is organized as follows. In Section \ref{sec:classical_groups_and_their_representations}, we review the notion of classical groups, the Langlands classification, Tadi\'c's formula, and several notions of semi-simplified and homological derivatives. We also establish a vanishing criterion for derivatives and a specialized Leibniz rule. In Section \ref{sec:l_packets_a_packets_and_representations_of_arthur_type}, we recall the local Langlands correspondence, local A-packets, and the current status of Arthur's endoscopic classification. In Section \ref{sec:certain_jacquet_modules_and_parabolic_inductions}, we collect some results of Atobe on derivatives, and study the decompositions of certain parabolic inductions needed in later proofs. In Section \ref{sec:bessel_periods_conjectures_and_some_inequalities}, we introduce Bessel periods and (homological) multiplicities, prove some inequalities using orbit analysis, and state the tempered and non-tempered GGP conjectures precisely. Section \ref{sec:reduction_to_the_good_parity_case} is devoted to the proof of Theorem \ref{T:ReductionToGoodParity}, which allows us to reduce the conjecture to the good parity case. In Sections \ref{sec:relevance_implies_period} and \ref{sec:period_implies_relevance}, we prove parts (1) and (2.b) of Conjecture \ref{Conj:Non-temperedGGP-intro}, respectively, in the good parity case. Finally, Section \ref{sec:counterexample} presents the counterexample to part (2.a) and formulates modified conjectures. In Apppendix \ref{App:ProofTechnicalCountingMulti}, we prove a technical lemma to supply the proof of Theorem \ref{T:ReductionToGoodParity}(2).

\vskip 5pt

\subsection*{Acknowledgements}
The authors would like to thank Wee Teck Gan for suggesting this project and for many inspiring conversations. We also thank Hiraku Atobe, Kei Yuen Chan, Maxim Gurevich and Chi-Heng Lo for helpful discussions and comments. Part of this work was carried out during the conference ``Wavefront Sets of Representations and Endoscopic Classification'' held in Hangzhou in October 2026; we thank the organizers and the Institute for Advanced Study in Mathematics (IASM) for their warm hospitality. R. Chen is partially supported by the National Natural Science Foundation of China (Project No. 12526208). J-Q. Chong is supported by the NUS Research Scholarship. C-J. Wang is supported by the Research Grants Council of the Hong Kong Special Administrative Region, China (Project No. 17305223, 17308324).

\vskip 5pt

\subsection*{AI declaration}
The first draft of this manuscript was prepared entirely by the three authors without the use of AI tools. During subsequent revisions, we used AI tools to help with proofreading. In the course of this review, we identified a gap in the proof of Theorem \ref{T:ReductionToGoodParity}(2), which we have corrected in the present version with the assistance of ChatGPT. The arguments in Appendix \ref{App:ProofTechnicalCountingMulti} are largely due to ChatGPT. We also used Gemini to help polish the language of the introduction. The authors take full responsibility for the content of the manuscript.

\vskip 10pt

\section{Classical groups and their representations}\label{sec:classical_groups_and_their_representations}

In this section, we review the basic notions of classical groups, the Langlands classification and the notions of derivatives.

\subsection{Formed spaces and isometry groups}\label{subsec:FormedSpaceNIsomGroup}

Let us start by giving the definition of classical groups. Fix a vector space \(V\) over a non-Archimedean local field \(E\) of characteristic $0$ with dimension \(n\). The \emph{general linear group} \(\GL(V)\) is the set of all invertible linear transformations on \(V\). To simplify notation, we shall often omit the vector space \(V\) and write \(\GL_n\) instead.

\vskip 5pt

Let \(\iota\) be an involution of \(E\) with fixed field \(F\), so that \(E=F\) if \(\iota\) is the identity map, or \(E\) is a quadratic extension of \(F\) otherwise. Endow \(V\) with a non-degenerate \(\iota\)-sesquilinear form \(\langle\cdot,\cdot\rangle\), such that for all \(u,v,w\in V\) and \(\alpha,\beta\in E,\) one has
\begin{itemize}
    \item \(\langle\alpha v+\beta w,u\rangle=\alpha\langle v,u\rangle+\beta\langle w,u\rangle\);

    \vskip 5pt

    \item \(\langle u,v\rangle=\varepsilon\cdot\langle v,u\rangle^{\iota}\) for a fixed choice of \(\varepsilon\in\{\pm 1\}\).
\end{itemize}
\vskip 5pt
Let \(G(V)\) be the subgroup of \(\GL(V)\) consisting of elements \(g\) such that \(\langle gv,gw\rangle=\langle v,w\rangle\) for all \(v,w\in V\). Depending on the choices of \(\iota\) and \(\varepsilon\), we have the following possibilities.

\vskip 5pt

\begin{center}
\begin{tabular}{|c|c|c|c|}
\hline
         & $\iota={\rm id}$, $\varepsilon=1$  & $\iota={\rm id}$, $\varepsilon=-1$  & $\iota\neq{\rm id}$, $\varepsilon=\pm1$ \\ \hline
\(G(V)\) & orthogonal group \(\rO(V)\) & symplectic group \(\Sp(V)\) & unitary group \(\U(V)\)                     \\ \hline
\end{tabular}
\end{center}
\vskip 5pt
As in the general linear group case, when there is no confusion we will opt to write \(G_n\) instead of \(G(V)\), where \(n=\dim V\). The following table lists out the complex Langlands group of the classical groups. 
\begin{center}
\begin{tabular}{|c|c|c|c|c|}
\hline
\(G(V)\) & \(\rO(V)\), $n$ odd & \(\rO(V)\), $n$ even & \(\Sp(V)\) & \(\U(V)\) \\\hline
\(^LG(V)\) & \(\Sp_{n-1}(\C)\) & $\begin{array}{cc}
\SO_n(\C) & \text{if } \disc(V) \in F^{\times 2}\\
\rO_n(\C) & \text{if } \disc(V) \notin F^{\times 2} 
\end{array}$ & \(\SO_{n+1}(\C)\) & \(\GL_n(\C)\rtimes\mathrm{Gal}(E/F)\)\\ \hline
\end{tabular}
\end{center}
\vskip 5pt
Here \(\mathrm{disc}(V)\) is the discriminant of \(V\). Since we focus on the Bessel case, we are only interested in the restriction problems on orthogonal groups and unitary groups. 
We should mention that when $\iota = {\rm id}$, originally the GGP problems were formulated for \emph{special} orthogonal groups rather than full orthogonal groups. For the discrepancy between special orthogonal groups and full orthogonal groups, we refer the readers to \cite[Sect. 2.5]{MR3788848} and \cite[Sect. 5.2]{MR3708200}. For convenience we set
\[
    G^0(V) \coloneqq \begin{cases}
     \SO(V) \quad & \text{ if $\iota={\rm id}$, $\varepsilon=1$ and $\dim V$ is odd};\\
     G(V) \quad & \text{ otherwise.}   
    \end{cases}
\]

\vskip 5pt

\subsection{Langlands classification} 

In this subsection, we recall the Langlands classification for general linear groups and classical groups.

\vskip 5pt

Let us first define some representations of general linear groups. Let \(P\) be a standard parabolic subgroup of $\GL_n$ such that its Levi subgroup is isomorphic to \(\GL_{n_1}\times \GL_{n_2}\times\cdots\times\GL_{n_r}\). Here ``standard'' means that $P$ is block upper triangular. For \(1\leq i\leq r\), let \(\tau_i\) be an irreducible representation of \(\GL_{n_i}.\) Then, we write
\[
    \tau_1\times \tau_2\times\cdots\times\tau_r\coloneqq\Ind^{\GL_n}_{P}\left(\tau_1\boxtimes\tau_2\boxtimes\cdots\boxtimes\tau_r\right)
\]
as the normalized parabolic induction. Define \(\calC(\GL_n)\) as the collection of irreducible supercuspidal representations of \(\GL_n\), and \(\calC_{\unit}(\GL_n)\) as the subset of $\calC(\GL_n)$ consisting of irreducible unitary supercuspidal representations. We set $\calC(\GL)\coloneqq \bigcup_{n\geq 0}\calC(\GL_n)$ and $\calC_{\unit}(\GL)\coloneqq \bigcup_{n\geq 0}\calC_{\unit}(\GL_n)$. A \emph{segment}, denoted by \([x,y]_{\rho},\) is a set of supercuspidal representations of the form
\[
    \{\rho|\cdot|^x,\rho|\cdot|^{x-1},\cdots,\rho|\cdot|^y\}
\]
where \(\rho\in\calC_{\unit}(\GL_n)\), and \(x,y\in\R\) such that \(x\geq y\) and \(x-y\) is an integer. To a segment \([x,y]_{\rho}\) we associate a \emph{Steinberg representation} \(\Delta_{\rho}[x,y]\), which is the unique irreducible subrepresentation of 
\[
    \rho|\cdot|^x\times\rho|\cdot|^{x-1}\times\cdots\times\rho|\cdot|^y,
\] 
or the unique irreducible quotient of \(\rho|\cdot|^y\times\rho|\cdot|^{y+1}\times\cdots\times\rho|\cdot|^x.\) On the other hand, the unique irreducible quotient of \(\rho|\cdot|^x\times\rho|\cdot|^{x-1}\times\cdots\times\rho|\cdot|^y\) (or the unique irreducible subrepresentation of \(\rho|\cdot|^y\times\rho|\cdot|^{y+1}\times\cdots\times\rho|\cdot|^x\)) is denoted as \(Z_{\rho}[x,y].\) By convention, we take \(\Delta_{\rho}[x,x+1]\) and \(Z_{\rho}[x,x+1]\) as the trivial representation of the trivial group.

\vskip 5pt

We say two segments \([x_1,y_1]_{\rho_1}\) and \([x_2,y_2]_{\rho_2}\) are \emph{linked} if one has \(\rho_1\cong\rho_2,\) and \([x_1,y_1]_{\rho_1}\not\subseteq [x_2,y_2]_{\rho_2},\) \([x_2,y_2]_{\rho_2}\not\subseteq [x_1,y_1]_{\rho_1}, \) but \([x_1,y_1]_{\rho_1}\cup[x_2,y_2]_{\rho_2}\) still forms a segment. This condition gives a criterion on the commutativity of the parabolic induction of two Steinberg representations, as described below.

\begin{Lem}[{\cite[Thm. 4.2]{ZeleII}}]\label{linkedlemma}
    If \([x_1,y_1]_{\rho_1}\) and \([x_2,y_2]_{\rho_2}\) are not linked, then \(\Delta_{\rho_1}[x_1,y_1]\times\Delta_{\rho_2}[x_2,y_2]\) is irreducible, and 
    \[
        \Delta_{\rho_1}[x_1,y_1]\times\Delta_{\rho_2}[x_2,y_2]\cong \Delta_{\rho_2}[x_2,y_2]\times\Delta_{\rho_1}[x_1,y_1].
    \]
\end{Lem}

\vskip 5pt

We now present the Langlands classification for \(\GL_n.\)

\begin{Thm}[Langlands Classification for \(\GL_n\)]\label{T:LanglandsClass4GL}
For every \(\tau\in\Irr(\GL_n(E)),\) there exist
\begin{itemize}
    \item \(\rho_1,\rho_2,\cdots,\rho_r\in\calC_{\unit}\) (which need not be distinct); and

    \vskip 5pt

    \item \(x_1,x_2,\cdots,x_r,y_1,y_2,\cdots,y_r\in\R\) satisfying \(x_1+y_1\leq\cdots\leq x_r+y_r\),
\end{itemize}
such that \(\tau\) is the unique irreducible subrepresentation of
\[
    \Delta_{\rho_1}[x_1,y_1]\times\cdots\times\Delta_{\rho_r}[x_r,y_r].
\]
In this case, we write 
\[
    \tau=\LS(\Delta_{\rho_1}[x_1,y_1],\cdots,\Delta_{\rho_r}[x_r,y_r]).
\]
\end{Thm}

\vskip 5pt

Let \(\left(x_{i,j}\right)_{1\leq i\leq s,1\leq j\leq t}\) be a matrix of real numbers such that \(x_{i,j}=x_{1,1}-i+j\). We define a \emph{(shifted) Speh representation} by
\[
\begin{pmatrix}
    x_{1,1} & \cdots & x_{1,t}\\
    \vdots & \ddots & \vdots\\
    x_{s,1} & \cdots & x_{s,t}
\end{pmatrix}_{\rho} \coloneqq \LS(\Delta_{\rho}[x_{1,1},x_{s,1}],\cdots,\Delta_{\rho}[x_{1,t},x_{s,t}]).
\]
We also write
\[
u_{\rho}(a,b):=\begin{pmatrix}
    \frac{a-1}{2}+\frac{1-b}{2} & \cdots & \frac{a-1}{2}+\frac{b-1}{2}\\
    \vdots & \ddots & \vdots\\
    \frac{1-a}{2}+\frac{1-b}{2} & \cdots & \frac{1-a}{2}+\frac{b-1}{2}
\end{pmatrix}_{\rho}
\]
for positive integers \(a\) and \(b\). Note that \(u_{\rho}(a,1)=\Delta_{\rho}\left[\frac{a-1}{2},\frac{1-a}{2}\right]\) and \(u_{\rho}(1,b)=Z_{\rho}\left[\frac{b-1}{2},\frac{1-b}{2}\right]\). As above, we have the commutativity of the parabolic induction of two Speh representations.

\begin{Lem}[{\cite[Thm. 1.1]{MR3224483}}]
Let \(\left(x_{i,j}\right)_{\rho}\) and \(\left(y_{i',j'}\right)_{\rho'}\) be two Speh representations, with \(1\leq i\leq s,1\leq j\leq t\) and \(1\leq i'\leq s',1\leq j'\leq t'\). If the two segments \([x_{1,t},x_{s,1}]_{\rho}\) and \([y_{1,t'},y_{s',1}]_{\rho'}\) are not linked, then the representation \((x_{i,j})_{\rho}\times (y_{i',j'})_{\rho'}\) is irreducible. Moreover, we have
\[
    (x_{i,j})_{\rho}\times (y_{i',j'})_{\rho'}\cong (y_{i',j'})_{\rho'}\times (x_{i,j})_{\rho}.
\]
\end{Lem}

\vskip 5pt

Every irreducible representation \(\tau\) of \(\GL_n\) admits a central character \(\omega_{\tau}\) on \(E^{\times}\). There is also a real number \(e_{\tau}\) such that \(|\omega_{\tau}(\lambda)|=|\lambda|^{n\cdot e_{\tau}}\) for every \(\lambda\in E^{\times}\). The number \(e_{\tau}\) is known as the \emph{exponent} of \(\tau\). If \(\tau\) is unitary, e.g. of the form \(\Delta_{\rho}[a,-a]\) for some half-integer \(a\) or \(u_{\rho}(a,b)\) for some positive integers \(a\) and \(b\), then for any $c\in\R$ the exponent of \(\tau|\cdot|^c\) is \(c\).

\vskip 5pt

We proceed to the Langlands classification for \(G_n= G(V)\). We fix a minimal parabolic subgroup $P_0$ of $G_n$, and call parabolic subgroups containing $P_0$ \emph{standard} parabolic subgroups. Let \(P\) be a parabolic subgroup of $G_n$ with the Levi subgroup
\[
    M \cong \GL_{n_1}\times\GL_{n_2}\times\cdots\times\GL_{n_r}\times G_{n_0}.
\]
Here $G_{n_0}=G(V_0)$ for some space $V_0$ in the Witt tower containing $V$ and of appropriate dimension. Note that \(G_n\) is quasi-split if and only if \(G_{n_0}\) is quasi-split. If \(\tau_i\) is an irreducible representation of \(\GL_{n_i}\) for \(1\leq i\leq r\) and \(\pi_0\) is an irreducible representation of \(G_{n_0}\), we set
\[
    \tau_1\times\tau_2\times\cdots\times\tau_r\rtimes\pi_0:=\Ind^{G_n}_{P}\left(\tau_1\boxtimes\tau_2\boxtimes\cdots\boxtimes\tau_r\boxtimes\pi_0\right)
\]
as the normalized parabolic induction. With this notation, we can state the Langlands classification for classical groups.

\begin{Thm}[Langlands classification for \(G_n\)]\label{lcgn}
    For every \(\pi\in\Irr(G_n),\) there exist
    \begin{itemize}
        \item \(\rho_1,\rho_2,\cdots,\rho_r\in\calC_{\unit}\) (which need not be distinct);

        \vskip 5pt

        \item \(x_1,x_2,\cdots,x_r,y_1,y_2,\cdots,y_r\in\R\) satisfying \(x_1+y_1\leq x_2+y_2\leq\cdots \leq x_r+y_r<0\) ($r\in\Z_{\geq0}$ can be zero); and

        \vskip 5pt

        \item an irreducible tempered representation \(\pi_0\) of \(G_{n_0}\),
    \end{itemize}
    such that \(\pi\) appears as the unique irreducible subrepresentation of 
    \[
        \Delta_{\rho_1}[x_1,y_1]\times\cdots\times\Delta_{\rho_r}[x_r,y_r]\rtimes\pi_0.
    \] 
    In this case, we write 
    \[
        \pi=\LS(\Delta_{\rho_1}[x_1,y_1],\cdots,\Delta_{\rho_r}[x_r,y_r];\pi_0).
    \]
\end{Thm}

\vskip 5pt

It is often very useful to realize \(\pi\) as a quotient of some representation. To this end, we can use the \emph{\(\MVW\)-functor} \(\pi\mapsto\pi^{\MVW}\) defined on the category of smooth representations of \(G_n\) \cite{MR1041060}. This is an \emph{involutive} covariant functor satisfying the following properties:
\begin{itemize}
    \item If \(\pi\in\Irr(G_n)\) is irreducible, then \(\pi^{\MVW}\cong\pi^{\vee}\);

    \vskip 5pt

    \item For any representation \(\tau\) of \(\GL_d\) and \(\sigma\) of \(G_{n_0}\) with \(n=n_0+2d\), one has
    \[
        (\tau\rtimes\sigma)^{\MVW}\cong \prescript{\iota}{}\tau\rtimes\sigma^{\MVW}.
    \]
\end{itemize}
\vskip 5pt
Here, the representation \(\prescript{\iota}{}\tau\) is the \(\iota\)-conjugate of \(\tau\). It is defined as follows. Let \(X\) be a vector space over \(E\) with dimension \(d\), spanned by the vectors \(\{x_1,x_2,\cdots,x_d\}.\) Set
\[
    j:\GL_d(E)\to \GL(X),\quad g\mapsto \left((x_1,\cdots,x_d)\mapsto (x_1,\cdots,x_d)g\right),
\]
i.e. \(j\) is the map identifying invertible matrices of size \(d\) with the corresponding linear operators on \(X\). Then, we define
\[
    \prescript{\iota}{}\tau(h):=\tau\left(j\circ \iota\circ j^{-1}(h)\right)
\]
for all \(h\in\GL(X)\). It can be shown that \(\prescript{\iota}{}\tau\) is independent of the choice of a basis of \(X\). To simplify the notation, for any representation $\tau$ of a general linear group, we define
\begin{equation}\label{E:TildeRepDef}
    \widetilde{\tau} =\begin{cases}
\tau^\vee \quad&\text{ if }E=F;\\
\prescript{\iota}{}\tau^\vee \quad&\text{ otherwise.}
\end{cases}
\end{equation}
Now, if $\pi= \LS\left(\Delta_{\rho_1}[x_1,y_1],\cdots,\Delta_{\rho_r}[x_r,y_r];\pi_0\right)$ is an irreducible representation of $G_n$, by dualizing and applying the MVW-functor, we obtain a surjection
\[
    \Delta_{\widetilde{\rho}_1}[-y_1,-x_1]\times\cdots\times\Delta_{\widetilde{\rho}_r}[-y_r,-x_r]\rtimes\pi_0\twoheadrightarrow\pi.
\]
We call the representation on the left-hand side the \emph{standard module} of \(\pi\), and denote it by \(\calM(\pi)\). We also define the \emph{co-standard module} of \(\pi\) as
\[
    \overline{\calM}(\pi):=\Delta_{\rho_1}[x_1,y_1]\times\cdots\times\Delta_{\rho_r}[x_r,y_r]\rtimes\pi_0.
\]
Note that by definition, there is an embedding \(\pi\hookrightarrow\overline{\calM}(\pi)\).

\vskip 5pt

\subsection{Tadi\'c's Formula and derivatives}\label{subsec:TadicFormulaDerivative}

In this subsection we shall work with $\GL_n$ and $G^0(V)$, namely, symplectic groups, special odd orthogonal groups, even orthogonal groups and unitary groups. For a group \(G\), write \(\calR(G)\) for the Grothendieck group of the category of smooth complex representations of \(G\) of finite length. We also set
\[
    \calR:=\bigoplus_{n\geq 0}\calR\left(\GL_n\right).
\]
Given a formed space $V$, for each $n\in \Z_{\geq 0}$ such that $n\equiv \dim V \pmod2$, we denote by $V_n$ the space of dimension $n$ in the Witt tower containing $V$, and 
\[
    \calR_V:=\bigoplus_{n\geq 0}\calR\left(G^0(V_n)\right).
\]
Obviously $\calR_V$ only depends on the Witt tower containing $V$. When there is no confusion, we write $G^0_n$ for $G^0(V_n)$ for brevity.

\vskip 5pt

For a parabolic subgroup \(P\) of \(\GL_n\) or \(G^0_n\), we denote the (normalized) Jacquet functor along \(P\) by \(\Jac_{P}\). We denote by \([\pi]\) the semi-simplification of a finite length representation $\pi$. The Jacquet functors on general linear groups give rise to a co-multiplication structure on \(\calR\):
\begin{align*}
    m^*:\calR&\longrightarrow\calR\otimes\calR\\
    \tau&\longmapsto\sum_{n_1+n_2=n}\left[\Jac_{P_{n_1,n_2}}(\tau)\right],
\end{align*}
where $\tau\in\calR(\GL_n)$, and $P_{n_1,n_2}$ is the maximal standard parabolic subgroup of \(\GL_n\) with Levi component isomorphic to \(\GL_{n_1}\times\GL_{n_2}\). On the other hand, the Jacquet functors on classical groups induce a map
\begin{align*}
    \mu^*:\calR_V&\longrightarrow\calR\otimes\calR_V\\
    \pi&\longmapsto\sum^n_{k=0}\left[\Jac_{P_k}(\pi)\right],
\end{align*}
where $\pi\in \calR\left(G^0_n\right)$, and \(P_k\) is the maximal standard parabolic subgroup with the Levi component \(\GL_k\times G_{n-2k}\).
The following result, known as \emph{Tadi\'c's formula}, helps us to compute \(\mu^*\) inductively. It is proved in \cite{MR1356358} for symplectic and special odd orthogonal groups, in \cite{MR1739616} for even orthogonal groups, and in \cite{MR1896238} for unitary groups.

\begin{Thm}[Tadi\'c's formula]\label{tadic}
Let \(M^*:\calR\rightarrow \calR\otimes \calR\) be the composition of maps
\[
    M^*=(m\otimes\mathrm{id})\circ( \sim\otimes m^*)\circ s\circ m^*,
\] 
where \(m:\calR\otimes\calR\rightarrow\calR\) is the multiplication map on \(\calR\) induced by parabolic inductions, and \(s:\calR\otimes\calR\rightarrow\calR\otimes\calR\) is the switching map, namely, $s\left(\tau\otimes \tau'\right) = \tau'\otimes\tau$. Then, we have
\[
    \mu^*(\tau\rtimes\sigma)=M^*(\tau)\rtimes\mu^*(\sigma).
\]
\end{Thm}


Next we discuss the notion of \emph{derivatives}. For general linear groups, one has the notion of \emph{Bernstein--Zelevinsky derivatives}, which can be used to describe the restrictions of representations to mirabolic subgroups \cite[Prop. 3.2, 3.5]{MR579172}.

\begin{Lem}\label{L:RestrictionToMirabolic}
Let $\tau$ be an irreducible representation of $\GL_n$, and $M_n$ the mirabolic subgroup of $\GL_n$. Then there is an $M_n$-equivariant filtration
\[
    \tau\big|_{M_n} = \mu_1 \supset \mu_2 \supset \cdots\supset 0,
\]
with successive quotients
\[
    \mu_\alpha/\mu_{\alpha+1} \cong \prescript{u}{}{\ind}_{T_\alpha\rtimes N_\alpha}^{M_n}\left(|\cdot|^{\frac{\alpha}{2}}\tau^{(\alpha)}\boxtimes\psi_\alpha\right),
\]
where:
\begin{itemize}
    \item $\prescript{u}{}{\ind}$ is the un-normalized compactly supported induction;
    \vskip 5pt
    
    \item $T_\alpha$ is a subgroup of $M_n$ isomorphic to $\GL_{n-\alpha}\times U_\alpha$, with $U_\alpha$ the subgroup of upper uni-triangular matrices in $\GL_\alpha$;
    \vskip 5pt

    \item $N_\alpha$ is the unipotent radical of $P_\alpha$, where $P_\alpha$ is the standard parabolic subgroup of $\GL_{n}$ with Levi component $\GL_{n-\alpha}\times \GL_\alpha$;
    \vskip 5pt

    \item $\tau^{(\alpha)}$ is a representation of $\GL_{n-\alpha}$, called the $\alpha$-th Bernstein-Zelevinsky derivative of $\tau$;
    \vskip 5pt

    \item $\psi_\alpha$ is a generic character of $U_\alpha$.
\end{itemize}
\end{Lem}

\vskip 5pt

For example, for a fixed $\rho\in\calC_{\unit}(\GL_d)$, if $\tau = \Delta_\rho[x,y]$, then $\tau^{(\alpha)}=\Delta_\rho[x,y+k]$ when $\alpha=kd$ for some $k\in\Z_{>0}$ and $\tau^{(\alpha)}=0$ otherwise. Similarly, if $\tau = Z_\rho[x,y]$, then $\tau^{(\alpha)}=Z_\rho[x-1,y]$ when $\alpha=d$ and $\tau^{(\alpha)}=0$ otherwise. 

\vskip 5pt

We will also make use of another notion of derivatives defined using (semi-simplified) Jacquet functors. Fix an irreducible representation $\tau$ of $\GL_d$.

\begin{Def}
Retain the notation from the above discussion. Let $\pi$ be a representation of $\GL_n$. The \(k\)-th \emph{left} (resp. \emph{right}) \(\tau\)-\emph{(semi-simplified) derivative} \(L^{(k)}_{\tau}(\pi)\) (resp. \(R^{(k)}_{\tau}(\pi)\)) is the unique semi-simple representation of \(\GL_{n-dk}\) such that in the Grothendieck group $\calR(\GL_{dk}\times\GL_{n-dk})$ (resp. $\calR(\GL_{n-dk}\times\GL_{dk})$), $\left[\Jac_{P_{dk,n-dk}}(\pi)\right]$ (resp. $\left[\Jac_{P_{n-dk,dk}}(\pi)\right]$) can be written as
\[
    \left[\Jac_{P_{dk,n-dk}}(\pi)\right] = \tau^k\boxtimes L^{(k)}_{\tau}(\pi)+ \sum_{i}\tau_i\boxtimes\pi_i
\]
\[
    \left(\text{ resp. } \left[\Jac_{P_{n-dk,dk}}(\pi)\right] = R^{(k)}_{\tau}(\pi)\boxtimes\tau^k + \sum_{i}\pi_i\boxtimes\tau_i\right),
\]
where for each index $i$, $\tau_i$ is an irreducible representation of $\GL_{dk}$ and $\pi_i$ is an irreducible representation of $\GL_{n-dk}$, such that \(\tau_i\not\cong\tau^k\).
\end{Def} 

\vskip 5pt

For classical groups, one can introduce the notion of derivatives in the same manner.

\begin{Def}
Retain the notation from the above discussion. Let $\pi$ be a representation of $G^0_n$. The \(k\)-th \(\tau\)-\emph{(semi-simplified) derivative} \(D^{(k)}_{\tau}(\pi)\) is the unique semi-simple representation of \(G^0_{n-2dk}\) such that in the Grothendieck group $\calR(\GL_{dk}\times G^0_{n-2dk})$, $\left[\Jac_{P_{dk}}(\pi)\right]$ can be written as
\[
    \left[\Jac_{P_{dk}}(\pi)\right]= \tau^k\boxtimes D^{(k)}_{\tau}(\pi) + \sum_i\tau_i\boxtimes\pi_i,
\]
where for each index $i$, $\tau_i$ is an irreducible representation of $\GL_{dk}$ and $\pi_i$ is an irreducible representation of $G^0_{n-2dk}$, such that \(\tau_i\not\cong\tau^k\). We say that \(D^{(k)}_{\tau}(\pi)\) is the highest derivative if \(D^{(k)}_{\tau}(\pi)\) is non-zero but \(D^{(k+1)}_{\tau}(\pi)\) is zero.
\end{Def}

\vskip 5pt

When $k=1$, we shall suppress the superscript ``$(k)$'' from the notation. In applications, the primary case of interest is when $\tau\cong\rho|\cdot|^x$ is supercuspidal. In the following, we set \(\calC^{\perp}(\GL_n)\) to be the subset of \(\calC_{\unit}(\GL_n)\) consisting of irreducible representations that are (conjugate) self-dual.
The following proposition lists some basic formulas and properties when \(\tau\) is supercuspidal.

\begin{Prop}[{\cite[Sect. 5]{MR3713922}}]\label{leibniz}
Let \(\rho\in\calC_{\unit}(\GL_d)\) and \(x\in\R\).
\begin{enumerate}
        \item For each segment $[a,b]_{\rho'}$, we have
        \[
            L_{\rho|\cdot|^x}(\Delta_{\rho'}[a,b])=\begin{cases}
                \Delta_{\rho'}[a-1,b] \quad &\text{ if }\rho|\cdot|^x\cong\rho'|\cdot|^a;\\
                0 \quad &\text{ otherwise.}
            \end{cases}
        \]
        Similarly, we have
        \[
            R_{\rho|\cdot|^x}(\Delta_{\rho'}[a,b])=\begin{cases}
                \Delta_{\rho'}[a,b+1] \quad &\text{ if }\rho|\cdot|^x\cong\rho'|\cdot|^b;\\
                0 \quad &\text{ otherwise.}
            \end{cases}
        \]

        \vskip 5pt
        
        \item (Leibniz rule) For finite length representations \(\pi\) of \(G^0_n\) and \(\tau,\tau_1,\tau_2\) of general linear groups, we have
        \begin{itemize}
            \item \(L_{\rho|\cdot|^x}(\tau_1\times\tau_2)=L_{\rho|\cdot|^x}(\tau_1)\times\tau_2+\tau_1\times L_{\rho|\cdot|^x}(\tau_2)\);

            \vskip 5pt

            \item \(R_{\rho|\cdot|^x}(\tau_1\times\tau_2)=R_{\rho|\cdot|^x}(\tau_1)\times\tau_2+\tau_1\times R_{\rho|\cdot|^x}(\tau_2)\);

            \vskip 5pt

            \item \(D_{\rho|\cdot|^x}(\tau\rtimes\pi)=L_{\rho|\cdot|^x}(\tau)\rtimes\pi + R_{\widetilde{\rho}|\cdot|^{-x}}(\tau)\rtimes\pi+\tau\rtimes D_{\rho|\cdot|^x}(\pi)\). 
        \end{itemize}
    \end{enumerate}
\end{Prop}
\vskip 5pt

We note that the Leibniz rule is a consequence of Theorem \ref{tadic}. For classical groups, there is an algorithm to compute highest $\tau$-(semi-simplified) derivatives when $\tau$ is non (conjugate) self-dual supercuspidal or $\tau\cong \Delta_\rho[0,-1]$ for some $\rho\in\calC^\perp(\GL_d)$. We refer the readers to \cite{MR4549708} and \cite{clw26} for details. For convenience we introduce the following notation: if $\tau\hookrightarrow \rho_1\times\rho_2\times\cdots\times\rho_r$ for some supercuspidal representations $\rho_i$ of general linear groups, we set
\[
    D_{\csupp(\tau)} \coloneqq D_{\rho_r}\circ \cdots \circ D_{\rho_2}\circ D_{\rho_1},
\]
and set $L_{\csupp(\tau)}$ and $R_{\csupp(\tau)}$ similarly.

\vskip 5pt

For the purposes of subsequent sections, we also introduce another notion of derivatives, which provides more refined information than the previous semi-simplified one.

\begin{Def}
Retain the notation from the above discussion. Let $\pi$ be a representation of $G^0_n$. For each $i\in\Z_{\geq 0}$, we set the $i$-th homological $\tau$-\emph{derivative} of $\pi$ to be
\[
    \D^i_\tau (\pi) \coloneqq {\rm H}_i\left(\GL_d,\tau^\vee\otimes \Jac_{P_d}(\pi)\right),
\] 
namely, the $i$-th (left) derived functor of $\left(\tau^\vee\otimes \Jac_{P_d}(-)\right)_{\GL_d}$, where the subscript $\GL_d$ means taking $\GL_d$-coinvariants. To simplify notation, when $i=0$ we denote $\D_\tau(\pi) \coloneqq \D_\tau^0(\pi)$.
\end{Def}

\vskip 5pt

We remark that in \cite{chan2026}, Chan and Pattanayak also introduced similar notions of (left and right) derivatives $\mathbb{L}_{\Delta_{\rho}[x,y]}$ and $\mathbb{R}_{\Delta_{\rho}[x,y]}$ for general linear groups. They also provided an algorithm to compute these derivatives. By using these notions, Chan was able to extend the non-tempered GGP conjecture for general linear groups to \emph{all} irreducible representations, see \cite{chan2022quotient} for more details. The following relationship between the two notions of derivatives, $\D_\tau(\pi)$ and $D_\tau(\pi)$, is worth noting.

\begin{Lem}
Let $\pi$ be a finite length representation of $G^0_n$. Then we have inclusions 
\[
    \left[\D_\tau(\pi)\right] \subset \left[D_\tau(\pi)\right] \subset \left[D_{\csupp(\tau)}(\pi)\right].
\]
In particular, $\D_\tau(\pi)\neq 0$ implies that $D_\tau(\pi)\neq 0$ and $D_{\csupp(\tau)}(\pi)\neq 0$. Conversely, if $\pi$ is irreducible and $\tau$ is supercuspidal, then $D_\tau(\pi)\neq 0$ implies that $\D_\tau(\pi)\neq 0$.   
\end{Lem}

\vskip 5pt

\begin{proof}
Let $M_d\cong \GL_d\times G^0_{n-2d}$ be the standard Levi subgroup of $P_d$. If $\D_\tau(\pi)\neq 0$, then there is a natural surjective $M_d$-equivariant map
\[
    \Jac_{P_d}(\pi) \twoheadrightarrow \tau\boxtimes \D_\tau(\pi).
\]
This implies that $\left[\D_\tau(\pi)\right] \subset \left[D_\tau(\pi)\right]$. To see that $\left[D_\tau(\pi)\right] \subset \left[D_{\csupp(\tau)}(\pi)\right]$, one can simply take Jacquet modules in stages. Finally, when $\pi$ is irreducible and $\tau$ is supercuspidal, by \cite[Lem. 5.3]{MR3713922} $D_\tau(\pi)\neq 0$ implies that there is some irreducible representation $\pi_0$ of an appropriate classical group $G^0_{n_0}$, such that
\[
    \pi \hookrightarrow \tau \rtimes \pi_0.
\]
Therefore $\pi_0$ is a quotient of $\D_\tau(\pi)$, and hence $\D_\tau(\pi)\neq 0$.

\end{proof}

\vskip 5pt

We conclude this section with an estimation of exponents in Jacquet modules, which follows from the Langlands classification. We also give two applications of this estimation.

\begin{Lem}\label{L:ExpInJacModLowerBound}
Let $\pi=\LS\left(\Delta_{\rho_1}[x_1,y_1],\cdots,\Delta_{\rho_r}[x_r,y_r];\pi_0\right)$ be as in Theorem \ref{lcgn} (here we allow $r=0$, in which case $x_1+y_1$ should be understood as $0$). Regard $\pi$ as a representation of $G^0_n$ by restriction; it is still irreducible. If we write
\[
    \left[\Jac_{P_d}(\pi)\right] = \sum_\alpha \tau_\alpha\boxtimes\pi_\alpha,
\]
where for each index $\alpha$, $\tau_\alpha$ is an irreducible representation of $\GL_{d}$ and $\pi_\alpha$ is an irreducible representation of $G^0_{n-2d}$. Then the exponents $e_{\tau_\alpha}$ of these $\tau_\alpha$ satisfy $e_{\tau_\alpha}\geq \frac{x_1+y_1}{2}$. 
\end{Lem}

\vskip 5pt

\begin{proof}
We prove by induction on the non-negative integer $r$. The base case is that $r=0$, namely $\pi$ is tempered. In this case the assertion in this lemma is exactly the same as Casselman's criterion. Now suppose that $r>0$. We have an embedding
\[
    \pi\hookrightarrow \Delta_{\rho_1}[x_1,y_1] \rtimes \pi_1^-,
\]
where $\pi^-_1\coloneqq\LS\left(\Delta_{\rho_2}[x_2,y_2],\cdots,\Delta_{\rho_r}[x_r,y_r];\pi_0\right)$. Apply Theorem \ref{tadic} to the parabolically induced representation $\Delta_{\rho_1}[x_1,y_1] \rtimes \pi_1^-$. Similar to \cite[Cor. 2.9]{MR4055180}, we have
\begin{align*}
\mu^*\left(\Delta_{\rho_1}[x_1,y_1] \rtimes \pi_1^-\right) = \sum_{\substack{k,\ell\geq 0;\\ k+\ell\leq x_1-y_1+1}} &\Big(\left(\Delta_{\widetilde{\rho}_1}\left[-y_1,-y_1-k+1\right]\times\Delta_{\rho_1}\left[x_1,x_1-\ell+1\right]\right) \\
&\otimes \Delta_{\rho_1}\left[x_1-\ell, y_1+k\right]\Big)\rtimes \mu^*\left(\pi^-_1\right).   
\end{align*}
Note that
\[
    \frac{-y_1 + (-y_1-k+1)}{2} > \frac{x_1 + y_1}{2} \quad \text{and} \quad \frac{x_1 + (x_1-\ell+1)}{2} \geq \frac{x_1 + y_1}{2},
\] 
and by the induction hypothesis, any irreducible summand $\tau_1^-\otimes\sigma_1^-$ of $\mu^*\left(\pi^-_1\right)$ has exponent $e_{\tau_1^-}\geq \frac{x_2+y_2}{2}\geq\frac{x_1 + y_1}{2}$. Since for any index $\alpha$, there exists an irreducible summand $\tau_1^-\otimes\sigma_1^-$ of $\mu^*\left(\pi^-_1\right)$ such that $\tau_\alpha\subset \left[\Delta_{\widetilde{\rho}_1}\left[-y_1,-y_1-k+1\right]\times\Delta_{\rho_1}\left[x_1,x_1-\ell+1\right]\times \tau^-_1\right]$, we also have $e_{\tau_\alpha}\geq \frac{x_1+y_1}{2}$ as desired.

\end{proof}

\vskip 5pt

The first application is a vanishing criterion for $D_\tau(\pi)$ and $\D^i_\tau(\pi)$.

\begin{Cor}\label{C:VanishingLemma}
In the context of Lemma \ref{L:ExpInJacModLowerBound}, let $\tau$ be an irreducible representation of $\GL_d$. If the exponent $e_\tau$ of $\tau$ satisfies $e_\tau<\frac{x_1+y_1}{2}$, then $D_\tau(\pi)=0$ and $\D^i_{\tau}\pi = 0$ for any $i\geq 0$.
\end{Cor}

\vskip 5pt

\begin{proof}
By Lemma \ref{L:ExpInJacModLowerBound}, one knows that
\[
    \left[\Jac_{P_d}(\pi)\right] = \sum_\alpha \tau_\alpha\boxtimes\pi_\alpha,
\]
where for each index $\alpha$, $\tau_\alpha$ is an irreducible representation of $\GL_{d}$ and $\pi_\alpha$ is an irreducible representation of $G^0_{n-2d}$, such that their exponents satisfy $e_{\tau_\alpha}\geq \frac{x_1+y_1}{2}$. Then by definition $D_\tau(\pi)=0$. Moreover, for each index $\alpha$, one has
\[
    {\rm H}_i\left(\GL_d, \tau^\vee\otimes \left(\tau_\alpha\boxtimes\pi_\alpha\right)\right) \cong {\rm H}_i\left(\GL_d, \tau^\vee\otimes\tau_\alpha\right)\otimes\pi_\alpha =0.
\]
Here in the last equality we have made use of the fact that $e_\tau\neq e_{\tau_\alpha}$. It follows that
\[
    \D_\tau^i(\pi) \coloneqq {\rm H}_i\left(\GL_d, \tau^\vee\otimes \Jac_{P_d}(\pi)\right)=0
\]
as desired. This completes the proof.

\end{proof}

\vskip 5pt

The second application is a Leibniz rule for Steinberg derivatives in a very special situation.

\begin{Cor}\label{leibniz-st}
In the context of Lemma \ref{L:ExpInJacModLowerBound}, let $\tau\cong\Delta_\rho[x,y]$ be an irreducible Steinberg representation of $\GL_d$, and $\delta\cong\Delta_{\rho_0}[x_0,y_0]$ an essentially discrete series representation of a general linear group. If
\[
    x+y\leq x_0+y_0\leq x_1+y_1 (\leq 0),
\]
then we have
\[
    D_\tau\left(\delta\rtimes\pi\right) = L_\tau\left(\delta\right)\rtimes \pi + R_{\widetilde{\tau}}\left(\delta\right)\rtimes \pi + \delta\rtimes D_\tau\left(\pi\right).
\]
\end{Cor}

\vskip 5pt

\begin{proof}
Obviously the right-hand side of the formula is contained in the left-hand side, so if $D_\tau\left(\delta\rtimes\pi\right)=0$ the formula holds automatically. Suppose that $D_\tau\left(\delta\rtimes\pi\right)\neq 0$. Apply Theorem \ref{tadic} to the parabolically induced representation $\delta\rtimes\pi$. Similar to \cite[Cor. 2.9]{MR4055180}, we have
\begin{align*}
\mu^*\left(\delta\rtimes\pi\right) = \sum_{\substack{k,\ell\geq 0;\\ k+\ell\leq x_0-y_0+1}} &\left(\left(\Delta_{\widetilde{\rho}_0}\left[-y_0,-y_0-k+1\right]\times\Delta_{\rho_0}\left[x_0,x_0-\ell+1\right]\right)\otimes \Delta_{\rho_0}\left[x_0-\ell, y_0+k\right]\right) \\
&\rtimes \mu^*\left(\pi\right).   
\end{align*}
We are interested in those irreducible summands of $\mu^*\left(\delta\rtimes\pi\right)$ of the form $\tau\otimes \sigma$. If  
\[
    \tau\otimes \sigma \subset \left(\left(\Delta_{\widetilde{\rho}_0}\left[-y_0,-y_0-k+1\right]\times\Delta_{\rho_0}\left[x_0,x_0-\ell+1\right]\right)\otimes \Delta_{\rho_0}\left[x_0-\ell, y_0+k\right]\right) \rtimes \mu^*\left(\pi\right),
\]
then there is an irreducible summand $\tau_0\otimes\sigma_0$ of $\mu^*(\pi)$, such that
\begin{equation*}\label{E:LeibnizSteinbergTadicTermkl}
    \tau\subset \left[\Delta_{\widetilde{\rho}_0}\left[-y_0,-y_0-k+1\right]\times\Delta_{\rho_0}\left[x_0,x_0-\ell+1\right]\times\tau_0\right].
\end{equation*}
Note that $\frac{-y_0 + (-y_0-k+1)}{2} \geq \frac{x_0 + y_0}{2}\geq e_\tau$ and $\frac{x_0 + (x_0-\ell+1)}{2} \geq \frac{x_0 + y_0}{2}\geq e_\tau$. Also, by Corollary \ref{C:VanishingLemma} one has
\[
    e_{\tau_0}\geq \frac{x_1+y_1}{2} \geq e_\tau.
\]
The existence of such a summand $\tau\otimes \sigma$ of $\mu^*(\delta\rtimes\pi)$ implies that these inequalities of exponents must all be equalities. In particular, $k$ and $\ell$ are either $0$ or $x_0-y_0+1$. We now analyze case by case the $\tau$-isotypic component of each term for each possible pair $(k,\ell)$.

\vskip 5pt

\textit{Case 1:} $k=0$ and $\ell=0$. Then the $\tau$-isotypic component of the corresponding term is exactly $\delta\rtimes D_\tau(\pi)$.

\vskip 5pt

\textit{Case 2:} $k=0$ and $\ell=x_0-y_0+1$. We claim that the $\tau$-isotypic component of the corresponding term vanishes whenever $[x,y]_\rho\neq [x_0,y_0]_{\rho_0}$. Indeed, given that $x+y=x_0+y_0$, there are three possibilities: $[x,y]_\rho\cap [x_0,y_0]_{\rho_0}=\varnothing$, $[x,y]_\rho\subsetneq [x_0,y_0]_{\rho_0}$, or $[x_0,y_0]_{\rho_0}\subsetneq [x,y]_\rho$. If $[x,y]_\rho\cap [x_0,y_0]_{\rho_0}=\varnothing$ or $[x,y]_\rho\subsetneq [x_0,y_0]_{\rho_0}$, the claim is immediate. If $[x_0,y_0]_{\rho_0}\subsetneq [x,y]_\rho$, then we must have $y<y_0\leq x_0<x$. The cuspidal support of $\tau_0$ in the above discussion is contained in $[x,x_0+1]_\rho\sqcup[y_0-1,y]_\rho$. It follows that there exist two representations $\tau_0^+$ and $\tau_0^-$, such that the cuspidal support of $\tau_0^+$ is contained in $[x,x_0+1]_\rho$, the cuspidal support of $\tau_0^-$ is contained in $[y_0-1,y]_\rho$, and
\[
    \tau_0 \hookrightarrow \tau_0^- \times \tau_0^+.
\]
By taking Jacquet modules in stages, it follows that $D_{\tau_0^-}(\pi)\neq 0$. This contradicts Corollary \ref{C:VanishingLemma}. Therefore our claim holds. Note that when $[x,y]_\rho\neq [x_0,y_0]_{\rho_0}$, one also has $L_\tau(\delta)=0$. On the other hand, when $[x,y]_\rho = [x_0,y_0]_{\rho_0}$, the $\tau$-isotypic component of the corresponding term is $\pi\cong L_\tau(\delta)\rtimes \pi$.

\vskip 5pt

\textit{Case 3:} $k=x_0-y_0+1$ and $\ell=0$. By using the same argument as in the previous case, one can conclude that the $\tau$-isotypic component of the corresponding term is $R_{\widetilde{\tau}}\left(\delta\right)\rtimes\pi$.

\vskip 5pt

The desired formula follows from the combination of all these three cases.

\end{proof}

\vskip 10pt

\section{L-packets, A-packets and representations of Arthur type} 
\label{sec:l_packets_a_packets_and_representations_of_arthur_type}

In this section, we briefly recall some background on the local Langlands correspondence (``LLC'' for short) and Arthur's endoscopic classification. As in Section \ref{subsec:TadicFormulaDerivative}, we shall work with $\GL_n$ and $G^0(V)$.

\subsection{Local Langlands correspondence}

Let $W_E$ be the Weil group of $E$ and $WD_E \coloneqq W_E\times\SL_2(\C)$ the Weil--Deligne group of $E$. We say that a continuous homomorphism $\phi: WD_E\ra \GL_n(\C)$ is a representation of $WD_E$, if
\begin{itemize}
    \item $\phi$ is semi-simple, and

    \vskip 5pt

    \item $\phi\big|_{\SL_2(\C)}$ is algebraic.
\end{itemize}
\vskip 5pt
We say that $\phi$ is \textit{tempered} if the image of $W_E$ is bounded. Given a representation $\phi$ of $WD_E$, one can define the notion of contragredient $\phi^\vee$, $\iota$-conjugate $\prescript{\iota}{}{\phi}$, (conjugate) self-duality and parity as in \cite[Sect. 3]{MR3202556}. Typically, we shall denote an irreducible representation of $WD_E$ by $\rho|\cdot|^x S_a$, where $\rho$ is a tempered irreducible representation of $W_E$, $x\in\R$ and $S_a$ the unique irreducible $a$-dimensional algebraic representation of $\SL_2(\C)$. As in (\ref{E:TildeRepDef}), for simplicity of notation we set
\[
    \widetilde{\phi} \coloneqq \begin{cases}
     \phi^\vee \quad&\text{ if }E=F;\\
\prescript{\iota}{}\phi^\vee \quad&\text{ otherwise.}   
    \end{cases}
\]

\vskip 5pt

Let $\widetilde{\Phi}\left(n\right)$ be the set of equivalence classes of $n$-dimensional representations of $WD_E$. The LLC for general linear groups, which was established in \cite{MR1876802, MR1738446}, is a bijection
\[
    \Irr\left(\GL_n\right) \longleftrightarrow \widetilde{\Phi}\left(n\right)
\]
satisfying some natural properties. For example, if $\tau\in\Irr\left(\GL_n\right)$ and $\phi$ is the $n$-dimensional representation of $WD_E$ associated to $\tau$, then $\tau$ is supercuspidal if and only if $\phi\big|_{W_E}$ is irreducible, $\tau$ is an essentially discrete series representation if and only if $\phi$ is irreducible, and $\tau$ is tempered if and only if $\phi$ is tempered. By abuse of notation, we shall not distinguish between $\rho\in\calC_{\unit}\left(\GL_n\right)$ and its L-parameter, i.e. the irreducible representation of $W_E$ associated to $\rho$. For a representation $\phi$ of $WD_E$, we recall that its \emph{infinitesimal character} is the $W_E$-representation defined by the formula
\[
    \lambda_\phi(x)\coloneqq \phi\left(x,\left(\begin{array}{cc}
                |x|^{\frac{1}{2}} & \\
                 & |x|^{-\frac{1}{2}}
            \end{array}\right)\right).
\]

\vskip 5pt

Now let $V$ be an $n$-dimensional formed space over $E$ as described in Section \ref{subsec:FormedSpaceNIsomGroup}, and consider the group $G^0(V)$. The LLC (as enhanced by Vogan) for classical groups, which was established in \cite{MR3135650, MR3338302, MR3708200, MR4308934, MR4324358, MR4776199}, asserts that there is a finite-to-one surjection
\[
    \calL: \bigsqcup_{V'} \Irr\left(G^0(V')\right) \lra \Phi\left(G^0(V)\right),
\]
where $V'$ runs over all $n$-dimensional $\varepsilon$-Hermitian spaces over $E$ such that $G^0(V')$ is a pure inner form of $G^0(V)$, and $\Phi\left(G^0(V)\right)$ is the set of equivalence classes of L-parameters of $G^0(V)$. We recall that an L-parameter of $G^0(V)$ is a continuous homomorphism $\phi: WD_F \ra \prescript{L}{}{G^0(V)}$ satisfying certain conditions. Following \cite[Sect. 8]{MR3202556}, by composing $\phi$ with the standard representation of $\prescript{L}{}{G^0(V)}$ and restricting to $WD_E$, we can regard an L-parameter $\phi$ of $G^0(V)$ as a (conjugate) self-dual representation of $WD_E$ of dimension
\[
    \begin{cases}
    n+1 \quad & \text{if $G^0(V)$ is a symplectic group},\\[5pt]
    n-1 \quad & \text{if $G^0(V)$ is a special odd orthogonal group},\\[5pt]
    n \quad & \text{if $G^0(V)$ is an even orthogonal group or unitary group}
    \end{cases}
\]
and parity
\[
    \begin{cases}
    +1 \quad & \text{if $G^0(V)$ is a symplectic group or even orthogonal group},\\[5pt]
    -1 \quad & \text{if $G^0(V)$ is a special odd orthogonal group},\\[5pt]
    (-1)^{n-1} \quad & \text{if $G^0(V)$ is a unitary group}.
    \end{cases}
\]
In the case that $G^0(V)$ is a symplectic group, $\phi$ also satisfies $\det\left(\phi\right) = \mathbbm{1}$; whereas in the case that $G^0(V)$ is a (special) orthogonal group, $\phi$ satisfies $\det\left(\phi\right) = \left(\disc(V),-\right)_E$. For each $\phi\in\Phi\left(G^0(V)\right)$, we put $\LP_\phi\left(G^0(V)\right)\coloneqq \calL^{-1}\left(\phi\right)\cap \Irr\left(G^0(V)\right)$; this is the local L-packet of $G^0\left(V\right)$ associated to the L-parameter $\phi$. The fiber 
\begin{equation*}
    \LP_\phi \coloneqq \calL^{-1}\left(\phi\right) = \bigsqcup_{V'} \LP_\phi\left(G^0(V')\right)
\end{equation*}
of $\calL$ at $\phi$ is called the Vogan L-packet of $G^0(V)$ associated to the L-parameter $\phi$. Since $\phi$ is (conjugate) self-dual, we can write $\phi$ as
\begin{equation}\label{E:DecompLparameter}
    \phi = \varphi + \sum_i m_i\phi_i + \widetilde{\varphi},
\end{equation}
where for each index $i$, $\phi_i$ is an irreducible (conjugate) self-dual representation of $WD_E$ of the same parity as $\phi$, these $\phi_i$ are pairwise distinct, $m_i\in \Z_{>0}$ is the multiplicity of $\phi_i$ in $\phi$, and any irreducible constituent of $\varphi$ is not (conjugate) self-dual of the same parity as $\phi$. For any irreducible tempered representation $\rho S_a$ of $WD_E$, set
\[
    m_\phi\left(\rho S_a\right) \coloneqq \dim \Hom_{WD_E}\left(\rho S_a, \phi\right).
\]
We say that:
\begin{itemize}
    \item $\phi$ is of \emph{good parity} if $\varphi=0$;

    \vskip 5pt

    \item $\phi$ is \emph{discrete} if it is of good parity and $m_i=1$ for all $i$.
\end{itemize}
Similar to the case of general linear groups, this map $\calL$ enjoys many good properties. For example, if $\pi\in\Irr\left(G^0(V)\right)$ and $\phi$ is the L-parameter associated to $\pi$, then $\pi$ is a discrete series if and only if $\phi$ is discrete, and $\pi$ is tempered if and only if $\phi$ is tempered (as a representation of $WD_E$). 

\vskip 5pt

For an L-parameter $\phi$ of $G^0(V)$, following Arthur \cite[Sect. 1.4]{MR3135650} we define its component group to be
\[
    \calS_\phi \coloneqq \pi_0\left({\rm Cent}_{\widehat{G^0(V)}}\left({\rm Im}\left(\phi\right)\right)\right).
\]
Here $\widehat{G^0(V)}$ is the Langlands dual group of $G^0(V)$; when $G^0(V)$ is an even orthogonal group, we understand $\widehat{G^0(V)}$ as $\rO_n(\C)$. Then there is a canonical injection (depending on the choice of a Whittaker datum $\mathscr{W}$ of a quasi-split pure inner form of $G^0(V)$)
\[
    \calJ^{\rm L}_{\mathscr{W}}: \LP_\phi \lra \Irr\left(\calS_\phi\right),
\]
and it is a bijection unless $G^0(V)$ is a symplectic group. As explicated in \cite[Sect. 8]{MR3202556}, this component group $\calS_\phi$ can be combinatorially described. Suppose that $\phi$ can be decomposed as in (\ref{E:DecompLparameter}). Let
\[
    A_\phi \coloneqq \prod_i \left(\Z/2\Z\right) \alpha_i
\]
be the elementary $2$-group with a canonical basis $\left\{\alpha_i\right\}_i$, where $\alpha_i$ corresponds to the irreducible constituent $\phi_i$. 
Then we have $\calS_\phi\cong A_\phi$. Having fixed the Whittaker datum $\mathscr{W}$, we shall write $\pi\left(\phi,\eta\right)$ for the irreducible representation in $\LP_\phi$ corresponding to $\eta\in\Irr\left(A_\phi\right)$. The following properties of the LLC are worth noting.

\begin{Thm}\label{T:LLCProperties}
The LLC (as enhanced by Vogan) for $G^0(V)$ enjoys:
\begin{enumerate}
    \item Local intertwining relation (``LIR'' for short): Let $\phi$ be a tempered L-parameter of $G^0(V)$. If $\phi = \phi_\tau + \phi_0 + \widetilde{\phi}_\tau$, where $\phi_\tau$ is an irreducible tempered representation of $WD_E$ corresponding to an irreducible unitary discrete series $\tau$ under the LLC for general linear groups, then $\LP_\phi = \left\{\pi~\big|~\exists\,\eta_0\in\Irr\left(A_{\phi_0}\right),~\pi\subset\tau\rtimes\pi\left(\phi_0,\eta_0\right)\right\}$, and
    \[
        \tau\rtimes\pi\left(\phi_0,\eta_0\right) = \bigoplus_{\eta|_{A_{\phi_0}}=\eta_0} \pi\left(\phi,\eta\right).
    \]
    Moreover, if $\phi_\tau$ is (conjugate) self-dual of the same parity as $\phi$, then $\eta(\alpha_\tau)$ can be determined by the eigenvalue of a certain normalized intertwining operator on $\pi\left(\phi,\eta\right)$, where $\alpha_\tau$ is the element in the canonical basis of $A_\phi$ corresponding to $\phi_\tau$.

    \vskip 5pt

    \item Compatibility with the Langlands subrepresentation: Let $\phi$ be an L-parameter of $G^0(V)$. If
    \[
        \phi = \left(\phi_1+\phi_2+\cdots+\phi_r\right) + \phi_0 + \left(\widetilde{\phi}_1+\widetilde{\phi}_2+\cdots+\widetilde{\phi}_r\right),
    \]
    where for each $1\leq i\leq r$, $\phi_i$ is an irreducible representation of $WD_E$ which corresponds to an essentially discrete series $\Delta_{\rho_i}[x_i,y_i]$, such that
    \[
        x_1+y_1\leq x_2+y_2\leq \cdots \leq x_r+y_r<0;
    \]
    $\phi_0$ is a tempered L-parameter. Then for any $\eta\in\Irr\left(A_\phi\right)$, we have
    \[
        \pi\left(\phi,\eta\right) = \LS\left(\Delta_{\rho_1}[x_1,y_1],\Delta_{\rho_2}[x_2,y_2],\cdots,\Delta_{\rho_r}[x_r,y_r];\pi\left(\phi_0,\eta\big|_{A_{\phi_0}}\right)\right).
    \]
\end{enumerate}    
\end{Thm}

\vskip 5pt

\subsection{A-packets and representations of Arthur type}

Next we turn to Arthur's endoscopic classification. For global purposes, Arthur introduced the notion of A-parameters and A-packets.

\vskip 5pt

For the general linear group $\GL_n$, a local A-parameter is an $n$-dimensional representation of $WD_E\times \SL_2(\C)$ with \emph{bounded image} on the Weil group $W_E$. To distinguish the two copies of $\SL_2(\C)$, we shall call the one in the Weil--Deligne group ``\emph{Deligne} $\SL_2(\C)$'', and the extra one ``\emph{Arthur} $\SL_2(\C)$''. The set of equivalence classes of local A-parameters of $\GL_n$ will be denoted by $\widetilde{\Psi}(n)$. Given a local A-parameter $\psi\in\widetilde{\Psi}(n)$, we write
\[
    \psi = \sum_i m_i \rho_i S_{a_i}\boxtimes S_{b_i},
\]
where for each index $i$, $\rho_i$ is an irreducible tempered representation of $W_E$, $S_a$ is the unique irreducible $a$-dimensional algebraic representation of Deligne $\SL_2(\C)$, and $S_b$ is the unique irreducible $b$-dimensional algebraic representation of Arthur $\SL_2(\C)$. Then one can associate an irreducible unitary representation
\[
    \tau_\psi \coloneqq \bigtimes_i u_{\rho_i}\left(a_i,b_i\right)^{\times m_i} 
\]
to $\psi$, where the superscript ``$\times m_i$'' is the abbreviation for $\underbrace{u_{\rho_i}\left(a_i,b_i\right)\times\cdots\times u_{\rho_i}\left(a_i,b_i\right)}_{m_i\text{-times}}$.

\vskip 5pt

Now let $V$ be an $n$-dimensional formed space over $E$ as described in Section \ref{subsec:FormedSpaceNIsomGroup}. Similar to the case of general linear groups, for the group $G^0(V)$, a local A-parameter is a homomorphism $\psi: WD_F\times \SL_2(\C) \ra \prescript{L}{}{G^0(V)}$ satisfying certain conditions. Denote by $\Psi\left(G^0(V)\right)$ the set of equivalence classes of A-parameters of $G^0(V)$. Again, by composing $\psi$ with the standard representation of $\prescript{L}{}{G^0(V)}$ and restricting to $WD_E$, we can regard $\psi$ as a (conjugate) self-dual representation of $WD_E\times \SL_2(\C)$, subject to the same dimension and parity conditions as for L-parameters. Therefore, we can write $\psi$ as
\[
    \psi = \psi_{\np} + \sum_i m_i\rho_iS_{a_i}\boxtimes S_{b_i} + \widetilde{\psi}_{\np},
\]
where for each index $i$, $\rho_i$ is an irreducible (conjugate) self-dual representation of $W_E$, $a_i,b_i\in\Z_{>0}$ such that $\rho_iS_{a_i}\boxtimes S_{b_i}$ is (conjugate) self-dual of the same parity as $\psi$, $m_i\in \Z_{>0}$ is the multiplicity of $\rho_iS_{a_i}\boxtimes S_{b_i}$ in $\psi$, and any irreducible constituent of $\psi_{\np}$ is not (conjugate) self-dual of the same parity as $\psi$. For any irreducible representation $\rho S_a\boxtimes S_b$ of $WD_E\times \SL_2(\C)$ with bounded image on the Weil group, set
\[
    m_\psi\left(\rho S_a\boxtimes S_b\right) \coloneqq \dim \Hom_{WD_E\times \SL_2(\C)}\left(\rho S_a\boxtimes S_b, \psi\right).
\]We say that:
\begin{itemize}
    \item $\psi$ is of \emph{good parity} if $\psi_{\np}=0$;

    \vskip 5pt

    \item $\psi$ is \emph{tempered} if the Arthur $\SL_2(\C)$ acts trivially; in this case $\psi$ can be naturally viewed as an L-parameter by forgetting the Arthur $\SL_2(\C)$.
\end{itemize}
\vskip 5pt
For such an A-parameter $\psi$, its component group is defined to be
\[
    \calS_\psi \coloneqq \pi_0\left({\rm Cent}_{\widehat{G^0(V)}}\left({\rm Im}\left(\psi\right)\right)\right).
\]
Once again, this component group $\calS_\psi$ can be combinatorially described, following \cite[Sect. 8]{MR3202556}. Let
\[
    A_\psi \coloneqq \prod_i \left(\Z/2\Z\right) \alpha_i
\]
be the elementary $2$-group with a canonical basis $\left\{\alpha_i\right\}_i$, where $\alpha_i$ corresponds to the irreducible constituent $\rho_iS_{a_i}\boxtimes S_{b_i}$. Then we have $\calS_\psi\cong A_\psi$. The following is the main local conjecture of Arthur's endoscopic classification for classical groups.

\begin{Conj}\label{Conj:APacketsNProperties}
For each local A-parameter $\psi$ of $G^0(V)$, there is a finite set of irreducible unitary representations $\AP_\psi\left(G^0(V)\right)$, together with a canonical map (depending on the choice of a Whittaker datum $\mathscr{W}$ of a quasi-split pure inner form of $G^0(V)$)
\[
    \calJ^{\rm A}_{\mathscr{W}}: \AP_\psi\left(G^0(V)\right) \lra \Irr\left(\calS_\psi\right),
\]
characterized by the following properties:

\begin{enumerate}
    \item $\AP_\psi\left(G^0(V)\right)$ satisfies the endoscopic character relation (``ECR'' for short), as well as the twisted ECR if $G^0(V)$ is quasi-split.

    \vskip 5pt

    \item LIR: If $\psi = \psi_\tau + \psi_0 + \widetilde{\psi}_\tau$, where $\psi_\tau$ is an irreducible A-parameter of a general linear group, corresponding to an irreducible unitary generalized Speh representation $\tau$, then $\AP_\psi\left(G^0(V)\right) = \left\{\pi~\big|~\exists\,\eta_0\in\Irr\left(A_{\psi_0}\right),~\pi\subset\tau\rtimes\pi\left(\psi_0,\eta_0\right)\right\}$, and
    \[
        \tau\rtimes\pi\left(\psi_0,\eta_0\right) = \bigoplus_{\eta|_{A_{\psi_0}}=\eta_0} \pi\left(\psi,\eta\right).
    \]
    Here $V_0$ is a formed space in the Witt tower containing $V$ of appropriate dimension, $\pi\left(\psi,\eta\right)$ is the direct sum of all irreducible representations in $\AP_\psi\left(G^0(V)\right)$ mapped to $\eta\in\Irr\left(A_\psi\right)$, and likewise for $\pi\left(\psi_0,\eta_0\right)$. 
    Moreover, if $\psi_\tau$ is (conjugate) self-dual of the same parity as $\psi$, then $\eta(\alpha_\tau)$ can be determined by the eigenvalue of a certain normalized intertwining operator on $\pi\left(\psi,\eta\right)$, where $\alpha_\tau$ is the element in the canonical basis of $A_\psi$ corresponding to $\psi_\tau$.

    \vskip 5pt

    \item Let $\phi_\psi: WD_F \ra \prescript{L}{}{G^0(V)}$ be the L-parameter of $G^0(V)$ defined by the formula
    \[
        \phi_\psi\left(w,x\right) = \psi\left(w,x,\left(\begin{array}{cc}
            |w|^{\frac{1}{2}} & \\
             & |w|^{-\frac{1}{2}}
        \end{array}\right)\right),
    \]
    where $w\in W_F$ and $x\in \SL_2(\C)$. Then there is a commutative diagram
    \[
    \begin{CD}
{\LP_{\phi_\psi}\left(G^0(V)\right)} @>>> {\AP_\psi\left(G^0(V)\right)}\\
@V{\calJ^{\rm L}_{\mathscr{W}}}VV @VV{\calJ^{\rm A}_{\mathscr{W}}}V\\
{\Irr\left(A_{\phi_\psi}\right)} @>>> {\Irr\left(A_\psi\right)}
\end{CD},
    \]
    where the top horizontal arrow is an inclusion, and the bottom horizontal arrow is induced by the natural surjection $A_\psi\twoheadrightarrow A_{\phi_\psi}$.
\end{enumerate}
 \end{Conj} 

\vskip 5pt

These conjectural properties of local A-packets imply that if $\psi=\phi$ is tempered, then $\AP_\psi\left(G^0(V)\right)$ should coincide with $\LP_\phi\left(G^0(V)\right)$. Therefore one can view the LLC for tempered representations as a consequence of this conjecture. We shall call irreducible representations in $\LP_{\phi_\psi}\left(G^0(V)\right)$ ``\emph{representations of Arthur type}'', and set $\lambda_\psi\coloneqq\lambda_{\phi_\psi}$. We list the status of Conjecture \ref{Conj:APacketsNProperties} below:

\begin{itemize}
    \item For symplectic groups and quasi-split special orthogonal groups, Arthur established the existence of $\AP_\psi\left(G^0(V)\right)$ as a \emph{multi-set}, and showed all these three properties \cite[Thm. 1.5.1, 2.2.1, 2.4.1, Prop. 7.4.1]{MR3135650}. Together with M{\oe}glin's multiplicity free result \cite{MR2767522}, one knows that Conjecture \ref{Conj:APacketsNProperties} holds for these groups.

    \vskip 5pt

    \item In \cite{MR3708200}, Atobe and Gan pointed out that Arthur's result and M{\oe}glin's result implicitly imply Conjecture \ref{Conj:APacketsNProperties} for quasi-split full even orthogonal groups. 

    \vskip 5pt

    \item For quasi-split unitary groups, following Arthur's method, Mok established the existence of $\AP_\psi\left(G^0(V)\right)$ as a multi-set, and showed all these three properties \cite[Thm. 2.5.1, 3.2.1, 3.4.3, Prop. 8.4.1]{MR3338302}. It is also expected that local A-packets for unitary groups are multiplicity free. M{\oe}glin and Renard proved this for real unitary groups \cite[Thm. 1.3]{MR3947270}. For unitary groups over $p$-adic fields, to the best of our knowledge, this result has not yet appeared in the literature, although one expects that the methods used for symplectic/ orthogonal groups should adapt to this setting.

    \vskip 5pt

    \item For non-quasi-split unitary groups, Kaletha, M{\'i}nguez, Shin and White established the LLC and verified both the ECR and the LIR \cite[Thm. 1.6.1, 2.6.2]{kaletha2014endoscopic}. In particular, their results imply Conjecture \ref{Conj:APacketsNProperties} when $\psi=\phi$ is tempered.

    \vskip 5pt

    \item For non-quasi-split special odd orthogonal groups, Ishimoto established the LLC and verified both the ECR and the LIR \cite[Thm. 1.2, 4.13]{MR4776199} following \cite{kaletha2014endoscopic}. In particular, his results imply Conjecture \ref{Conj:APacketsNProperties} when $\psi=\phi$ is tempered.

    \vskip 5pt

    \item In a series of papers \cite{MR2209850, MR2533005, MR2760663, MR2767522}, M{\oe}glin developed a theory to explicitly construct the elements in local A-packets for symplectic and orthogonal groups. This theory of explicit construction of local A-packets was further developed by Bin Xu \cite{MR3679701, MR4293795}, Atobe \cite{MR4472864, MR4661536} and Hazeltine--Liu--Lo--Zhang \cite{MR4912257}. In a recent paper of Hazeltine, Baiying Liu and Chi-Heng Lo \cite{HLLEvenOrthogonal26}, they develop the explicit construction of local A-packets for even orthogonal groups (not necessarily quasi-split) based on these previous works, and in particular prove Conjecture \ref{Conj:APacketsNProperties} for non-quasi-split even orthogonal groups. They also explicate in their paper how one can adapt their method to develop this explicit construction for non-quasi-split special odd orthogonal groups, based on the result for metaplectic group recently done by Jiahe Chen \cite{chen2026construction}. 

    \vskip 5pt

    \item For those A-parameters $\psi$ such that the Deligne $\SL_2(\C)$ acts trivially (i.e. co-tempered case), it is expected that the associated A-packet $\AP_\psi\left(G^0(V)\right)$ can be constructed from tempered L-packets through the \emph{Aubert--Zelevinsky involution}. Indeed, for such a $\psi$, one can take
    \[
        \AP_\psi\left(G^0(V)\right) \coloneqq \left\{\widehat{\pi}~\big|~\pi\in\LP_{\widehat{\psi}}\left(G^0(V)\right)\right\},
    \]
    where $\widehat{\pi}$ is the Aubert--Zelevinsky dual of $\pi$, and $\widehat{\psi}$ is the tempered L-parameter obtained by switching the Deligne $\SL_2(\C)$ and Arthur $\SL_2(\C)$. We refer the readers to the recent work of Liu--Lo--Shahidi \cite[Sect. 5.3]{liu2024anti} for more details. Since the Aubert--Zelevinsky involution is compatible with the (twisted) endoscopic transfer \cite[Thm. 1.5]{MR2098406} \cite[Appendix]{MR3679701}, (twisted) ECR holds for this construction. In a forthcoming paper of the second author, Chi-Heng Lo and Brian Wen \cite{clw26}, they develop the theory of derivatives for pure inner forms of classical groups. We expect that once the theory of derivatives is established for pure inner forms, Conjecture \ref{Conj:APacketsNProperties}(2)(3) for co-tempered A-packets of pure inner forms can be verified following the methods of \cite{atobe2024local} and \cite{MR4549708}.
\end{itemize}
\vskip 5pt
To summarize, for all groups considered in this paper, Conjecture \ref{Conj:APacketsNProperties} is either established or widely expected. Throughout the rest of the paper, we assume Conjecture \ref{Conj:APacketsNProperties}.

\vskip 5pt

It will be convenient to introduce the notion of enhanced component groups. Let
\[
    \psi = \varphi + \sum_j \psi_j + \widetilde{\varphi}
\]
be either an L-parameter or A-parameter of $G^0(V)$. Here $\varphi$ is the bad parity part, and each $\psi_j$ is an irreducible representation of $WD_E$ or $WD_E\times\SL_2(\C)$ of the same parity as $\psi$, not necessarily distinct. We set
\[
    \calA_\psi \coloneqq \prod_{j} \left(\Z/2\Z\right) \alpha_j
\]
be the elementary $2$-group with a canonical basis $\left\{\alpha_j\right\}_j$, where $\alpha_j$ corresponds to the irreducible constituent $\psi_j$. Then there is a natural quotient map $\calA_\psi\twoheadrightarrow A_\psi$, and this induces an inclusion $\Irr\left(A_\psi\right)\hookrightarrow \Irr\left(\calA_\psi\right)$. For $\eta\in\Irr\left(\calA_\psi\right)$, we extend the definition
\[
    \pi\left(\psi,\eta\right) \coloneqq \begin{cases}
        \pi\left(\psi,\eta\right) \quad & \text{if }\eta\in\Irr\left(A_\psi\right);\\
        0 \quad & \text{otherwise}.
    \end{cases}
\]

\vskip 10pt


\section{Certain Jacquet modules and parabolic inductions} 
\label{sec:certain_jacquet_modules_and_parabolic_inductions}

In this section, we recall some results of Atobe on Jacquet modules \cite{MR4055180}, and then we will study certain parabolically induced representations, which will appear in the proofs of our main results.

\vskip 5pt

\emph{Caveat:} We should mention that, although Atobe only treated ${\rm Sp}_{2n}$ and split ${\rm SO}_{2n+1}$ in his paper, most of his results in \cite{MR4055180} are expected to hold for all classical groups considered in Section \ref{subsec:FormedSpaceNIsomGroup} (not necessarily quasi-split) as well. Some results in \cite{MR4055180} have already been extended to quasi-split even orthogonal and unitary groups in \cite[App. C]{atobe2024local}. 
To extend Atobe's results to all classical groups considered in Section \ref{subsec:FormedSpaceNIsomGroup}, one needs to extend two lemmas of Xu \cite[Lem. 7.2, Lem. 7.3]{MR3713922} first. As in Xu's proofs, these two lemmas can be deduced from the compatibility of Jacquet modules and (twisted) endoscopic transfer, together with some computations for general linear groups. With these two lemmas at hand, Atobe's proofs extend to all orthogonal groups and unitary groups with only mild modifications. Some of these results for pure inner forms will also be covered in a forthcoming paper \cite{clw26} of the second author and Lo and Wen. 

\subsection{Jacquet modules and local Langlands correspondence} 
\label{sub:jacquet_modules_and_local_langlands_correspondence}

Let $V$ be an $n$-dimensional formed space over $E$ as described in Section \ref{subsec:FormedSpaceNIsomGroup}, and consider the group $G^0(V)$. Suppose that $\phi$ is a good parity tempered L-parameter of $G^0(V)$. 
We would like to describe derivatives of $\pi\left(\phi,\eta\right)$ in terms of $\left(\phi,\eta\right)$ for $\eta\in\Irr\left(\calA_\phi\right)$. The first important lemma we need is the following.

\begin{Lem}[{\cite[Lem. 4.1]{MR4055180}}]\label{L:AtobeJacLLCLem.4.1}
Let $\rho\in\calC_{\unit}\left(\GL\right)$ and $x\in\R$. If $D_{\rho|\cdot|^x}\left(\pi\right)\neq 0$ for some $\pi\in\LP_\phi$, then $x\in\frac{1}{2}\Z_{\geq0}$ and
\[
    \rho S_{2x+1} \subset \phi.
\]
\end{Lem}

\vskip 5pt

This gives a necessary condition for the derivative $D_{\rho|\cdot|^x}\left(\pi\right)$ to be non-zero. The next theorem allows us to determine this derivative completely.

\begin{Thm}[{\cite[Thm. 4.2]{MR4055180}}]\label{T:AtobeJacLLCMainThm4.2}
Let $\rho\in\calC_{\unit}\left(\GL\right)$ and $x\in\frac{1}{2}\Z_{\geq 0}$. Suppose that 
\[
    \phi = \phi_0 + m \rho S_{2x+1},
\]
where $m = m_\phi\left(\rho S_{2x+1}\right)$ is the multiplicity of $\rho S_{2x+1}$ in $\phi$. Let $\eta\in\Irr\left(\calA_\phi\right)$. 
\begin{itemize}
        \item Suppose that $m=1$. If $x=0$, then $D_\rho\left(\pi\left(\phi,\eta\right)\right)=0$. If $x>0$, set
\[
    \phi' = \phi - \rho S_{2x+1} + \rho S_{2x-1}.
\]
Then there is a natural inclusion $\calA_{\phi'}\hookrightarrow \calA_{\phi}$, namely it is the identity on $\calA_{\phi_0}$, and in the case $x>\frac{1}{2}$ it sends the basis element $\alpha_{\rho,{2x-1}}\in\calA_{\phi'}$ (corresponding to $\rho S_{2x-1}$) to $\alpha_{\rho,2x+1}\in \calA_\phi$ (corresponding to $\rho S_{2x+1}$). Let $\eta'\coloneqq \eta\big|_{\calA_{\phi'}}$. We have
\[
    D_{\rho|\cdot|^x}\left(\pi\left(\phi,\eta\right)\right) = \pi\left(\phi',\eta'\right).
\]

        \vskip 5pt

        \item Suppose that $m=2$. Let $\eta_0\coloneqq \eta\big|_{\calA_{\phi_0}}$. If $x=0$, then $D_\rho\left(\pi\left(\phi,\eta\right)\right)=\pi\left(\phi_0,\eta_0\right)$. If $x>0$, set $\eta_+\coloneqq \eta$ and let $\eta_-\in\Irr\left(\calA_\phi\right)$ be the unique character of $\calA_\phi$ which coincides with $\eta_+$ on $\calA_{\phi_0}$ and takes the opposite values of $\eta_+$ on other basis elements. Let $\phi'$ be as in the previous case, $\eta'_\pm\coloneqq\eta_\pm\big|_{\calA_{\phi'}}$. We have
        \[
            D_{\rho|\cdot|^x}\left(\pi\left(\phi,\eta\right)\right) = \Delta_{\rho}\left[x,-(x-1)\right]\rtimes\pi\left(\phi_0,\eta_0\right) + \pi\left(\phi',\eta'_+\right) - \pi\left(\phi',\eta'_-\right).
        \]

        \vskip 5pt

        \item Suppose that $m\geq 3$. Take $\kappa\in\left\{1,2\right\}$ such that $m\equiv \kappa\pmod2$. Let $\phi_{m-2} \coloneqq \phi -2\rho S_{2x+1}$ and $\phi_\kappa = \phi_0 + \kappa \rho S_{2x+1}$. Then there are natural inclusions $\calA_{\phi_{m-2}} \hookrightarrow \calA_\phi$ and $\calA_{\phi_\kappa} \hookrightarrow \calA_\phi$. Set $\eta_{m-2} \coloneqq \eta\big|_{\calA_{\phi_{m-2}}}$ and $\eta_{\kappa} \coloneqq \eta\big|_{\calA_{\phi_{\kappa}}}$. We have
        \begin{align*}
            D_{\rho|\cdot|^x}\left(\pi\left(\phi,\eta\right)\right) = (m-\kappa)\Delta_{\rho}\left[x,-(x-1)\right]&\rtimes\pi\left(\phi_{m-2},\eta_{m-2}\right) \\
            &+ \Delta_\rho\left[x,-x\right]^{\frac{m-\kappa}{2}}\rtimes D_{\rho|\cdot|^x}\left(\pi\left(\phi_\kappa,\eta_\kappa\right)\right). 
        \end{align*}
    \end{itemize}    
\end{Thm}
\vskip 5pt

By using this theorem, we can determine the shape of the L-parameters for two special classes of irreducible tempered representations.

\begin{Cor}\label{C:TempShapeofL-parameter}
In the context of Theorem \ref{T:AtobeJacLLCMainThm4.2}, let $\pi\in\LP_\phi$ such that $D_{\rho|\cdot|^x}\left(\pi\right)\neq 0$. 
\begin{itemize}
       \item If $D_{\rho|\cdot|^{x'}}\left(\pi\right)=0$ for any real number $x'>x$, then the L-parameter $\phi$ is of the form
       \[
            \phi = \phi_{\rho,x} + m\rho S_{2x+1} + \rho S_{2x+3} + \cdots + \rho S_{2x+2j-1}.
       \]
       Here $\phi_{\rho,x}$ satisfies the property that for any integer $a\geq 2x+1$, $\rho S_a\not\subset \phi_{\rho,x}$; $m\coloneqq m_\phi\left(\rho S_{2x+1}\right)$, and $j\in\Z_{> 0}$.  

       \vskip 5pt

       \item If $D_{\rho|\cdot|^{x'}}\left(\pi\right)=0$ for any real number $x'\neq x$, then the L-parameter $\phi$ is of the form
       \[
            \phi = \phi^{\rho} + \rho S_{\kappa} + \rho S_{\kappa+2} + \cdots + \rho S_{\kappa + 2i-2} + m\rho S_{2x+1} + \rho S_{2x+3} + \cdots + \rho S_{2x+2j-1}.
       \]
       Here $\phi^{\rho}$ satisfies the property that for any integer $a\in\Z_{>0}$, $\rho S_a\not\subset \phi^{\rho}$; $m\coloneqq m_\phi\left(\rho S_{2x+1}\right)$, $\kappa\in\left\{1,2\right\}$ such that $\kappa\equiv 2x+1\pmod 2$, $i\in\Z_{\geq 0}$ such that $\kappa+2i-2<2x+1$, and $j\in\Z_{> 0}$. We allow $i=0$, in which case $\rho S_{\kappa} + \rho S_{\kappa+2} + \cdots + \rho S_{\kappa + 2i-2}$ is understood to be zero.
    \end{itemize}   
\end{Cor}

\vskip 5pt

\begin{proof}
We only prove the first statement here, the proof of the second is similar. We write
\[
    \phi = \phi_{\rho,x} + m_1\rho S_{2x+1} + m_2\rho S_{2x+3} + \cdots + m_j\rho S_{2x+2j-1},
\] 
where $\phi_{\rho,x}$ is as in the statement, $m_i\coloneqq m_\phi\left(\rho S_{2x+2i-1}\right)$ for all $1\leq i\leq j$, and $m_j\neq 0$. By Theorem \ref{T:LLCProperties}(1), we must have $m_i\leq 1$ for all $2\leq i\leq j$. On the other hand, if $m_i=0$ for some $i\in\left\{2,\cdots,j-1\right\}$, then by Theorem \ref{T:AtobeJacLLCMainThm4.2} we must have $D_{\rho|\cdot|^{x+i'-1}}\left(\pi\right)\neq 0$, where $i'$ is the smallest integer greater than $i$ such that $m_{i'}\neq 0$. This contradicts our assumption. Therefore $m_i=1$ for all $2\leq i\leq j$.

\end{proof}

\vskip 5pt

As another consequence, by combining Theorem \ref{T:AtobeJacLLCMainThm4.2} with Corollary \ref{leibniz-st}, we can generalize Lemma \ref{L:AtobeJacLLCLem.4.1} to certain Steinberg derivatives.

\begin{Cor}\label{C:UnitarySteinbergDerivative}
Let $\tau\cong\Delta_\rho\left[x,-x\right]$ be an irreducible unitary Steinberg representation of $\GL_d$. If $D_\tau\left(\pi\left(\phi,\eta\right)\right)\neq 0$ for some $\eta\in\Irr\left(\calA_\phi\right)$, then
\[
    m\coloneqq m_\phi\left(\rho S_{2x+1}\right) \geq 2.
\]
If this is the case, take $\kappa\in\left\{1,2\right\}$ such that $m\equiv \kappa\pmod2$. Set $\left(\phi_0,\eta_0\right)$, $\left(\phi_{m-2},\eta_{m-2}\right)$ and $\left(\phi_\kappa,\eta_\kappa\right)$ as in Theorem \ref{T:AtobeJacLLCMainThm4.2}. Then we have
\[
    D_\tau\left(\pi\left(\phi,\eta\right)\right) = \left(m-\kappa\right)\pi\left(\phi_{m-2},\eta_{m-2}\right) + \tau^{\frac{m-\kappa}{2}}\rtimes D_\tau\left(\pi\left(\phi_\kappa,\eta_\kappa\right)\right),
\]
and (when $\kappa=2$)
\[
    D_\tau\left(\pi\left(\phi_2,\eta_2\right)\right) = \pi\left(\phi_0,\eta_0\right).
\]
\end{Cor}

\vskip 5pt

\begin{proof}
We shall only prove the first statement, i.e. $m\coloneqq m_\phi\left(\rho S_{2x+1}\right) \geq 2$. Once this is established, the second statement follows easily by induction on $m$ using Theorem \ref{T:LLCProperties}(1) and Corollary \ref{leibniz-st}. 

\vskip 5pt

If $x=0$, it follows from Theorem \ref{T:AtobeJacLLCMainThm4.2} that $m\geq 2$. So assume that $x>0$. We prove by contradiction. Suppose on the contrary that $m_\phi\left(\rho S_{2x+1}\right)=1$. Note that there is an embedding
\[
    \tau \hookrightarrow \rho|\cdot|^x\times \Delta_\rho\left[x-1,-x\right].
\]
Thus $D_\tau\left(\pi\left(\phi,\eta\right)\right)\neq 0$ implies that $D_{\Delta_\rho\left[x-1,-x\right]}\circ D_{\rho|\cdot|^x}\left(\pi\left(\phi,\eta\right)\right)\neq 0$. By Theorem \ref{T:AtobeJacLLCMainThm4.2}, $D_{\rho|\cdot|^x}\left(\pi\left(\phi,\eta\right)\right)$ is irreducible and tempered; but this contradicts Casselman's criterion. Therefore we must have $m\coloneqq m_\phi\left(\rho S_{2x+1}\right) \geq 2$ as desired.

\end{proof}

\vskip 5pt

Since in later sections we will study co-tempered representations, we would like to convert these results to the co-tempered setting by using the Aubert--Zelevinsky involution. Let $\psi$ be a good parity co-tempered A-parameter of $G^0(V)$.

\begin{Cor}\label{C:co-temperedAtobeLem.4.1NStDer}
Let $\rho\in\calC_{\unit}\left(\GL\right)$ and $x\in\R$. If $D_{\rho|\cdot|^{-x}}\left(\pi\right)\neq 0$ for some $\pi\in\AP_\psi\left(G^0(V)\right)$, then $x\in\frac{1}{2}\Z_{\geq0}$ and
\[
    \rho \boxtimes S_{2x+1} \subset \psi.
\]
Moreover:
\begin{itemize}
    \item if $m_\psi\left(\rho\boxtimes S_{2x+1}\right)=1$, then $D_{\rho|\cdot|^{-x}}\left(\pi\right)$ is irreducible and lies in the co-tempered A-packet associated to the A-parameter
    \[
        \psi' \coloneqq \psi - \rho \boxtimes S_{2x+1} + \rho \boxtimes S_{2x-1};
    \]
    in this case $D_{Z_\rho[x,-x]}\left(\pi\right)=0$;

\vskip 5pt

\item if $m_\psi\left(\rho\boxtimes S_{2x+1}\right)\geq 2$, then $D_{Z_\rho[x,-x]}\left(\pi\right)\neq 0$, and any irreducible constituent of $D_{Z_\rho[x,-x]}\left(\pi\right)$ lies in the co-tempered A-packet associated to the A-parameter
\[
    \psi_{m-2} \coloneqq \psi - 2\rho\boxtimes S_{2x+1}.
\]
\end{itemize} 
\end{Cor}

\vskip 5pt

\begin{Cor}\label{C:co-temperedShapeofAPara}
Let $\pi\in\AP_\psi\left(G^0(V)\right)$ such that $D_{\rho|\cdot|^{-x}}\left(\pi\right)\neq 0$. 
\begin{itemize}
       \item If $D_{\rho|\cdot|^{-x'}}\left(\pi\right)=0$ for any real number $x'>x$, then the co-tempered A-parameter $\psi$ is of the form
       \[
            \psi = \psi_{\rho,x} + m\rho \boxtimes S_{2x+1} + \rho\boxtimes S_{2x+3} + \cdots + \rho\boxtimes S_{2x+2j-1}.
       \]
       Here $\psi_{\rho,x}$ satisfies the property that for any integer $a\geq 2x+1$, $\rho\boxtimes S_a\not\subset \psi_{\rho,x}$; $m\coloneqq m_\psi\left(\rho\boxtimes S_{2x+1}\right)$, and $j\in\Z_{> 0}$.  

       \vskip 5pt

       \item If $D_{\rho|\cdot|^{-x'}}\left(\pi\right)=0$ for any real number $x'\neq x$, then the co-tempered A-parameter $\psi$ is of the form
       \[
            \psi = \psi^{\rho} + \rho\boxtimes S_{\kappa} + \cdots + \rho\boxtimes S_{\kappa + 2i-2} + m\rho\boxtimes S_{2x+1} + \rho\boxtimes S_{2x+3} + \cdots + \rho\boxtimes S_{2x+2j-1}.
       \]
       Here $\psi^{\rho}$ satisfies the property that for any integer $a\in\Z_{>0}$, $\rho\boxtimes S_a\not\subset \psi^{\rho}$; $m\coloneqq m_\psi\left(\rho\boxtimes S_{2x+1}\right)$, $\kappa\in\left\{1,2\right\}$ such that $\kappa\equiv 2x+1\pmod 2$, $i\in\Z_{\geq 0}$ such that $\kappa+2i-2<2x+1$, and $j\in\Z_{> 0}$. We allow $i=0$, in which case $\rho\boxtimes S_{\kappa} + \cdots + \rho\boxtimes S_{\kappa + 2i-2}$ is understood to be zero.
    \end{itemize}   
\end{Cor}

\vskip 5pt


\subsection{Certain non-unitary inductions} 
\label{sub:certain_non_unitary_inductions}

Next we consider some non-unitary inductions. Let us start with a simple case. The following lemma has been established in \cite{clw26}. For the sake of completeness, we also include a sketch of proof here.

\begin{Lem}\label{L:non-unitaryInd-0}
Let $\phi$ be a good parity tempered L-parameter, $\rho\in\calC_{\unit}\left(\GL\right)$ and $x\in\frac{1}{2}\Z_{> 0}$, such that $\rho S_{2x+1}$ is (conjugate) self-dual of the same parity as $\phi$. Then for any $\eta\in\Irr\left(\calA_\phi\right)$, the co-standard module $\Delta_\rho\left[x-1,-x\right] \rtimes \pi\left(\phi,\eta\right)$ fits into the following short exact sequence
\[
    0 \lra \LS\left(\Delta_\rho\left[x-1,-x\right];\pi\left(\phi,\eta\right)\right) \lra \Delta_\rho\left[x-1,-x\right] \rtimes \pi\left(\phi,\eta\right) \lra \bigoplus_{\eta'|_{\calA_\phi}=\eta}\pi\left(\phi',\eta'\right) \lra 0,
\]
where $\phi'\coloneqq \phi + \rho S_{2x-1} + \rho S_{2x+1}$, and we view $\calA_\phi$ as a subgroup of $\calA_{\phi'}$ through the natural embedding $\calA_\phi\hookrightarrow\calA_{\phi'}$.
\end{Lem}

\vskip 5pt

\begin{proof}
In the case that $\rho S_{2x+1}\not\subset \phi$, this is proved by Atobe \cite[Prop. 5.2]{MR4055180}. So it only remains to consider the case that $m\coloneqq m_\phi\left(\rho S_{2x+1}\right)\geq 1$. We divide the proof into four steps.

\vskip 5pt

\emph{Step 1:} Firstly we prove that $\LS\left(\Delta_\rho\left[x-1,-x\right];\pi\left(\phi,\eta\right)\right)$ is the only non-tempered subquotient of $\Delta_\rho\left[x-1,-x\right] \rtimes \pi\left(\phi,\eta\right)$. Indeed, suppose that 
\[
    \sigma = \LS\left(\Delta_{\rho_1}[x_1,y_1],\cdots,\Delta_{\rho_r}[x_r,y_r];\sigma_0\right)
\]
is an irreducible non-tempered subquotient of $\Delta_\rho\left[x-1,-x\right] \rtimes \pi\left(\phi,\eta\right)$. By Theorem \ref{tadic} and Lemma \ref{L:AtobeJacLLCLem.4.1}, all exponents of Jacquet modules of $\sigma$ must be half integers. In particular, $x_1+y_1\leq -1$. Then by Corollary \ref{leibniz-st} one knows that
\begin{align*}
    D_{\Delta_{\rho_1}[x_1,y_1]}\left(\Delta_\rho\left[x-1,-x\right] \rtimes \pi\left(\phi,\eta\right)\right) = &L_{\Delta_{\rho_1}[x_1,y_1]}\left(\Delta_\rho\left[x-1,-x\right]\right)\rtimes \pi\left(\phi,\eta\right)\\
    & \qquad + R_{\Delta_{\widetilde{\rho}_1}[-y_1,-x_1]}\left(\Delta_\rho\left[x-1,-x\right]\right)\rtimes \pi\left(\phi,\eta\right)\\
    & \qquad \qquad + \Delta_\rho\left[x-1,-x\right] \rtimes D_{\Delta_{\rho_1}[x_1,y_1]}\left(\pi\left(\phi,\eta\right)\right).
\end{align*}
Note that $R_{\Delta_{\widetilde{\rho}_1}[-y_1,-x_1]}\left(\Delta_\rho\left[x-1,-x\right]\right)=0$ due to exponents, and Casselman's criterion implies that $D_{\Delta_{\rho_1}[x_1,y_1]}\left(\pi\left(\phi,\eta\right)\right)=0$. Therefore to have $D_{\Delta_{\rho_1}[x_1,y_1]}\left(\Delta_\rho\left[x-1,-x\right] \rtimes \pi\left(\phi,\eta\right)\right)\neq 0$ we must have $\Delta_\rho\left[x-1,-x\right]\cong \Delta_{\rho_1}[x_1,y_1]$, and it follows that
\begin{align*}
    D_{\Delta_{\rho_1}[x_1,y_1]}\left(\Delta_\rho\left[x-1,-x\right] \rtimes \pi\left(\phi,\eta\right)\right) &= \pi\left(\phi,\eta\right)\\
    & = D_{\Delta_{\rho_1}[x_1,y_1]}\left(\LS\left(\Delta_\rho\left[x-1,-x\right] ; \pi\left(\phi,\eta\right)\right)\right).
\end{align*}
This implies that $\sigma \cong \LS\left(\Delta_\rho\left[x-1,-x\right]; \pi\left(\phi,\eta\right)\right)$ as desired.

\vskip 5pt

\emph{Step 2:} Secondly we claim that if $\Delta_\rho\left[x-1,-x\right]\rtimes \pi\left(\phi,\eta\right)$ has an irreducible tempered subquotient $\pi'$, then $\pi'$ must belong to the L-packet associated to $\phi'$ as defined in the lemma. To see this, one just needs to note that $\pi'$ and $\LS\left(\Delta_\rho\left[x-1,-x\right];\pi\left(\phi,\eta\right)\right)$ share the same infinitesimal character, which determines a unique tempered L-parameter $\phi'$.

\vskip 5pt

\vskip 5pt

\emph{Step 3:} Thirdly we prove the case $m_\phi\left(\rho S_{2x+1}\right)=1$. Note that in this case $m_{\phi'}\left(\rho S_{2x+1}\right)=2$. Therefore any irreducible tempered subquotient $\pi'$ of $\Delta_\rho\left[x-1,-x\right]\rtimes \pi\left(\phi,\eta\right)$ satisfies $D_{\Delta_\rho[x,-x]}\left(\pi'\right)\neq 0$. On the other hand, by Theorem \ref{tadic} and a similar argument as in the proof of Corollary \ref{leibniz-st}, we have
\[
    D_{\Delta_\rho[x,-x]}\left(\Delta_\rho\left[x-1,-x\right]\rtimes \pi\left(\phi,\eta\right)\right) = D_{\rho|\cdot|^x}\left(\pi\left(\phi,\eta\right)\right).
\]
In particular: 
\begin{itemize}
    \item if $D_{\rho|\cdot|^x}\left(\pi\left(\phi,\eta\right)\right)=0$, then there is no such irreducible tempered subquotient $\pi'$; by Theorem \ref{T:AtobeJacLLCMainThm4.2} this derivative vanishes if and only if $\rho S_{2x-1}\subset \phi$ and
    \[
        \eta\left(\alpha_{\rho,2x-1}\right)\eta\left(\alpha_{\rho,2x+1}\right)=-1,
    \]
    where $\alpha_{\rho,2x-1}$ and $\alpha_{\rho,2x+1}$ are basis elements in $\calA_\phi$ corresponding to $\rho S_{2x-1}$ and $\rho S_{2x+1}$ respectively. In this subcase $\Delta_\rho\left[x-1,-x\right]\rtimes \pi\left(\phi,\eta\right)$ is irreducible.

    \vskip 5pt

    \item if $D_{\rho|\cdot|^x}\left(\pi\left(\phi,\eta\right)\right)\neq0$, then by Theorem \ref{T:AtobeJacLLCMainThm4.2} we have
    \[
        D_{\rho|\cdot|^x}\left(\pi\left(\phi,\eta\right)\right) = \pi\left(\phi_0,\eta\big|_{\calA_{\phi_0}}\right),
    \]
    where $\phi_0 = \phi - \rho S_{2x+1} + \rho S_{2x-1}$. It then follows from Corollary \ref{C:UnitarySteinbergDerivative} that there is at most one such irreducible tempered subquotient $\pi'$. On the other hand, let $\eta'\in\Irr\left(\calA_{\phi'}\right)$ be the unique character such that $\eta'\big|_{\calA_\phi} = \eta$, and takes the value $\eta\left(\alpha_{\rho,2x+1}\right)$ on other basis elements. Then, by Theorem \ref{T:LLCProperties}(1) and Theorem \ref{T:AtobeJacLLCMainThm4.2}, we know that there is an embedding
    \[
        \pi\left(\phi',\eta'\right)\hookrightarrow \rho|\cdot|^x \times \Delta_\rho\left[x-1,-x+1\right]\rtimes \pi\left(\phi,\eta\right).
    \]
    Note that $D^{(2)}_{\rho|\cdot|^x}\left(\pi\left(\phi',\eta'\right)\right)\neq 0$, and by Proposition \ref{leibniz}(2) one has
    \begin{align*}
        D^{(2)}_{\rho|\cdot|^x}\left(\rho|\cdot|^x \times \Delta_\rho\left[x-1,-x+1\right]\rtimes \pi\left(\phi,\eta\right)\right) &= \Delta_\rho\left[x-1,-x+1\right]\rtimes\pi\left(\phi_0,\eta\big|_{\calA_{\phi_0}}\right)\\
        & = D^{(2)}_{\rho|\cdot|^x}\left(\Delta_\rho\left[x,-x+1\right]\rtimes \pi\left(\phi,\eta\right)\right).
    \end{align*}
    This implies that the above embedding factors through $\Delta_\rho\left[x,-x+1\right]\rtimes \pi\left(\phi,\eta\right)$. Hence by using the MVW-contragredient trick one knows that $\pi\left(\phi',\eta'\right)$ is a quotient of $\Delta_\rho\left[x-1,-x\right]\rtimes \pi\left(\phi,\eta\right)$.
\end{itemize}
\vskip 5pt
This completes the proof of the case $m_\phi\left(\rho S_{2x+1}\right)=1$.

\vskip 5pt

\emph{Step 4:} Finally we prove the general case by using the case $m_\phi\left(\rho S_{2x+1}\right)\leq 1$. By using the same argument as in the previous step, one can easily show that all $\pi\left(\phi',\eta'\right)$ predicted by the lemma appear as quotients of the co-standard module $\Delta_\rho\left[x-1,-x\right]\rtimes\pi\left(\phi,\eta\right)$. 
Take $\kappa\in\left\{0,1\right\}$ such that $m\equiv \kappa\pmod2$. Then there is a natural inclusion
\[
    \Delta_\rho\left[x-1,-x\right]\rtimes\pi\left(\phi,\eta\right) \hookrightarrow \Delta_\rho\left[x-1,-x\right] \times \Delta_\rho[x,-x]^{\frac{m-\kappa}{2}}\rtimes\pi\left(\phi_\kappa,\eta_\kappa\right),
\]
where $\phi_\kappa \coloneqq \phi - (m-\kappa)\rho S_{2x+1}$, and $\eta_\kappa \coloneqq \eta\big|_{\calA_{\phi_\kappa}}$. When $\kappa=1$, by Theorem \ref{T:LLCProperties}(1) this inclusion is an isomorphism. When $\kappa=0$, consider the unique companion character $\eta_-\in\Irr\left(\calA_\phi\right)$ such that $\eta_-$ coincides with $\eta$ on $\calA_{\phi_0}$ and takes the opposite values of $\eta$ on other basis elements. Then
\[
    \Delta_\rho\left[x-1,-x\right]\rtimes\left(\pi\left(\phi,\eta\right)\oplus \pi\left(\phi,\eta_-\right)\right) \cong \Delta_\rho\left[x-1,-x\right] \times \Delta_\rho[x,-x]^{\frac{m-\kappa}{2}}\rtimes\pi\left(\phi_\kappa,\eta_\kappa\right).
\]
In either case, since $[x-1,-x]_\rho\subset[x,-x]_\rho$, we can exchange the induction order
\[
    \Delta_\rho\left[x-1,-x\right] \times \Delta_\rho[x,-x]^{\frac{m-\kappa}{2}}\rtimes\pi\left(\phi_\kappa,\eta_\kappa\right)\cong\Delta_\rho[x,-x]^{\frac{m-\kappa}{2}}\times \Delta_\rho\left[x-1,-x\right]\rtimes\pi\left(\phi_\kappa,\eta_\kappa\right).
\]
Applying this lemma to $\phi_\kappa$, the right-hand side of the preceding isomorphism can be decomposed as
\[
    \left[\Delta_\rho[x,-x]^{\frac{m-\kappa}{2}}\times \Delta_\rho\left[x-1,-x\right]\rtimes\pi\left(\phi_\kappa,\eta_\kappa\right)\right] = 
\]
\[
    \left[\Delta_\rho[x,-x]^{\frac{m-\kappa}{2}}\rtimes \LS\left(\Delta_\rho\left[x-1,-x\right];\pi\left(\phi_\kappa,\eta_\kappa\right)\right)\right] + \left[\Delta_\rho[x,-x]^{\frac{m-\kappa}{2}}\rtimes \left(\bigoplus_{\eta'_\kappa}\pi\left(\phi'_\kappa,\eta'_\kappa\right)\right)\right].
\]
Comparing this expression and the assertion of the lemma, it only remains to show that there is no irreducible tempered subquotient of $\Delta_\rho[x,-x]^{\frac{m-\kappa}{2}}\rtimes \LS\left(\Delta_\rho\left[x-1,-x\right];\pi\left(\phi_\kappa,\eta_\kappa\right)\right)$.

\vskip 5pt

Suppose on the contrary there is an irreducible tempered subquotient $\pi'$ of $\Delta_\rho[x,-x]^{\frac{m-\kappa}{2}}\rtimes \LS\left(\Delta_\rho\left[x-1,-x\right];\pi\left(\phi_\kappa,\eta_\kappa\right)\right)$. Note that $\LS\left(\Delta_\rho\left[x-1,-x\right];\pi\left(\phi_\kappa,\eta_\kappa\right)\right)$ lies in an Arthur type L-packet $\LP_{\phi_\psi}$, where
\[
    \psi = \phi_\kappa + \rho S_{2x}\boxtimes S_2.
\] 
In particular, it is unitary; hence the parabolically induced representation $\Delta_\rho[x,-x]^{\frac{m-\kappa}{2}}\rtimes \LS\left(\Delta_\rho\left[x-1,-x\right];\pi\left(\phi_\kappa,\eta_\kappa\right)\right)$ is semi-simple. This implies that there is an embedding
\[
    \pi' \hookrightarrow \Delta_\rho[x,-x]^{\frac{m-\kappa}{2}}\rtimes \LS\left(\Delta_\rho\left[x-1,-x\right];\pi\left(\phi_\kappa,\eta_\kappa\right)\right),
\]
and it follows that $\LS\left(\Delta_\rho\left[x-1,-x\right];\pi\left(\phi_\kappa,\eta_\kappa\right)\right)$ is an irreducible constituent of $D_{\Delta_\rho[x,-x]}^{(\frac{m-\kappa}{2})}\left(\pi'\right)$. However, Corollary \ref{C:UnitarySteinbergDerivative} asserts that any irreducible constituent of $D_{\Delta_\rho[x,-x]}^{(\frac{m-\kappa}{2})}\left(\pi'\right)$ should be tempered. This leads to a contradiction. Hence there is no such $\pi'$, and the lemma holds in the case $m_{\phi}\left(\rho S_{2x+1}\right)>1$.

\end{proof}

\vskip 5pt

For later proofs, we need to generalize this lemma to certain A-parameters. We thank Chi-Heng Lo for helpful discussions on this result.

\begin{Prop}\label{P:non-unitaryInd-1}
Suppose that
\[
    \psi = \sum_{i=1}^r \rho_i S_{a_i}\boxtimes S_{b_i} 
\]
is a good parity local A-parameter, such that
\begin{itemize}
    \item $b_i\leq 2$ for all $1\leq i\leq r$;

    \vskip 5pt

    \item if $b_i=b_j=2$ for some $i\neq j$, then $\rho_i S_{a_i}\not\cong \rho_j S_{a_j}$.
\end{itemize}
Let $\rho\in\calC_{\unit}\left(\GL\right)$ and $x\in\frac{1}{2}\Z_{>0}$, such that $\rho S_{2x}\boxtimes S_2$ is (conjugate) self-dual of the same parity as $\psi$, and $\rho S_{2x}\not\cong \rho_i S_{a_i}$ for any $1\leq i\leq r$. Then for any $\eta\in\Irr\left(\calA_{\phi_\psi}\right)$, the parabolically induced representation $\Delta_\rho\left[x-1,-x\right] \rtimes \pi\left(\phi_\psi,\eta\right)$ fits into the following short exact sequence
\[
    0 \lra \pi\left(\phi_{\psi_+},\eta\right) \lra \Delta_\rho\left[x-1,-x\right] \rtimes \pi\left(\phi_\psi,\eta\right) \lra \bigoplus_{\eta'|_{\calA_{\phi_{\psi}}}=\eta}\pi\left(\phi_{\psi'},\eta'\right) \lra 0,
\]
where $\psi_+ = \psi + \rho S_{2x}\boxtimes S_2$ and $\psi'\coloneqq \psi + \rho S_{2x-1} + \rho S_{2x+1}$; we identify $\calA_{\phi_\psi}$ with $\calA_{\phi_{\psi_+}}$, and view $\calA_{\phi_\psi}$ as a subgroup of $\calA_{\phi_{\psi'}}$ through the natural embedding $\calA_{\phi_\psi}\hookrightarrow\calA_{\phi_{\psi'}}$.
\end{Prop}

\vskip 5pt

The proof of this proposition relies on additional ingredients and is deferred to the end of this section.

\vskip 5pt

\subsection{Certain unitary inductions} 
\label{sub:certain_unitary_inductions}

We also need to study certain unitary inductions. For symplectic groups and split special odd orthogonal groups, the following result is a special case of \cite[Thm. 4.4]{MR4412277}. The main ingredients of the proof are the explicit construction of A-packets and the notion of \emph{extended multi-segments}. Thanks to the recent work \cite{HLLEvenOrthogonal26}, the proof can now be extended to even orthogonal groups. Once the explicit construction of A-packets via extended multi-segments has been extended to other classical groups, we expect the proof of \cite[Thm. 4.4]{MR4412277} will go through with only slight modifications. 

\begin{Prop}\label{P:unitaryInd}
Suppose that
\[
    \psi = \sum_{i=1}^r \rho_i S_{a_i}\boxtimes S_{b_i} 
\]
is a good parity local A-parameter, such that $b_i\leq 2$ for all $1\leq i\leq r$. Let $\rho\in\calC_{\unit}\left(\GL\right)$ and $a,b\in\Z_{>0}$, such that $b\leq 2$ and $\rho S_{a}\boxtimes S_b$ is (conjugate) self-dual of the same parity as $\psi$. Set $\psi_+ = \psi + 2\rho S_a\boxtimes S_b$, and we regard $\calA_{\phi_\psi}$ as a subgroup of $\calA_{\phi_{\psi_+}}$ through the natural embedding $\calA_{\phi_\psi}\hookrightarrow\calA_{\phi_{\psi_+}}$.
Then for any $\eta\in\Irr\left(A_{\phi_\psi}\right)$, the parabolically induced representation $u_\rho\left(a,b\right) \rtimes \pi\left(\phi_\psi,\eta\right)$ can be decomposed as follows:
\begin{itemize}
    \item If $\rho S_a\boxtimes S_b \subset \psi$, then
\[
    u_\rho\left(a,b\right) \rtimes \pi\left(\phi_\psi,\eta\right) = \pi\left(\phi_{\psi_+},\eta\right)
\]
    is irreducible. Here we note that $A_{\phi_\psi}\cong A_{\phi_{\psi_+}}$, hence $\eta$ can be regarded as a character of $A_{\phi_{\psi_+}}$.

    \vskip 5pt

    \item If $\rho S_a\boxtimes S_b \not\subset \psi$ and $b=1$, then
    \[
        u_\rho\left(a,b\right) \rtimes \pi\left(\phi_\psi,\eta\right) = \bigoplus_{\eta_+|_{\calA_{\phi_\psi}}=\eta} \pi\left(\phi_{\psi_+},\eta_+\right).
    \] 

    \vskip 5pt

    \item If $\rho S_a\boxtimes S_b \not\subset \psi$ and $b=2$, then
    \[
        u_\rho\left(a,b\right) \rtimes \pi\left(\phi_\psi,\eta\right) = \pi\left(\phi_{\psi_+},\eta\right) \oplus \bigoplus_{\eta'|_{\calA_{\phi_\psi}}=\eta} \pi\left(\phi_{\psi'},\eta'\right),
    \]
    where $\psi' = \psi + \rho S_a\boxtimes S_b + \rho S_{a+1} + \rho S_{a-1}$, and we regard $\calA_{\phi_\psi}$ as a subgroup of $\calA_{\phi_{\psi'}}$ through the natural embedding $\calA_{\phi_\psi}\hookrightarrow\calA_{\phi_{\psi'}}$.
\end{itemize}
\end{Prop}

\vskip 5pt



Finally we complete the proof of Proposition \ref{P:non-unitaryInd-1}.

\begin{proof}[Proof of Proposition \ref{P:non-unitaryInd-1}]
Let $\rho_0\coloneqq \rho$ and $a_0\coloneqq2x$. Assume that 
\[
    2= b_1=\cdots = b_k > b_{k+1} = \cdots =b_r=1.
\]
We shall prove by induction on the non-negative integer $k$. The base case is that $k=0$, i.e. $b_i=1$ for all $1\leq i\leq r$. In this case the desired conclusion follows from Lemma \ref{L:non-unitaryInd-0}. 

\vskip 5pt

Now suppose that the desired conclusion has been proved for all $k<k_0$. Similar to the argument in Step 1 and Step 2 of the proof of Lemma \ref{L:non-unitaryInd-0}, it is not hard to see that any irreducible subquotient of $\Delta_\rho\left[x-1,-x\right] \rtimes \pi\left(\phi_\psi,\eta\right)$ must be of the form $\pi\left(\phi_{\psi_I},\eta_I\right)$, where $I\subset\left\{0,1,2,\cdots,k\right\}$,
\[
    \psi_I = \sum_{i\in I} \rho_i S_{a_i}\boxtimes S_2 + \sum_{\substack{0\leq i \leq k \\ i\notin I}} \left(\rho_i S_{a_i+1} + \rho_i S_{a_i-1}\right) + \sum_{i=k+1}^r \rho_i S_{a_i},
 \]
and $\eta_I\in\Irr\left(\calA_{\phi_{\psi_I}}\right)$. We claim that $I$ must be non-empty, i.e. any irreducible subquotient of $\Delta_\rho\left[x-1,-x\right] \rtimes \pi\left(\phi_\psi,\eta\right)$ must be non-tempered. Assuming this claim, we prove the case $k=k_0$. Set
\[
    D_I \coloneqq \circ_{i\in I} D_{\Delta_{\rho_i}[x_i-1,-x_i]},
\]
where $x_i=\frac{a_i}{2}$, and the composition is taken in ascending order of the indices $i \in I$. Then by Corollary \ref{leibniz-st}, one has $D_I\left(\pi\left(\phi_{\psi_I},\eta_I\right)\right) = \pi\left(\phi_I, \eta_I\right)$, where
\[
    \phi_I = \sum_{\substack{0\leq i \leq k \\ i\notin I}} \left(\rho_i S_{a_i+1} + \rho_i S_{a_i-1}\right) + \sum_{i=k+1}^r \rho_i S_{a_i},
\]
and we identify $\calA_{\phi_I}$ with $\calA_{\phi_{\psi_I}}$. On the other hand, again by Corollary \ref{leibniz-st}, one also has 
\[
    D_I\left(\Delta_\rho\left[x-1,-x\right] \rtimes \pi\left(\phi_\psi,\eta\right)\right) = \begin{cases}
        \Delta_\rho\left[x-1,-x\right] \rtimes \pi\left(\phi_{\psi_I^-},\eta\right), \quad &\text{if }0\notin I;\\
        \pi\left(\phi_{\psi_I^-},\eta\right), \quad &\text{if }0\in I.
    \end{cases}
\]
Here
\[
    \psi_I^- = \psi - \sum_{i\in I\backslash\{0\}} \rho_i S_{a_i}\boxtimes S_{b_i}.
\]
It then follows from the induction hypothesis that $I\backslash\left\{0\right\} = \left\{1,2,\cdots,k\right\}$, and $\eta_I\big|_{\calA_{\phi_\psi}}=\eta$. The subquotient corresponding to $I=\left\{0,1,2,\cdots,k\right\}$ is $\pi\left(\phi_{\psi_+},\eta\right)$, and the others are irreducible summands of $\bigoplus_{\eta'|_{\calA_{\phi_{\psi}}}=\eta}\pi\left(\phi_{\psi'},\eta'\right)$. 

\vskip 5pt

It only remains to show our claim. Suppose that $\pi_+$ is an irreducible subquotient of $\Delta_\rho\left[x-1,-x\right] \rtimes \pi\left(\phi_\psi,\eta\right)$, then $\Delta_\rho\left[x,1-x\right]\rtimes \pi_+$ is a subquotient of
\begin{equation}\label{E:non-unitaryIndTounitaryInd}
\begin{aligned}
    \big[\Delta_\rho\left[x,1-x\right]&\rtimes \left(\Delta_\rho\left[x-1,-x\right] \rtimes \pi\left(\phi_\psi,\eta\right)\right)\big]\\
    & = \left[\left(u_\rho\left(2x,2\right) + u_\rho\left(2x+1,1\right) \times u_\rho\left(2x-1,1\right)\right)\rtimes \pi\left(\phi_\psi,\eta\right)\right].
\end{aligned}    
\end{equation}
If $\pi_+$ is tempered, then there is a unique irreducible quotient of $\Delta_\rho\left[x,1-x\right]\rtimes \pi_+$, which lies in the Arthur type L-packet associated to the A-parameter
\[
    \psi_t\coloneqq\sum_{i=0}^k \left(\rho_i S_{a_i+1} + \rho_i S_{a_i-1}\right) + \sum_{i=k+1}^r \rho_i S_{a_i} + \rho S_{2x}\boxtimes S_2.
\]
However, Proposition \ref{P:unitaryInd} implies that any irreducible subquotient of the right-hand side of (\ref{E:non-unitaryIndTounitaryInd}) should lie in one of the following Arthur type L-packets associated to A-parameters
\begin{itemize}
    \item $\psi+2\rho S_{2x}\boxtimes S_2$; or

    \vskip 5pt

    \item $\psi+\rho S_{2x}\boxtimes S_2 + \rho S_{2x+1}+ \rho S_{2x-1}$; or

    \vskip 5pt

    \item $\psi+ 2\rho S_{2x+1}+ 2\rho S_{2x-1}$.
\end{itemize}
\vskip 5pt
None of these A-packets coincides with $\psi_t$. This leads to a contradiction. Hence our claim holds.

\end{proof}

\vskip 10pt



\section{Bessel periods, some inequalities and conjectures} 
\label{sec:bessel_periods_conjectures_and_some_inequalities}

In this section, we briefly recall the definition of Bessel periods, establish some useful inequalities by doing some orbit analysis, and state the local GGP conjectures for Bessel periods, both for tempered and non-tempered L-packets. 

\subsection{Bessel periods, multiplicities and inequalities} 
\label{sub:bessel_subgroups_and_bessel_periods}

We recall the setup of \cite[Sect. 12]{MR3202556}. Fix a non-Archimedean local field $E$ of characteristic $0$ and an involution $\iota$ of $E$ with fixed field $F$. Then, we consider a pair of non-degenerate Hermitian spaces $W \subset V$ over $E$ such that
\vskip 5pt
\begin{itemize}
    \item $\dim W^{\perp}$ is odd;

    \vskip 5pt

    \item $W^{\perp}$ is a split space.
\end{itemize}
\vskip 5pt
Setting $G\left(V\right)$ (resp. $G\left(W\right)$) to be the isometry group of the space $V$ (resp. $W$) as in Section \ref{subsec:FormedSpaceNIsomGroup}, we assume further that $G\left(V\right)$ is quasi-split. As was explained in \cite[Sect. 2]{MR3202556}, one has a subgroup 
\[
    \Bes\left(V,W\right) \coloneqq G\left(W\right)^{\Delta} \cdot N \subset G\left(V\right) \times G\left(W\right)
\]
where $G\left(W\right)^{\Delta}$ is the diagonal embedding of $G\left(W\right)$ into $G\left(V\right)\times G\left(W\right)$, and $N$ is the unipotent radical of the parabolic subgroup of $G\left(V\right)$ stabilizing a full flag in $W^\perp\subset V$, which is normalized by $G\left(W\right)^{\Delta}$. Moreover, a certain character $\nu$ of $\Bes\left(V,W\right)$ was defined in \cite[Sect. 12]{MR3202556}: it is generic on $N$ and is trivial on $G\left(W\right)^\Delta$. When $\dim W^{\perp} = 1$, for example, $\nu$ is simply the trivial representation. 

\vskip 5pt

Given an irreducible representation $\pi \boxtimes \sigma$ of $G\left(V\right) \times G\left(W\right)$, Gan, Gross and Prasad defined in \cite[Sect. 12]{MR3202556} a multiplicity
 \[
    m\left(\pi,\sigma\right) \coloneqq  \dim \Hom_{\Bes(V,W)}\left(\pi\boxtimes \sigma,\nu\right)
 \]
The celebrated multiplicity one theorem \cite{MR2529937, MR2680495, MR2874638} asserts that $m\left(\pi,\sigma\right)\leq 1$. 
Gan, Gross and Prasad also introduced in \cite[Sect. 2]{MR3202556} a notion of relevant pure inner forms $G\left(V'\right) \times G\left(W'\right)$ of $G(V) \times G(W)$, and one can likewise consider the multiplicity $m\left(\pi',\sigma'\right)$ for representations of the pure inner form. Our goal is to determine the multiplicity $m\left(\pi,\sigma\right)$ when $\pi \boxtimes \sigma$ belongs to local L-packets of Arthur type. 

\vskip 5pt

To present our main results and arguments more clearly, it seems better to extend the definition of multiplicities a bit. Let $V$ and $W$ be a pair of Hermitian spaces, such that their dimensions are of different parities; and we set $G\left(V\right)$ and $G\left(W\right)$ to be their isometry groups, i.e. (full) orthogonal groups or unitary groups.
For a representation $\pi$ of $G\left(V\right)$ and $\sigma$ of $G\left(W\right)$, we set their (homological Bessel) multiplicity as follows.
\begin{itemize}
    \item If $V$ and $W$ are not relevant, namely, there is no embedding of one into the other, or there is an embedding but the orthogonal complement space is not split, then
    \[
        m^i\left(\pi,\sigma\right) \coloneqq 0;
    \]

    \vskip 5pt

    \item If there is an embedding $W\hookrightarrow V$ such that $W^\perp$ is split, then
    \[
        m^i(\pi,\sigma)\coloneqq \dim \Ext^i_{\Bes(V,W)}(\pi\boxtimes \sigma,\nu);
    \]

    \vskip 5pt

    \item If there is an embedding $V\hookrightarrow W$ such that $V^\perp$ is split, then
    \[
        m^i(\pi,\sigma)\coloneqq m^i(\sigma,\pi).
    \]
\end{itemize}
\vskip 5pt
It follows that whenever $m^i(\pi,\sigma)$ is defined, we always have $m^i(\pi,\sigma) = m^i(\sigma,\pi)$. The next result is a generalization of \cite[Thm. 15.1]{MR3202556}, which will appear many times in later proofs of our main results.

\begin{Prop}\label{P:Frobenius-MackeyInequality}
Let $V$, $W$ be a pair of Hermitian spaces defined over $E$, and let $\pi$, $\sigma$ be representations of $G\left(V\right)$, $G\left(W\right)$ respectively. Suppose that $\tau$ is an irreducible representation of $\GL_d$. Then
\[
    m\left(\pi,\D_{\left( \tau^{[0]}\right)^\vee}\left(\sigma\right)\right)\leq m\left(\tau\rtimes\pi,\sigma\right)\leq \sum_{\alpha=0}^d m\left(\pi,\D_{\left(\tau^{[\alpha]}\right)^\vee}\left(\sigma\right)\right),
\]
where $\tau^{[\alpha]}\coloneqq\tau^{(\alpha)}|\cdot|^{\frac{1}{2}}$ is the $\alpha$-th shifted Bernstein--Zelevinsky derivative of $\tau$, and when $\alpha=0$ we set $\tau^{[0]}\coloneqq \tau|\cdot|^{\frac{1}{2}}$. Moreover, let $\alpha_0$ be the highest degree such that $\tau^{[\alpha_0]}\neq 0$:
\begin{itemize}
    \item If $\D_{\left(\tau^{[\alpha]}\right)^\vee}\left(\sigma\right)=0$ for all $\alpha<\alpha_0$, then 
    \[
        m\left(\tau\rtimes\pi,\sigma\right) \leq m\left(\pi,\D_{\left(\tau^{[\alpha_0]}\right)^\vee}\left(\sigma\right)\right);
    \]

    \vskip 5pt

    \item If for $i=0,1$, one has $\D^i_{\left(\tau^{[\alpha]}\right)^\vee}\left(\sigma\right)=0$ for all $\alpha<\alpha_0$, then
    \[
        m\left(\tau\rtimes\pi,\sigma\right) = m\left(\pi,\D_{\left(\tau^{[\alpha_0]}\right)^\vee}\left(\sigma\right)\right).
    \]
\end{itemize}
\end{Prop} 

\vskip 5pt

\begin{proof}
Let $\calH$ be the (Hermitian) hyperbolic plane, and $V_d\coloneqq V \oplus \calH^d$. Without loss of generality, by using \cite[Thm. 15.1]{MR3202556}, one can assume that $W$ is a non-degenerate codimension $1$ subspace of $V_d$. Then we can appeal to the Mackey theory to estimate the multiplicity $m\left(\tau\rtimes \pi,\sigma\right)$. 

\vskip 5pt

Let $P$ be the maximal parabolic subgroup of $G\left(V_d\right)$ stabilizing a $d$-dimensional isotropic subspace of $V_d$. There is a $G\left(W\right)$-equivariant filtration on $\left(\tau\rtimes \pi\right)\big|_{G\left(W\right)}$, whose successive quotients are parametrized by the $G\left(W\right)$-orbits of $G\left(V_d\right)/P_d$. As explicated in the proof of \cite[Thm. 15.1]{MR3202556}, there are at most two such orbits:
\begin{itemize}
    \item the closed orbit consists of $d$-dimensional isotropic subspaces of $V_d$ which are contained in $W$;

    \vskip 5pt

    \item the open orbit consists of $d$-dimensional isotropic subspaces of $V_d$ which are not contained in $W$.
\end{itemize}
\vskip 5pt
Let $\Pi_{\rm cl}$ be the quotient of $\tau\rtimes \pi$ corresponding to the closed orbit, and $Q$ the stabilizer of a representative in the closed orbit. We know that $Q$ is a maximal parabolic subgroup of $G\left(W\right)$ stabilizing a $d$-dimensional isotropic subspace of $W$. For any non-negative integer $k$, let $W_{-k}$ be a subspace of $W$ such that $W\cong W_{-k}\oplus\calH^k$. Then the Levi component of $Q$ is $L\cong\GL_d\times G\left(W_{-d}\right)$. We have
\begin{align*}
    \Pi_{\rm cl} &\cong \prescript{u}{}{\ind}_Q^{G(W)}\Big(\delta_P^{\frac{1}{2}}\cdot\tau\boxtimes\left(\pi\big|_{G(W_{-d})}\right)\Big)\\
    &\cong \tau|\cdot|^{\frac{1}{2}}\rtimes\left(\pi\big|_{G(W_{-d})}\right).
\end{align*}
Here we use the symbol ``$\prescript{u}{}{\ind}$'' to denote the un-normalized compact induction. On the other hand, let $\Pi_{\rm op}$ be the $G\left(W\right)$-subrepresentation of $\tau\rtimes \pi$ corresponding to the open orbit, and $S$ the stabilizer of a representative in the open orbit. By \cite[Lem. 15.2]{MR3202556}, there is a natural projection $S\twoheadrightarrow R\times G\left(V\right)$, where $R$ is a mirabolic subgroup of $\GL_d$, and one has
\[
    \Pi_{\rm op} \cong \prescript{u}{}{\ind}_S^{G\left(W\right)} \left(\left(\tau\big|_R\right)\boxtimes \pi\right).
\]
Note that $\left(\tau\big|_R\right)\boxtimes \pi$ is a representation of $R\times G\left(V\right)$, and we regard it as a representation of $S$ by inflation. Then, Lemma \ref{L:RestrictionToMirabolic} implies that there is an $R$-equivariant filtration on $\tau\big|_R$, indexed by $1\leq\alpha\leq d$, with successive quotients $\prescript{u}{}{\ind}_{T_\alpha\rtimes N_\alpha}^R\left(|\cdot|^{\frac{\alpha}{2}}\tau^{(\alpha)}\boxtimes\psi_\alpha\right)$; this filtration induces a filtration on $\Pi_{\rm op}$, with successive quotients
\[
    \prescript{u}{}{\ind}_S^{G\left(W\right)} \left(\prescript{u}{}{\ind}_{T_\alpha\rtimes N_\alpha}^R\left(|\cdot|^{\frac{\alpha}{2}}\tau^{(\alpha)}\boxtimes\psi_\alpha\right)\boxtimes \pi\right) \cong \tau^{[\alpha]} \rtimes \left(\prescript{u}{}{\ind}_{\Bes(W_{\alpha-d},V)}^{G(W_{\alpha-d})} \left(\pi\otimes \nu_{\alpha-d}\right)\right).
\]
Here $\nu_{\alpha-d}$ is the character of $\Bes\left(W_{\alpha-d},V\right)$ defined in \cite[Sect. 12]{MR3202556}.

\vskip 5pt

Next we consider the short exact sequence of $G\left(W\right)$-representations
\[
    0 \lra \Pi_{\rm op} \lra \tau\rtimes \pi \lra \Pi_{\rm cl} \lra 0,
\]
and apply the functor $\Hom_{\Bes\left(V_d,W\right)}\left(-\boxtimes \sigma, \nu\right)$ to it. Note that by our assumption, $\Bes\left(V_d,W\right) = G\left(W\right)$ and $\nu$ is simply the trivial representation. Then, one obtains a long exact sequence
\[
    0\ra \Hom_{G\left(W\right)}\left(\Pi_{\rm cl}\boxtimes \sigma, \C\right) \ra \Hom_{G\left(W\right)}\left(\left(\tau\rtimes \pi\right)\boxtimes \sigma, \C\right) \ra \Hom_{G\left(W\right)}\left(\Pi_{\rm op}\boxtimes \sigma, \C\right) \ra \cdots
\]
By Frobenius reciprocity, the second adjointness and duality, we have
\begin{align*}
\Hom_{G\left(W\right)}\left(\Pi_{\rm cl}\boxtimes \sigma, \C\right) &\cong \Hom_{\GL_d\times G\left(W_{-d}\right)}\left(\tau|\cdot|^{\frac{1}{2}}\boxtimes \pi, \Jac_{\overline{Q}}\left(\sigma^\vee\right)\right)\\
& \cong \Hom_{G(W_{-d})}\left(\D_{\left(\tau|\cdot|^{\frac{1}{2}}\right)^\vee}\left(\sigma\right)\boxtimes \pi,\C\right),
\end{align*}
where $\overline{Q}$ is the parabolic subgroup of $G\left(W\right)$ opposite to $Q$. Likewise, for each $1\leq\alpha\leq d$, we have 
\begin{align*}
\Hom_{G\left(W\right)}\Big(\tau^{[\alpha]} &\rtimes \left(\prescript{u}{}{\ind}_{\Bes(W_{\alpha-d},V)}^{G(W_{\alpha-d})} \left(\pi\otimes \nu_{\alpha-d}\right)\right)\boxtimes \sigma, \C\Big) \\
&\cong \Hom_{\GL_{d-\alpha}\times G\left(W_{\alpha-d}\right)}\left(\tau^{[\alpha]} \boxtimes \left(\prescript{u}{}{\ind}_{\Bes(W_{\alpha-d},V)}^{G(W_{\alpha-d})} \left(\pi\otimes \nu_{\alpha-d}\right)\right), \Jac_{\overline{Q}_{\alpha-d}}\left(\sigma^\vee\right)\right)\\
& \cong \Hom_{\Bes(W_{\alpha-d},V)}\left(\D_{\left(\tau^{[\alpha]}\right)^\vee}\left(\sigma\right)\boxtimes \pi,\nu_{\alpha-d}\right),
\end{align*}
where $Q_{\alpha-d}$ is a maximal parabolic subgroup of $G\left(W\right)$ stabilizing a $(d-\alpha)$-dimensional isotropic subspace of $W$, and $\overline{Q}_{\alpha-d}$ is the opposite parabolic. The desired (in)equalities follow from the long exact sequence, and a standard spectral sequence argument.

\end{proof}

\vskip 5pt


\subsection{Local Gan--Gross--Prasad conjecture: tempered case} 
\label{sub:local_gan_gross_prasad_conjecture_tempered_case}

To state the local GGP conjecture, we first recall the distinguished characters defined in \cite[Sect. 6]{MR3202556}. Let $W\subset V$ be a pair of non-degenerate Hermitian spaces defined over $E$ as in the previous subsection. Let $\phi_V$ and $\phi_W$ be (\emph{not} necessarily tempered) L-parameters of $G^0\left(V\right)$ and $G^0\left(W\right)$ respectively. The characters $\eta_{\phi_V,\phi_W}\in\Irr\left(A_{\phi_V}\right)$ and $\eta_{\phi_W,\phi_V}\in\Irr\left(A_{\phi_W}\right)$ are defined by formulas:

\begin{itemize}
    \item if $\iota={\rm id}$, $G^0\left(V\right)$ and $G^0\left(W\right)$ are orthogonal groups, then 
    \[
        \begin{cases}
         \eta_{\phi_V,\phi_W}\left(\alpha\right) = \epsilon\left(\phi_V^\alpha\otimes\phi_W\right)\cdot \det\left(\phi_V^\alpha\right)\left(-1\right)^{\frac{\dim \phi_W}{2}}\cdot \det\left(\phi_W\right)\left(-1\right)^{\frac{\dim \phi_V^\alpha}{2}};\\
         \eta_{\phi_W,\phi_V}\left(\beta\right) = \epsilon\left(\phi_V\otimes\phi_W^\beta\right)\cdot \det\left(\phi_V\right)\left(-1\right)^{\frac{\dim \phi_W^\beta}{2}}\cdot \det\left(\phi_W^\beta\right)\left(-1\right)^{\frac{\dim \phi_V}{2}}   
        \end{cases}
    \]
    for any $\alpha\in A_{\phi_V}$ and $\beta\in A_{\phi_W}$;

    \vskip 5pt

    \item if $\iota\neq{\rm id}$, $G^0\left(V\right)$ and $G^0\left(W\right)$ are unitary groups, then 
    \[
        \begin{cases}
         \eta_{\phi_V,\phi_W}\left(\alpha\right) = \epsilon\left(\phi_V^\alpha\otimes\phi_W\right);\\
         \eta_{\phi_W,\phi_V}\left(\beta\right) = \epsilon\left(\phi_V\otimes\phi_W^\beta\right)  
        \end{cases}
    \]
    for any $\alpha\in A_{\phi_V}$ and $\beta\in A_{\phi_W}$.
\end{itemize}
\vskip 5pt
Here the notation $\epsilon\left(-\right)$ stands for the root number, we refer the readers to \cite[Sect. 5]{MR3202556} for more details; $\phi_V^\alpha$ is the $(-1)$-eigenspace of any lift of $\alpha$ in ${\rm Cent}_{\widehat{G^0(V)}}\left({\rm Im}\left(\phi_V\right)\right)$, and likewise $\phi_W^\beta$ is the $(-1)$-eigenspace of any lift of $\beta$ in ${\rm Cent}_{\widehat{G^0(W)}}\left({\rm Im}\left(\phi_W\right)\right)$. The tempered local GGP conjecture, originally proposed by Gross and Prasad \cite{MR1186476, MR1295124} for special orthogonal groups, generalized to unitary groups by Gan, Gross, and Prasad \cite{MR3202556}, and extended to full orthogonal groups in \cite[Conj. 5.2]{MR3708200} by Atobe and Gan, determines the multiplicities $m\left(\pi, \sigma\right)$ as $\pi$ and $\sigma$ vary in a tempered L-packet of relevant pure inner forms. It has been proved in full generality \cite{MR2684300, MR3155344, MR3155345, MR3202558, MR3202559, MR3155346, MR3251763, MR3676153, MR4146145, MR4568696, MR5115665, luo2020local, chenluo2022local}.

\begin{Thm}\label{T:TemperedGGP}
Let $W\subset V$ be a pair of non-degenerate Hermitian spaces defined over $E$ as in Section \ref{sub:bessel_subgroups_and_bessel_periods}, and $\phi_V$, $\phi_W$ be tempered L-parameters of $G^0\left(V\right)$, $G^0\left(W\right)$ respectively. Then we have the strong multiplicity one property
\[
    \sum_{\left(\pi,\sigma\right)\in \left(\LP_{\phi_V}\times \LP_{\phi_W}\right)^\rel} m\left(\pi,\sigma\right)=1,
\]
where:
\begin{itemize}
    \item in the case that $\iota={\rm id}$, one of $G^0(V)$ and $G^0(W)$ is a special odd orthogonal group; for any irreducible representation of this special odd orthogonal group, we view it as an irreducible representation of the full odd orthogonal group by letting $\{\pm 1\}$ act trivially;

    \vskip 5pt

    \item $\left(\LP_{\phi_V}\times \LP_{\phi_W}\right)^\rel$ is the subset of the Vogan L-packet $\LP_{\phi_V}\times \LP_{\phi_W}$ consisting of all irreducible representations of relevant pure inner forms of $G\left(V\right)\times G\left(W\right)$.
\end{itemize}
\vskip 5pt
Moreover, the unique member $\left(\pi,\sigma\right)\in\left(\LP_{\phi_V}\times \LP_{\phi_W}\right)^\rel$ such that $m\left(\pi,\sigma\right)\neq 0$ is indexed by the distinguished character $\eta_{\phi_V,\phi_W}\boxtimes\eta_{\phi_W,\phi_V}$ under the LLC.
\end{Thm}

\vskip 5pt

Note that the internal parametrization of L-packets depends on the choice of a Whittaker datum of a quasi-split pure inner form. In this theorem the Whittaker datum is chosen as explicated in \cite[Sect. 17]{MR3202556}.

\vskip 5pt


\subsection{Local Gan--Gross--Prasad conjecture: non-tempered case} 
\label{sub:local_gan_gross_prasad_conjecture_non_tempered_case}

Next we describe how to extend the local GGP conjecture to the non-tempered case. Instead of tempered L-parameters, now we consider A-parameters. The notion of ``\emph{relevant pair of A-parameters}'' \cite[Sect. 3]{MR4190046} is a key ingredient.

\begin{Def}\label{D:relevance-definition}
Let $\left(\psi_1,\psi_2\right)$ be a pair of local A-parameters, and write
\[
    \psi_1 = \sum_k \phi_{1,k} \boxtimes S_k \quad \text{and} \quad \psi_2 = \sum_\ell \phi_{2,\ell}\boxtimes S_\ell,
\]
where $\phi_{1,k}$ and  $\phi_{2,\ell}$ are tempered representations of $WD_E$. The pair $\left(\psi_1,\psi_2\right)$ is said to be \emph{relevant} if there exist decompositions
\[  
    \phi_{1,k}=\phi_{1,k}^+ + \phi_{1,k}^- \quad \text{and} \quad \phi_{2,\ell}=\phi_{2,\ell}^+ + \phi_{2,\ell}^-,
\]
such that for all $i\geq 1$, one has
\[  
    \phi_{1,i}^+=\phi_{2,i+1}^-\quad \text{and} \quad \phi_{1,i+1}^-=\phi_{2,i}^+.
\]
\end{Def}

\vskip 5pt

It is worth noting that any two tempered A-parameters are relevant, and the decompositions of $\phi_{1,k}$ and $\phi_{2,\ell}$ into $\phi_{1,k}^\pm$ and $\phi_{2,\ell}^\pm$ are unique. We give an alternative description of the notion of relevance below. For a pair of local A-parameters $\psi_1$, $\psi_2$, an irreducible representation $\varphi$ of $WD_E$ (of bounded image on $W_E$) and a positive integer $b$, set
\[
    \Lambda\left(\varphi,b;\psi_1,\psi_2\right) \coloneqq m_{\psi_1}\left(\varphi\boxtimes S_b\right) - m_{\psi_2}\left(\varphi\boxtimes S_{b+1}\right) + m_{\psi_1}\left(\varphi\boxtimes S_{b+2}\right) - m_{\psi_2}\left(\varphi\boxtimes S_{b+3}\right) + \cdots
\]

\begin{Lem}\label{L:AlterDescriptionRelevance}
Let $V$ and $W$ be a pair of Hermitian spaces, such that their dimensions are of different parities. Let $\psi_V$ and $\psi_W$ be local A-parameters of $G^0\left(V\right)$ and $G^0\left(W\right)$ respectively. Suppose that
\[
    \psi_V = \psi_{V,\np} + \psi_{V,\gp} + \widetilde{\psi}_{V,\np} \quad \text{and} \quad \psi_W = \psi_{W,\np} + \psi_{W,\gp} + \widetilde{\psi}_{W,\np},
\]
where $\psi_{V,\gp}$ (resp. $\psi_{V,\np}$) is the good parity (resp. bad parity) part of $\psi_V$, and likewise for $\psi_{W,\gp}$, $\psi_{W,\np}$. Then:
\begin{enumerate}
    \item The pair $\left(\psi_V, \psi_W\right)$ is relevant if and only if both the pairs
\[
    \left(\psi_{V,\gp}, \psi_{W,gp}\right)\quad \text{and} \quad\left(\psi_{V,\np} + \widetilde{\psi}_{V,\np}, \psi_{W,\np} + \widetilde{\psi}_{W,\np}\right)
\]
    are relevant.

    \vskip 5pt

    \item The pair $\left(\psi_V, \psi_W\right)$ is relevant if and only if these alternating sums $\Lambda\left(\varphi,b;\psi_V,\psi_W\right)$ and $\Lambda\left(\varphi,b;\psi_W,\psi_V\right)$ are non-negative for any irreducible representation $\varphi$ of $WD_E$ and positive integer $b$. The same conclusion holds for the pairs $\left(\psi_{V,\gp}, \psi_{W,gp}\right)$ and $\left(\psi_{V,\np} + \widetilde{\psi}_{V,\np}, \psi_{W,\np} + \widetilde{\psi}_{W,\np}\right)$.
\end{enumerate}
\end{Lem}

\vskip 5pt

\begin{proof}
The first statement is a direct consequence of the definition of relevance. The second statement can be proved in the same way as \cite[Lem. 3.3]{chenchenNTGGP4ArchiGL}.

\end{proof}

\vskip 5pt

Finally we can state the non-tempered local GGP conjecture \cite[Conj. 6.1, 7.1]{MR4190046} and our main result.

\begin{Conj}\label{Conj:Non-temperedGGP}
Let $W\subset V$ be a pair of non-degenerate Hermitian spaces defined over $E$ as in Section \ref{sub:bessel_subgroups_and_bessel_periods}, and $\psi_V$, $\psi_W$ be A-parameters of $G^0\left(V\right)$, $G^0\left(W\right)$ respectively. Then:
\begin{enumerate}
       \item (``Relevance implies period'') If $\left(\psi_V,\psi_W^\vee\right)$ is relevant, then the strong multiplicity one property
\[
    \sum_{\left(\pi,\sigma\right)\in \left(\LP_{\phi_{\psi_V}}\times \LP_{\phi_{\psi_W}}\right)^\rel} m\left(\pi,\sigma\right)=1
\]
        still holds, where $\left(\LP_{\phi_{\psi_V}}\times \LP_{\phi_{\psi_W}}\right)^\rel$ is the subset of $\LP_{\phi_{\psi_V}}\times \LP_{\phi_{\psi_W}}$ consisting of irreducible representations of relevant pure inner forms, and the unique member $\left(\pi,\sigma\right)\in\left(\LP_{\phi_{\psi_V}}\times \LP_{\phi_{\psi_W}}\right)^\rel$ such that $m\left(\pi,\sigma\right)\neq 0$ is indexed by the distinguished character $\eta_{\phi_{\psi_V},\phi_{\psi_W}}\boxtimes\eta_{\phi_{\psi_W},\phi_{\psi_V}}$ under the LLC.

        \vskip 5pt

        \item (``Period implies relevance'' for L-packets of Arthur type) If there exists $\left(\pi,\sigma\right)\in\left(\LP_{\phi_{\psi_V}}\times \LP_{\phi_{\psi_W}}\right)^\rel$ such that $m\left(\pi,\sigma\right)\neq 0$, then $\left(\psi_V,\psi_W^\vee\right)$ must be relevant.

        \vskip 5pt

        \item (``Period implies relevance'' for co-tempered A-packets) Suppose that $\psi_V$, $\psi_W$ are co-tempered. If there exists $\left(\pi,\sigma\right)\in\left(\AP_{{\psi_V}}\times \AP_{{\psi_W}}\right)^\rel$ such that $m\left(\pi,\sigma\right)\neq 0$, then $\left(\psi_V,\psi_W^\vee\right)$ must be relevant.

        \vskip 5pt

        \item (``Period implies relevance'' for general A-packets) If there exists $\left(\pi,\sigma\right)\in\left(\AP_{{\psi_V}}\times \AP_{{\psi_W}}\right)^\rel$ such that $m\left(\pi,\sigma\right)\neq 0$, then there exist local A-parameters $\psi_{V,\pi}$ and $\psi_{W,\sigma}$, such that $\pi\in\AP_{\psi_{V,\pi}}$, $\sigma\in\AP_{\psi_{W,\sigma}}$, and $\left(\psi_{V,\pi},\psi_{W,\sigma}^\vee\right)$ is relevant.
      \end{enumerate}   
\end{Conj}

\vskip 5pt

We note that Conjecture \ref{Conj:Non-temperedGGP}(3) implies Conjecture \ref{Conj:Non-temperedGGP}(2) when both \(\psi_V\) and \(\psi_W\) are co-tempered. With that, our main result is as follows.

\begin{Thm}\label{T:mainresult}
Parts (1) and (3) of Conjecture \ref{Conj:Non-temperedGGP} hold, whereas part (2) fails as stated in full generality.    
\end{Thm}

\vskip 5pt

The remainder of this paper is devoted to the proof of our main result.


\vskip 10pt

\section{Reduction to the good parity case} 
\label{sec:reduction_to_the_good_parity_case}

The main purpose of this section is to reduce Conjecture \ref{Conj:Non-temperedGGP} to the case that both $\psi_V$ and $\psi_W$ are of good parity. In fact we prove a slightly more general result. 

\begin{Thm}\label{T:ReductionToGoodParity}
In the setting of Conjecture \ref{Conj:Non-temperedGGP}, suppose that the A-parameters $\psi_V$ and $\psi_W$ can be decomposed as
\[
    \psi_V = \varphi_V + \psi_{V,0} + \widetilde{\varphi}_V \quad \text{and} \quad \psi_W = \varphi_W + \psi_{W,0} + \widetilde{\varphi}_{W},
\]
satisfying the following conditions:
\begin{itemize}
    \item $\psi_{V,0}\subset \psi_{V,\gp}$ and $\psi_{W,0}\subset \psi_{W,\gp}$;

    \vskip 5pt

    \item if $\rho S_a\boxtimes S_b\subset \varphi_{V}$, then $\rho|\cdot|^x\not\subset\lambda_{\psi_{V,0}}$ and $\rho^\vee|\cdot|^{x+\frac{1}{2}} \not\subset \lambda_{\psi_{W,0}}$ for any $x\in\frac{a+b}{2}+\Z$;

    \vskip 5pt

    \item likewise, if $\rho S_a\boxtimes S_b\subset \varphi_{W}$, then $\rho|\cdot|^x\not\subset\lambda_{\psi_{W,0}}$ and $\rho^\vee|\cdot|^{x+\frac{1}{2}} \not\subset \lambda_{\psi_{V,0}}$ for any $x\in\frac{a+b}{2}+\Z$.
\end{itemize}
\vskip 5pt
Let $\tau_{\varphi_V}$ be the irreducible unitary representation of a general linear group corresponding to the A-parameter $\varphi_V$, and likewise $\tau_{\varphi_W}$ the irreducible unitary representation corresponding to $\varphi_W$.
\begin{enumerate}
    \item If $\left(\varphi_V + \widetilde{\varphi}_V,\left(\varphi_W + \widetilde{\varphi}_W\right)^\vee\right)$ is relevant, and $m\left(\pi_0,\sigma_0\right)\neq0$ for a pair of irreducible representations $\left(\pi_0,\sigma_0\right)\in \left(\AP_{\psi_{V,0}}\times\AP_{\psi_{W,0}}\right)^\rel$, then we have
    \[
        m\left(\tau_{\varphi_V}\rtimes\pi_0, \tau_{\varphi_W}\rtimes \sigma_0\right) \neq 0.
    \] 

    \vskip 5pt

    \item Let $\left(\pi_0,\sigma_0\right)\in \left(\AP_{\psi_{V,0}}\times\AP_{\psi_{W,0}}\right)^\rel$, and let $\left(\phi_{\pi_0},\eta_{\pi_0}\right)$ and $\left(\phi_{\sigma_0},\eta_{\sigma_0}\right)$ be the corresponding enhanced L-parameters. If there exists a pair of irreducible representations $\pi\subset\tau_{\varphi_V}\rtimes\pi_0$ and $\sigma\subset\tau_{\varphi_W}\rtimes\sigma_0$, such that
    \vskip 5pt
    \begin{itemize}
        \item $m\left(\pi,\sigma\right)\neq 0$;

        \vskip 5pt

        \item the enhanced L-parameter $\left(\phi_\pi,\eta_\pi\right)$ of $\pi$ satisfies $\phi_\pi = \phi_{\varphi_V}+\phi_{\pi_0}+\phi_{\widetilde{\varphi}_V}$, and $\eta_\pi\big|_{\calA_{\phi_{\pi_0}}}=\eta_{\pi_0}$;

        \vskip 5pt

        \item the enhanced L-parameter $\left(\phi_\sigma,\eta_\sigma\right)$ of $\sigma$ satisfies $\phi_\sigma = \phi_{\varphi_W}+\phi_{\sigma_0}+\phi_{\widetilde{\varphi}_W}$, and $\eta_\sigma\big|_{\calA_{\phi_{\sigma_0}}}=\eta_{\sigma_0}$;
    \end{itemize}
    then $\left(\varphi_V + \widetilde{\varphi}_V,\left(\varphi_W + \widetilde{\varphi}_W\right)^\vee\right)$ is relevant, and $m\left(\pi_0,\sigma_0\right)\neq0$.
\end{enumerate}
\end{Thm}

\vskip 5pt

Note that for a pair of Hermitian spaces $W\subset V$ as in the setting of Conjecture \ref{Conj:Non-temperedGGP}, given any pair of A-parameters $\left(\psi_V,\psi_W\right)$, the decompositions 
\[
    \psi_V = \psi_{V,\np} + \psi_{V,\gp} + \widetilde{\psi}_{V,\np} \quad \text{and} \quad \psi_W = \psi_{W,\np} + \psi_{W,\gp} + \widetilde{\psi}_{W,\np} 
\]
satisfy all conditions in the above theorem. As a direct consequence, we deduce that:
\begin{Cor}
If Conjecture \ref{Conj:Non-temperedGGP}(1) (resp. (3), (4)) holds for good parity (resp. good parity co-tempered, good parity) A-parameters, then it holds in full generality.  
\end{Cor}

\vskip 5pt

We can also deduce a partial result on the direction ``period implies relevance'' beyond the co-tempered case.

\begin{Cor}\label{C:PeriodImpliesRelevanceEvenParameter}
In the setting of Conjecture \ref{Conj:Non-temperedGGP}, suppose that the A-parameters $\psi_V$ and $\psi_W$ can be written as
\[
    \psi_V = \varphi_V + \widetilde{\varphi}_V \quad \text{and} \quad \psi_W = \varphi_W + \widetilde{\varphi}_{W},
\] 
where $\varphi_V$ and $\varphi_W$ are representations of $WD_E\times \SL_2(\C)$. If there exists $\left(\pi,\sigma\right)\in\left(\LP_{\phi_{\psi_V}}\times \LP_{\phi_{\psi_W}}\right)^\rel$ such that $m\left(\pi,\sigma\right)\neq 0$, then $\left(\psi_{V},\psi_{W}^\vee\right)$ is relevant.
\end{Cor}

\subsection{Ladder representations and their Jacquet modules} 
\label{sub:ladder_representations_and_their_jacquet_modules}

In the proof of Theorem \ref{T:ReductionToGoodParity}, we need to compute some Jacquet modules of $\tau_{\varphi_V}$ and $\tau_{\varphi_W}$. Since $\tau_{\varphi_V}$ and $\tau_{\varphi_W}$ are products of generalized Speh representations $u_\rho\left(a,b\right)$, we are interested in determining $m^*\left(u_\rho\left(a,b\right)\right)$, where $m^*$ is the co-multiplication of $\calR$ defined in Section \ref{subsec:TadicFormulaDerivative}.

\vskip 5pt

This problem was settled by Kret and Lapid \cite[Thm. 2.1]{MR2996769} for a larger class of irreducible representations called ``\emph{ladder representations}''. Let
\[
    \tau=\LS\left(\Delta_{\rho_1}[x_1,y_1],\Delta_{\rho_2}[x_2,y_2],\cdots,\Delta_{\rho_r}[x_r,y_r]\right)
\]
be an irreducible representation as in Theorem \ref{T:LanglandsClass4GL}. We shall say that $\tau$ is a ladder representation if 
\begin{itemize}
    \item $\rho_1\cong\rho_2\cong\cdots\cong\rho_r$;

    \vskip 5pt

    \item $x_1\equiv x_2\equiv\cdots\equiv x_r \pmod \Z$ and $y_1\equiv y_2\equiv\cdots\equiv y_r \pmod \Z$;

    \vskip 5pt

    \item $x_1<x_2<\cdots< x_r$ and $y_1<y_2<\cdots< y_r$.
\end{itemize}
\vskip 5pt
For such representations, their Jacquet modules satisfy a very elegant formula.

\begin{Lem}[{\cite[Thm. 2.1]{MR2996769}}]\label{L:KretLapidJacFormula}
Let $\tau$ be a ladder representation as above. Then
\begin{align*}
    m^*\left(\tau\right) = \sum_{z_1<z_2<\cdots<z_r} \LS\big(\Delta_{\rho_1}[x_1,z_1],&\Delta_{\rho_2}[x_2,z_2],\cdots,\Delta_{\rho_r}[x_r,z_r]\big)\\ 
    &\otimes \LS\big(\Delta_{\rho_1}[z_1-1,y_1],\Delta_{\rho_2}[z_2-1,y_2],\cdots,\Delta_{\rho_r}[z_r-1,y_r]\big).
\end{align*}  
\end{Lem}

\vskip 5pt

It is worth noting that on the right-hand side of the formula, $\LS\left(\Delta_{\rho_1}[x_1,z_1],\cdots,\Delta_{\rho_r}[x_r,z_r]\right)$ and $\LS\left(\Delta_{\rho_1}[z_1-1,y_1],\cdots,\Delta_{\rho_r}[z_r-1,y_r]\right)$ are also ladder representations. 

\vskip 5pt


\subsection{Proof of Theorem \ref{T:ReductionToGoodParity}(1): periods of unitary inductions} 
\label{sub:proof_of_theorem_t:reductiontogoodparity_the_}

In this subsection we prove the first part of Theorem \ref{T:ReductionToGoodParity}, namely, we show that if $\left(\varphi_V + \widetilde{\varphi}_V,\left(\varphi_W + \widetilde{\varphi}_W\right)^\vee\right)$ is relevant and $m\left(\pi_0,\sigma_0\right)\neq0$, then $m\left(\tau_{\varphi_V}\rtimes\pi_0, \tau_{\varphi_W}\rtimes \sigma_0\right)\neq 0$.

\vskip 5pt

Firstly we claim that, without loss of generality, one can assume that $\left(\varphi_V,\varphi_W^\vee\right)$ is relevant. To see this, we decompose $\varphi_V + \widetilde{\varphi}_V$ and $\left(\varphi_W + \widetilde{\varphi}_W\right)^\vee$ as
\[
    \varphi_V + \widetilde{\varphi}_V = \sum_{k} \left(\Phi_{V,k}^+ + \Phi_{V,k}^-\right) \boxtimes S_k \quad \text{and} \quad \left(\varphi_W + \widetilde{\varphi}_W\right)^\vee = \sum_{\ell} \left(\Phi_{W,\ell}^+ + \Phi_{W,\ell}^-\right) \boxtimes S_\ell
\]
as in Definition \ref{D:relevance-definition}, so that for all $i\geq 1$, one has $\Phi_{V,i}^+ = \Phi_{W,i+1}^-$ and $\Phi_{V,i+1}^- = \Phi_{W,i}^+$. Since $\varphi_V + \widetilde{\varphi}_V$ and $\left(\varphi_W + \widetilde{\varphi}_W\right)^\vee$ are (conjugate) self-dual and the above pair of decompositions is unique, we have
\[
    \begin{cases}
    \Phi_{V,k}^+ = \widetilde{\Phi}_{V,k}^+,\\ 
     \Phi_{V,k}^- = \widetilde{\Phi}_{V,k}^-   
    \end{cases}
    \quad \text{and} \quad
        \begin{cases}
    \Phi_{W,k}^+ = \widetilde{\Phi}_{W,k}^+,\\ 
     \Phi_{W,k}^- = \widetilde{\Phi}_{W,k}^-.   
    \end{cases}
\]
Therefore we can decompose $\Phi_{V,k}^+$ as
\[
    \Phi_{V,k}^+ = \Big(\sum_{i\in I_{V,k}^+} \rho_i S_{a_i}\Big) + \Big(\sum_{j\in J_{V,k}^+} m_j\rho_j S_{a_j}\Big) + \Big(\sum_{i\in I_{V,k}^+} \widetilde{\rho}_i S_{a_i}\Big),
\]
where for each index $i\in I_{V,k}^+$, $\rho_i S_{a_i}$ is not (conjugate) self-dual; and for each index $j\in J_{V,k}^+$, $\rho_j S_{a_j}$ is (conjugate) self-dual, and $m_j \coloneqq m_{\Phi_{V,k}^+}\left(\rho_j S_{a_j}\right)$ is the multiplicity of $\rho_j S_{a_j}$ in $\Phi_{V,k}^+$. Likewise, we also have the decomposition
\[
    \Phi_{V,k}^- = \Big(\sum_{i\in I_{V,k}^-} \rho_i S_{a_i}\Big) + \Big(\sum_{j\in J_{V,k}^-} m_j\rho_j S_{a_j}\Big) + \Big(\sum_{i\in I_{V,k}^-} \widetilde{\rho}_i S_{a_i}\Big).
\]
We would like to show that for any $j\in J_{V,k}^\pm$, $m_j$ is even. Note that by the (conjugate) self-duality of $\rho_j S_{a_j}$, the multiplicities $m_{\varphi_V + \widetilde{\varphi}_V}\left(\rho_j S_{a_j}\boxtimes S_{k+2t}\right)$ and $m_{(\varphi_W + \widetilde{\varphi}_W)^\vee}\left(\rho_j S_{a_j}\boxtimes S_{k+2t+1}\right)$ are even for any $t\in\Z_{\geq 0}$. It follows that for $j\in J_{V,k}^-$, one knows that
\[
    m_j = \Lambda\left(\rho_j S_{a_j},k; \varphi_V + \widetilde{\varphi}_V,\left(\varphi_W + \widetilde{\varphi}_W\right)^\vee\right)
\]
is even. For $j\in J_{V,k}^+$, $m_j = m_{\varphi_V + \widetilde{\varphi}_V}\left(\rho_j S_{a_j}\right) - m_{\Phi_{V,k}^-}\left(\rho_j S_{a_j}\right)$ is also even. Hence, we can set
\[
    \left(\varphi'_{V}\right)_k^\pm \coloneqq \Big(\sum_{i\in I_{V,k}^\pm} \rho_i S_{a_i}\Big) + \Big(\sum_{j\in J_{V,k}^\pm} \frac{m_j}{2}\rho_j S_{a_j}\Big),
\]
and
\[
    \varphi'_V \coloneqq \sum_k \left(\left(\varphi'_{V}\right)_k^+ + \left(\varphi'_{V}\right)_k^-\right)\boxtimes S_k.
\]
It is easy to check that $\varphi'_V + \widetilde{\varphi}'_V = \varphi_V + \widetilde{\varphi}_V$. We define $\varphi'_W$ similarly. Then it follows from our construction that $\left(\varphi'_V, \left(\varphi'_W\right)^\vee\right)$ is relevant. By replacing $\varphi_V$ and $\varphi_W$ by $\varphi'_V$ and $\varphi'_W$ if necessary, we can assume $\left(\varphi_V, \left(\varphi_W\right)^\vee\right)$ is relevant.

\vskip 5pt

Next, under this stronger relevance assumption, we prove that $m\left(\tau_{\varphi_V}\rtimes\pi_0, \tau_{\varphi_W}\rtimes \sigma_0\right)\neq 0$ by induction. Since $\left(\varphi_V, \left(\varphi_W\right)^\vee\right)$ is relevant, we can write these two A-parameters as
\begin{align*}
    \varphi_V &= \left(\varphi^+_{V,1} + \varphi^+_{V,2} + \cdots + \varphi^+_{V,r}\right) + \left(\varphi^-_{V,1} + \varphi^-_{V,2} + \cdots + \varphi^-_{V,s}\right); \text{ and}\\
    \varphi^\vee_W &= \left(\varphi^+_{W,1} + \varphi^+_{W,2} + \cdots + \varphi^+_{W,s}\right)^\vee + \left(\varphi^-_{W,1} + \varphi^-_{W,2} + \cdots + \varphi^-_{W,r}\right)^\vee,
\end{align*}
where:
\begin{itemize}
    \item for each $1\leq i\leq r$, $\varphi^+_{V,i} = \rho_i S_{a_i}\boxtimes S_{b_i}$ is an irreducible representation of $WD_E\times \SL_2(\C)$, $\varphi^-_{W,i} = \rho^\vee_i S_{a_i}\boxtimes S_{b_i+1}$, and 
    \[
        a_1+b_1\leq a_2+b_2\leq \cdots\leq a_r+b_r;
    \]

    \vskip 5pt

    \item for each $1\leq j\leq s$, $\left(\varphi^+_{W,j}\right)^\vee = \rho'_j S_{c_j}\boxtimes S_{d_j}$ is an irreducible representation of $WD_E\times \SL_2(\C)$, $\varphi^-_{V,j} = \rho'_j S_{c_j}\boxtimes S_{d_j+1}$, and 
    \[
        c_1+d_1\leq c_2+d_2\leq \cdots\leq c_s+d_s.
    \]
\end{itemize}
\vskip 5pt
Here, we allow $b_i$ and $d_j$ to be zero; in which case $\rho_i S_{a_i}\boxtimes S_{0}$ and $\rho'_j S_{c_j}\boxtimes S_{0}$ are treated as zero. We re-order the set
\[
    \left\{\varphi^+_{V,1}, \varphi^+_{V,2}, \cdots, \varphi^+_{V,r}, \left(\varphi^+_{W,1}\right)^\vee, \left(\varphi^+_{W,2}\right)^\vee, \cdots, \left(\varphi^+_{W,s}\right)^\vee\right\}
\]
as $\left\{\phi_1,\phi_2,\cdots,\phi_{r+s}\right\}$, where for each $1\leq i\leq r+s$, $\phi_i = \varrho_i S_{p_i}\boxtimes S_{q_i}$ and
\[
    p_1+q_1\leq p_2+q_2 \leq \cdots\leq p_{r+s}+q_{r+s}.
\]
We set
\[
    \phi_{V,i} \coloneqq \begin{cases}
      \varphi^+_{V,i'}, \quad &\text{if }\phi_i= \varphi^+_{V,i'};\\
      \varphi^-_{V,j'}, \quad &\text{if }\phi_i= \left(\varphi^+_{W,j'}\right)^\vee  
    \end{cases}
    \quad\text{and}\quad
    \phi_{W,i} \coloneqq \begin{cases}
      \varphi^-_{W,i'}, \quad &\text{if }\phi_i= \varphi^+_{V,i'};\\
      \varphi^+_{W,j'}, \quad &\text{if }\phi_i= \left(\varphi^+_{W,j'}\right)^\vee  
    \end{cases},
\]
and set $\tau_{V,i}$, $\tau_{W,i}$ to be irreducible unitary representations corresponding to A-parameters $\phi_{V,i}$, $\phi_{W,i}$ respectively. Let
\[
    \pi_i\coloneqq \tau_{V,i}\rtimes \pi_{i-1} \quad \text{and} \quad \sigma_i\coloneqq \tau_{W,i}\rtimes \sigma_{i-1}.
\]
Then one has $\pi_{r+s} = \tau_{\varphi_V}\rtimes\pi_0$ and $\sigma_{r+s}=\tau_{\varphi_W}\rtimes \sigma_0$. We shall perform induction on this positive integer $i$. Namely, given that $m\left(\pi_{i-1},\sigma_{i-1}\right)\neq 0$, we aim to show $m\left(\pi_{i},\sigma_{i}\right)\neq 0$.

\vskip 5pt

In the case that $\phi_{W,i} = \phi_i^\vee$, i.e. $\tau_{V,i}=u_{\varrho}\left(p_i,q_i+1\right)$ and $\tau_{W,i}=u_{\varrho^\vee}\left(p_i,q_i\right)$, we apply Proposition \ref{P:Frobenius-MackeyInequality} to compute $m\left(\tau_{V,i}\rtimes \pi_{i-1},\sigma_i\right)$. Again by \cite[Thm. 15.1]{MR3202556} there is no harm to assume that the pair $\left(\pi_i,\sigma_i\right)$ is in the codimension $1$ situation. Note that for any $\alpha<p_i\cdot \dim \varrho$, $\varrho|\cdot|^{\frac{p_i+q_i}{2}}$ is in the cuspidal support of $\tau_{V,i}^{[\alpha]}$, whereas $\left(\varrho|\cdot|^{\frac{p_i+q_i}{2}}\right)^\vee$ is not in the cuspidal support of $\sigma_i$. This implies that $\D^*_{\left(\tau_{V,i}^{[\alpha]}\right)^\vee}\left(\sigma_i\right) = 0$ for any degree $*\in\Z_{\geq0}$. Therefore
\[
    m\left(\pi_i,\sigma_i\right) = \dim \Hom\left(u_{\varrho}\left(p_i,q_i\right)\rtimes{\prescript{u}{}{\ind}_{\Bes}\left(\pi_{i-1}\right)}, u_{\varrho}\left(p_i,q_i\right)\rtimes \sigma_{i-1}^\vee\right),
\]
which is non-zero since by using the parabolic induction functor $u_{\varrho}\left(p_i,q_i\right)\rtimes-$ one knows that the right-hand side of the preceding equality is at least $m\left(\pi_{i-1},\sigma_{i-1}\right)$.

\vskip 5pt

In the case that $\phi_{V,i} = \phi_i$, i.e. $\tau_{V,i}=u_{\varrho}\left(p_i,q_i\right)$ and $\tau_{W,i}=u_{\varrho^\vee}\left(p_i,q_i+1\right)$, the desired conclusion can be proved similarly by applying Proposition \ref{P:Frobenius-MackeyInequality} to $m\left(\tau_{W,i}\rtimes \sigma_{i-1},\pi_i\right)$. This completes the proof.


\subsection{Proof of Theorem \ref{T:ReductionToGoodParity}(2): partial relevance of A-parameters} 
\label{sub:proof_of_theorem_t:reductiontogoodparity_the_only_if_part}

In this subsection we prove the second part of Theorem \ref{T:ReductionToGoodParity}. More precisely, we show that if $m\left(\pi, \sigma\right)\neq 0$ for a pair of representations in $\left(\LP_{\phi_{\psi_{V}}}\times\LP_{\phi_{\psi_{W}}}\right)^\rel$, then the pair of A-parameters $\left(\varphi_V + \widetilde{\varphi}_V,\left(\varphi_W + \widetilde{\varphi}_W\right)^\vee\right)$ is relevant and $m\left(\pi_0,\sigma_0\right)\neq0$. Our argument in this subsection is largely inspired by an email correspondence of the first author and Maxim Gurevich on an analogous problem for general linear groups. We thank Maxim for helpful discussions.

\vskip 5pt

Suppose that
\[
    \varphi_V = \sum_{i=1}^r \rho_i S_{a_i}\boxtimes S_{b_i} \quad \text{and} \quad \varphi_W = \sum_{j=1}^s \rho'_j S_{c_j}\boxtimes S_{d_j},
\]
and these sums are ordered so that $a_1+b_1\geq a_2+b_2\geq \cdots\geq a_r+b_r$ and $c_1+d_1\geq c_2+d_2\geq \cdots\geq c_s+d_s$. For simplicity of notation we set 
\[
    \begin{cases}
       A_i \coloneqq \frac{a_i+b_i}{2}-1,\\
       B_i \coloneqq \frac{a_i-b_i}{2} 
    \end{cases}
    \quad\text{and}\quad
    \begin{cases}
      C_j \coloneqq \frac{c_j+d_j}{2}-1,\\
      D_j \coloneqq \frac{c_j-d_j}{2}.  
    \end{cases}
\]
We shall prove the desired conclusion by induction on $\dim \varphi_V + \dim \varphi_W$. The base case is when $\dim \varphi_V + \dim \varphi_W=0$, in which case the desired conclusion holds automatically.
Next assume that $\dim \varphi_V + \dim \varphi_W>0$. Since $m\left(\pi,\sigma\right)=m\left(\sigma,\pi\right)$, without loss of generality we may assume that $a_1+b_1\geq c_1+d_1$. It would be convenient to treat the cases $b_1>1$ and $b_1=1$ separately. 

\vskip 5pt

\emph{Case 1:} $b_1>1$. Let $\varphi_V^-\coloneqq \varphi_V - \rho_1 S_{a_1}\boxtimes S_{b_1} + \rho_1 S_{a_1}\boxtimes S_{b_1-2}$, and $\psi_V^- \coloneqq \varphi_V^- + \psi_{V,0} + \widetilde{\varphi}_V^-$. One can easily check that the pair of A-parameters $\left(\psi_V^-,\psi_W\right)$ still satisfies the conditions in Theorem \ref{T:ReductionToGoodParity}. Let $\pi^-\coloneqq\pi\Big(\phi_\pi^-,\eta_\pi\big|_{\calA_{\phi_\pi^-}}\Big)$, where $\phi_\pi^- \coloneqq \phi_{\varphi_V^-} + \phi_{\pi_0} + \phi_{\widetilde{\varphi}_V^-}$. Thanks to \cite[Prop. 5.3]{MR2504024}, we have $\pi^-\subset\tau_{\varphi_V^-}\rtimes\pi_0$. By Theorem \ref{T:LLCProperties}(2), one knows that there is a surjection
\[
    \Delta_{\rho_1}\left[A_1,-B_1\right] \times \Delta_{\widetilde{\rho}_1}\left[A_1,-B_1\right] \rtimes \pi^-\twoheadrightarrow \pi.
\]
Therefore we have
\[
    m\left(\Delta_{\rho_1}\left[A_1,-B_1\right] \times \Delta_{\widetilde{\rho}_1}\left[A_1,-B_1\right] \rtimes \pi^-,\sigma\right) \geq m\left(\pi,\sigma\right) > 0.
\]
Note that for any $1\leq j\leq s$, $A_1+\frac{1}{2}>C_j$. Hence by a standard cuspidal support argument, one can show that $\D_{\left(\Delta_{\rho_1}\left[A_1,-B_1\right] \times \Delta_{\widetilde{\rho}_1}\left[A_1,-B_1\right]\right)^{[\alpha],\vee}}\left(\sigma\right) = 0$
for any $\alpha<2d_{\rho_1}a_1$, where $d_{\rho_1}\coloneqq \dim \rho_1$. It then follows from Proposition \ref{P:Frobenius-MackeyInequality} that
\[
    m\left(\Delta_{\rho_1}\left[A_1,-B_1\right] \times \Delta_{\widetilde{\rho}_1}\left[A_1,-B_1\right] \rtimes \pi^-,\sigma\right) \leq m\left(\pi^-,\sigma\right).
\]
Since $\dim\varphi_V^-+\dim\varphi_W<\dim\varphi_V+\dim\varphi_W$, by the induction hypothesis, the pair of A-parameters $\left(\varphi_V^-+\widetilde{\varphi}_V^-, \left(\varphi_W+\widetilde{\varphi}_W\right)^\vee\right)$ is relevant and $m\left(\pi_0,\sigma_0\right)\neq 0$.

\vskip 5pt

We then appeal to Lemma \ref{L:AlterDescriptionRelevance} to show the relevance of the pair $\left(\varphi_V+\widetilde{\varphi}_V, \left(\varphi_W+\widetilde{\varphi}_W\right)^\vee\right)$. Note that for any irreducible representation $\varphi$ of $WD_E$ and any $b\in\Z_{>0}$, we have
\[
    \Lambda\left(\varphi,b;\varphi_V+\widetilde{\varphi}_V, \left(\varphi_W+\widetilde{\varphi}_W\right)^\vee\right) \geq \Lambda\left(\varphi,b;\varphi_V^-+\widetilde{\varphi}_V^-, \left(\varphi_W+\widetilde{\varphi}_W\right)^\vee\right) \geq 0.
\]
On the other hand, for any $\left(\varphi,b\right)\notin \left\{\left(\rho_1 S_{a_1}, b_1-1\right), \left(\widetilde{\rho}_1 S_{a_1}, b_1-1\right)\right\}$, the alternating sums
\[
    \Lambda\left(\varphi,b;\left(\varphi_W+\widetilde{\varphi}_W\right)^\vee,\varphi_V+\widetilde{\varphi}_V\right) = \Lambda\left(\varphi,b;\left(\varphi_W+\widetilde{\varphi}_W\right)^\vee,\varphi_V^-+\widetilde{\varphi}_V^-\right)
\]
are non-negative; for the case that $\left(\varphi,b\right) = \left(\rho_1 S_{a_1}, b_1-1\right)$, by definition
\[
    \Lambda\left(\rho_1 S_{a_1}, b_1-1;\left(\varphi_W+\widetilde{\varphi}_W\right)^\vee,\varphi_V+\widetilde{\varphi}_V\right) = m_{\left(\varphi_W+\widetilde{\varphi}_W\right)^\vee}\left(\rho_1 S_{a_1}\boxtimes S_{b_1-1}\right) - m_{\varphi_V+\widetilde{\varphi}_V}\left(\rho_1 S_{a_1}\boxtimes S_{b_1}\right);
\]
and for the case that $\left(\varphi,b\right) = \left(\widetilde{\rho}_1 S_{a_1}, b_1-1\right)$, by the (conjugate) self-duality of $\left(\varphi_W+\widetilde{\varphi}_W\right)^\vee$ and $\varphi_V+\widetilde{\varphi}_V$ we have
\[
    \Lambda\left(\widetilde{\rho}_1 S_{a_1}, b_1-1;\left(\varphi_W+\widetilde{\varphi}_W\right)^\vee,\varphi_V+\widetilde{\varphi}_V\right) = \Lambda\left(\rho_1 S_{a_1}, b_1-1;\left(\varphi_W+\widetilde{\varphi}_W\right)^\vee,\varphi_V+\widetilde{\varphi}_V\right).
\]
Thus it only remains to show that
\begin{equation}\label{E:ReductionToGoodParityOnlyIfInequality-1}
m_{\left(\varphi_W+\widetilde{\varphi}_W\right)^\vee}\left(\rho_1 S_{a_1}\boxtimes S_{b_1-1}\right) \geq m_{\varphi_V+\widetilde{\varphi}_V}\left(\rho_1 S_{a_1}\boxtimes S_{b_1}\right).   
\end{equation}
Consider the positive integer $\beta\coloneqq\max_{i,j}\left\{b_i,d_j\right\}$. There are three subcases.

\vskip 5pt

\emph{Subcase 1.1:} $b_1=\beta$, and $d_j<\beta$ for any $1\leq j\leq s$. We set
\begin{equation}\label{E:DefIndMultiInA-parameter}
    m\coloneqq \begin{cases}
        m_{\varphi_V+\widetilde{\varphi}_V}\left(\rho_1 S_{a_1}\boxtimes S_{b_1}\right) \quad &\text{if }\rho_1\not\cong\widetilde{\rho}_1;\\
        \frac{1}{2}\cdot m_{\varphi_V+\widetilde{\varphi}_V}\left(\rho_1 S_{a_1}\boxtimes S_{b_1}\right) \quad &\text{if }\rho_1\cong\widetilde{\rho}_1.
    \end{cases}
\end{equation}
Then by the construction of non-discrete A-packets (see \cite[Prop. 2.4.3]{MR3135650}), one knows that $\pi$ is an irreducible summand of the unitary induction
\[
    u_{\rho_1}\left(a_1,b_1\right)^{\times m}\rtimes \pi^\flat,
\]
where $\pi^\flat$ belongs to the A-packet associated to $\psi_V^\flat = \psi_V - m\left(\rho_1 S_{a_1}\boxtimes S_{b_1} + \widetilde{\rho}_1 S_{a_1}\boxtimes S_{b_1}\right)$. We appeal to Proposition \ref{P:Frobenius-MackeyInequality} to analyze $m\left(u_{\rho_1}\left(a_1,b_1\right)^{\times m}\rtimes \pi^\flat,\sigma\right)$. For any $\alpha<md_{\rho_1}a_1$, the cuspidal representation $\left(\rho_1|\cdot|^{A_1+\frac{1}{2}}\right)^\vee$ is contained in the cuspidal support of $\left(u_{\rho_1}\left(a_1,b_1\right)^{\times m}\right)^{[\alpha],\vee}$ but not in the cuspidal support of $\sigma$; this implies that $\D_{\left(u_{\rho_1}\left(a_1,b_1\right)^{\times m}\right)^{[\alpha],\vee}}\left(\sigma\right)=0$. It then follows from Proposition \ref{P:Frobenius-MackeyInequality} that
\[
    m\left(u_{\rho_1}\left(a_1,b_1\right)^{\times m}\rtimes \pi^\flat,\sigma\right) \leq m\left(\pi^\flat, \D_{u_{\rho_1^\vee}\left(a_1,b_1-1\right)^{\times m}}\left(\sigma\right)\right).
\]
In particular, we have $\D_{u_{\rho_1^\vee}\left(a_1,b_1-1\right)^{\times m}}\left(\sigma\right)\neq0$. This implies that there is an irreducible representation $\sigma^\flat$ of an appropriate classical group, such that we have an embedding $\sigma \hookrightarrow u_{\rho_1^\vee}\left(a_1,b_1-1\right)^{\times m}\rtimes \sigma^\flat$.
Since $\sigma$ is a direct summand of the unitary induction $\tau_{\varphi_W}\rtimes\sigma_0$, by the second adjointness, there is a non-zero equivariant map
\[
    \tau_{\varphi_W}\boxtimes\sigma_0 \longrightarrow \Jac_{\overline{Q}}\left(u_{\rho_1^\vee}\left(a_1,b_1-1\right)^{\times m}\rtimes \sigma^\flat\right).
\]
Here $Q$ is the maximal parabolic subgroup of $G\left(W\right)$ stabilizing a $\dim\varphi_W$-dimensional totally isotropic subspace $Y$, and $\overline{Q}$ is the parabolic subgroup opposite to $Q$. Then we use the geometric lemma to compute
\[
    \Jac_{\overline{Q}}\left(u_{\rho_1^\vee}\left(a_1,b_1-1\right)^{\times m}\rtimes \sigma^\flat\right)\cong \left(\Jac_Q\left(u_{\rho_1}\left(a_1,b_1-1\right)^{\times m}\rtimes \left(\sigma^\flat\right)^\vee\right)\right)^\vee.
\]

\vskip 5pt

For simplicity of notation, we let $k\coloneqq m a_1\left(b_1-1\right)\dim\rho_1$ and $\ell\coloneqq \dim\varphi_W$. We also set $A'\coloneqq A_1-\frac{1}{2}$ and $B'\coloneqq B_1+\frac{1}{2}$.
Denote by $M_Q\cong \GL_\ell\times G\left(W^\flat\right)$ the Levi subgroup of $Q$, where $W^\flat$ is the space in the same Witt tower as $W$ and of appropriate dimension; also denote by $Q'$ the maximal parabolic subgroup of $G\left(W\right)$ stabilizing a $k$-dimensional totally isotropic subspace, so its Levi subgroup is $M_{Q'}\cong \GL_k\times G\left(W_0\right)$ for some space $W_0$. There is an $M_Q$-stable filtration on $\Jac_Q\left(u_{\rho_1}\left(a_1,b_1-1\right)^{\times m}\rtimes \left(\sigma^\flat\right)^\vee\right)$ indexed by $Q$-orbits on $G\left(W\right)/Q'$. Such $Q$-orbits are parametrized by triples $\left(\alpha,\beta,\gamma\right)$, where
\begin{itemize}
    \item $\alpha$, $\beta$, $\gamma\in\Z_{\geq 0}$ such that $\alpha+\beta+\gamma=k$;

    \vskip 5pt

    \item $\alpha+\gamma\leq \ell$, and $\beta\leq r_W-\ell$; here $r_W$ is the Witt index of $W$.
\end{itemize}
\vskip 5pt
For a $k$-dimensional totally isotropic subspace $X\subset W$, the triple $\left(\alpha_X,\beta_X,\gamma_X\right)$ associated to the $Q$-orbit of $X$ is determined by
\[
    \alpha_X \coloneqq \dim Y\cap X, \quad \beta_X\coloneqq \dim Y^\perp\cap X - \alpha_X, \quad \gamma_X\coloneqq k-\alpha_X-\beta_X.
\]
The stabilizer of $X$ in the Levi subgroup $M_Q$ is a parabolic subgroup of $M_Q$ with Levi subgroup isomorphic to
\[
    \underbrace{\Big(\GL_{\alpha_X}\times\GL_{\ell-\alpha_X-\gamma_X}\times\GL_{\gamma_X}\Big)}_{\subset \GL_\ell} \times \underbrace{\left(\GL_{\beta_X}\times G\left(W^\flat\right)\right)}_{\subset G\left(W_0\right)}.
\]
Therefore, given the parabolically induced representation $u_{\rho_1}\left(a_1,b_1-1\right)^{\times m}\rtimes \left(\sigma^\flat\right)^\vee$, by using the geometric lemma and the result of Kret--Lapid (Lemma \ref{L:KretLapidJacFormula}) one knows that there is a filtration on $\Jac_Q\left(u_{\rho_1}\left(a_1,b_1-1\right)^{\times m}\rtimes \left(\sigma^\flat\right)^\vee\right)$, with each graded piece of the form
\begin{align*}
        J\Big(\vec{z}_1,\cdots,\vec{z}_m;\vec{w}_1,\cdots,\vec{w}_m;\Xi\boxtimes\sigma^\sharp\Big)=&\\
    \Big({\rm L}^-\left(\vec{z}_1\right)\times{\rm L}^-\left(\vec{z}_2\right)\times\cdots\times {\rm L}^-\left(\vec{z}_m\right)\times &\Xi \times {\rm L}^+\left(\vec{w}_m\right)^\sim\times\cdots\times{\rm L}^+\left(\vec{w}_2\right)^\sim\times{\rm L}^+\left(\vec{w}_1\right)^\sim\Big)\\
    \boxtimes&\Big({\rm L}\left(\vec{z}_1,\vec{w}_1\right)\times{\rm L}\left(\vec{z}_2,\vec{w}_2\right)\times\cdots\times {\rm L}\left(\vec{z}_m,\vec{w}_m\right)\rtimes \sigma^\sharp\Big).
\end{align*}
Here for each $1\leq j\leq m$, $\vec{z}_j = \left(z_{j,1},z_{j,2},\cdots,z_{j,b_1-1}\right)$ and $\vec{w}_j = \left(w_{j,1},w_{j,2},\cdots,w_{j,b_1-1}\right)$ are two strictly ascending sequences of half-integers, with coordinates satisfying $z_{j,k}, w_{j,k}\equiv A' \pmod \Z$ and $w_{j,k}<z_{j,k}$; and
\[
        \begin{cases}
            {\rm L}^+\left(\vec{w}_j\right)\coloneqq \LS\left(\Delta_{\rho_1}[w_{j,1},-A'], \Delta_{\rho_1}[w_{j,2},-A'+1],\cdots, \Delta_{\rho_1}[w_{j,b_1-1},-B']\right);\\[5pt]

            {\rm L}^-\left(\vec{z}_j\right)\coloneqq \LS\left(\Delta_{\rho_1}[B',z_{j,1}], \Delta_{\rho_1}[B'+1,z_{j,2}],\cdots, \Delta_{\rho_1}[A',z_{j,b_1-1}]\right);\\[5pt]

            {\rm L}\left(\vec{z}_j,\vec{w}_j\right)\coloneqq \LS\left(\Delta_{\rho_1}[z_{j,1}-1,w_{j,1}+1], \Delta_{\rho_1}[z_{j,2}-1,w_{j,2}+1],\cdots, \Delta_{\rho_1}[z_{j,b_1-1}-1,w_{j,b_1-1}+1]\right).
        \end{cases}
     \]
The representation $\Xi\boxtimes \sigma^\sharp$ is an irreducible subquotient of a Jacquet module of $\left(\sigma^\flat\right)^\vee$ along a maximal parabolic subgroup. If there is some $1\leq j\leq m$, such that
\[
    \left(w_{j,1},w_{j,2},\cdots,w_{j,b_1-1}\right) \neq \left(z_{j,1}-1,z_{j,2}-1,\cdots,z_{j,b_1-1}-1\right),
\]
i.e. ${\rm L}\left(\vec{z}_j,\vec{w}_j\right)$ is not the trivial representation of the trivial group $\GL_0$, then by our assumption on infinitesimal characters, we have
\[
    \Hom_{G\left(W_0\right)}\left(\sigma_0,\left({\rm L}\left(\vec{z}_1,\vec{w}_1\right)\times{\rm L}\left(\vec{z}_2,\vec{w}_2\right)\times\cdots\times {\rm L}\left(\vec{z}_m,\vec{w}_m\right)\rtimes \sigma^\sharp\right)^\vee\right)=0.
\]
Hence, there is a graded piece $J\left(\vec{z}_1,\cdots,\vec{z}_m;\vec{w}_1,\cdots,\vec{w}_m;\Xi\boxtimes\sigma^\sharp\right)$, with the property that $\left(w_{j,1},w_{j,2},\cdots,w_{j,b_1-1}\right) =\left(z_{j,1}-1,z_{j,2}-1,\cdots,z_{j,b_1-1}-1\right)$ for all $1\leq j\leq m$, such that
\[
    \Hom_{\GL_{\ell}\times G\left(W_0\right)}\left(\tau_{\varphi_W}\boxtimes\sigma_0,J\left(\vec{z}_1,\cdots,\vec{z}_m;\vec{w}_1,\cdots,\vec{w}_m;\Xi\boxtimes\sigma^\sharp\right)^\vee\right)\neq0.
\]
In particular, by using the MVW-involution this implies that
\[
    \Hom_{\GL_\ell}\left(\tau_{\varphi_W},{\rm L}^-\left(\vec{z}_1\right)^\vee\times\cdots\times {\rm L}^-\left(\vec{z}_m\right)^\vee\times \Xi^\vee \times \prescript{\iota}{}{{\rm L}^+\left(\vec{w}_m\right)}\times\cdots\times\prescript{\iota}{}{{\rm L}^+\left(\vec{w}_1\right)}\right)\neq0.
\]
By Proposition \ref{P:TechnicalCountingMultiApp}, the non-vanishing of this Hom space implies that
\[
    m\leq\begin{cases}
        m_{\varphi_W}\left(\rho_1^\vee S_{a_1}\boxtimes S_{b_1-1}\right) + m_{\varphi_W}\left(\widetilde{\rho}_1^\vee S_{a_1}\boxtimes S_{b_1-1}\right), \quad &\text{if } \rho_1\not\cong\widetilde{\rho}_1;\\
        m_{\varphi_W}\left(\rho_1^\vee S_{a_1}\boxtimes S_{b_1-1}\right), \quad &\text{if }\rho_1\cong\widetilde{\rho}_1.
    \end{cases}
\]
To streamline the discussion, we defer the proof of this technical result (Proposition \ref{P:TechnicalCountingMultiApp}) to Appendix \ref{App:ProofTechnicalCountingMulti}. This completes the proof of the desired inequality (\ref{E:ReductionToGoodParityOnlyIfInequality-1}) in this subcase.

\vskip 5pt

\emph{Subcase 1.2:} $b_1<\beta$, and $d_j<\beta$ for any $1\leq j\leq s$. In this subcase, there is some $i_0\in\left\{1,2,\cdots,r\right\}$, such that $b_{i_0}=\beta$. Let $\varphi_V^\dagger\coloneqq \varphi_V - \rho_{i_0} S_{a_{i_0}}\boxtimes S_{b_{i_0}} + \rho_{i_0} S_{a_{i_0}}\boxtimes S_{b_{i_0}-2}$, and $\psi_V^\dagger \coloneqq \varphi_V^\dagger + \psi_{V,0} + \widetilde{\varphi}_V^\dagger$. One can easily check that the pair of A-parameters $\left(\psi_V^\dagger,\psi_W\right)$ still satisfies the conditions in Theorem \ref{T:ReductionToGoodParity}. Let $\pi^\dagger\coloneqq\pi\Big(\phi_\pi^\dagger,\eta_\pi\big|_{\calA_{\phi_\pi^\dagger}}\Big)$, where $\phi_\pi^\dagger\coloneqq \phi_{\varphi_V^\dagger} + \phi_{\pi_0} + \phi_{\widetilde{\varphi}_V^\dagger}$. Again, thanks to \cite[Prop. 5.3]{MR2504024}, we have $\pi^\dagger\subset\tau_{\varphi_V^\dagger}\rtimes\pi_0$. By Theorem \ref{T:LLCProperties}(2), one knows there is a surjection
\[
    \Delta_{\rho_{i_0}}\left[A_{i_0},-B_{i_0}\right] \times \Delta_{\widetilde{\rho}_{i_0}}\left[A_{i_0},-B_{i_0}\right] \rtimes \pi^\dagger\twoheadrightarrow \pi.
\]
Therefore we have
\[
    m\left(\Delta_{\rho_{i_0}}\left[A_{i_0},-B_{i_0}\right] \times \Delta_{\widetilde{\rho}_{i_0}}\left[A_{i_0},-B_{i_0}\right] \rtimes \pi^\dagger,\sigma\right) \geq m\left(\pi,\sigma\right) > 0.
\]
Since for any $1\leq j\leq s$, $b_{i_0}>d_j$. By an argument similar to the proof of Lemma \ref{L:ExpInJacModLowerBound}, one can show that $\D_{\left(\Delta_{\rho_{i_0}}\left[A_{i_0},-B_{i_0}\right] \times \Delta_{\widetilde{\rho}_{i_0}}\left[A_{i_0},-B_{i_0}\right]\right)^{[\alpha],\vee}}\left(\sigma\right) = 0$
for any $\alpha<2d_{\rho_{i_0}}a_{i_0}$, where $d_{\rho_{i_0}}\coloneqq \dim\rho_{i_0}$. It then follows from Proposition \ref{P:Frobenius-MackeyInequality} that
\[
    m\left(\Delta_{\rho_{i_0}}\left[A_{i_0},-B_{i_0}\right] \times \Delta_{\widetilde{\rho}_{i_0}}\left[A_{i_0},-B_{i_0}\right] \rtimes \pi^\dagger,\sigma\right) \leq m\left(\pi^\dagger,\sigma\right).
\]
Since $\dim\varphi_V^\dagger+\dim\varphi_W<\dim\varphi_V+\dim\varphi_W$, by the induction hypothesis, the pair of A-parameters $\left(\varphi_V^\dagger+\widetilde{\varphi}_V^\dagger, \left(\varphi_W+\widetilde{\varphi}_W\right)^\vee\right)$ is relevant. Note that $\rho_1 S_{a_1}\not\cong \rho_{i_0} S_{a_{i_0}}$. By Lemma \ref{L:AlterDescriptionRelevance}, we have
\[
    0\leq \Lambda\left(\rho_1 S_{a_1}, b_1-1;\left(\varphi_W+\widetilde{\varphi}_W\right)^\vee,\varphi_V^\dagger+\widetilde{\varphi}_V^\dagger\right) = \Lambda\left(\rho_1 S_{a_1}, b_1-1;\left(\varphi_W+\widetilde{\varphi}_W\right)^\vee,\varphi_V+\widetilde{\varphi}_V\right).
\]
This implies the desired inequality (\ref{E:ReductionToGoodParityOnlyIfInequality-1}).

\vskip 5pt

\emph{Subcase 1.3:} There is some $j_0\in\left\{1,2,\cdots,s\right\}$, such that $d_{j_0}=\beta$. The proof of this subcase is similar to the previous subcase. Let $\varphi_W^-\coloneqq \varphi_W - \rho'_{j_0} S_{c_{j_0}}\boxtimes S_{d_{j_0}} + \rho'_{j_0} S_{c_{j_0}}\boxtimes S_{d_{j_0}-2}$, and $\psi_W^- \coloneqq \varphi_W^- + \psi_{W,0} + \widetilde{\varphi}_W^-$. One can easily check that the pair of A-parameters $\left(\psi_V,\psi_W^-\right)$ still satisfies the conditions in Theorem \ref{T:ReductionToGoodParity}. Let $\sigma^-\coloneqq\pi\Big(\phi_\sigma^-,\eta_\sigma\big|_{\calA_{\phi_\sigma^-}}\Big)$, where $\phi_\sigma^-\coloneqq \phi_{\varphi_W^-} + \phi_{\sigma_0} + \phi_{\widetilde{\varphi}_W^-}$. It follows from \cite[Prop. 5.3]{MR2504024} that $\sigma^-\subset\tau_{\varphi_W^-}\rtimes\sigma_0$. By Theorem \ref{T:LLCProperties}(2), one knows there is a surjection
\[
    \Delta_{\rho'_{j_0}}\left[C_{j_0},-D_{j_0}\right] \times \Delta_{\widetilde{\rho}'_{j_0}}\left[C_{j_0},-D_{j_0}\right] \rtimes \sigma^-\twoheadrightarrow \sigma.
\]
Using Proposition \ref{P:Frobenius-MackeyInequality} again, one can show that
\[
    m\left(\sigma,\pi\right)\leq m\left(\Delta_{\rho'_{j_0}}\left[C_{j_0},-D_{j_0}\right] \times \Delta_{\widetilde{\rho}'_{j_0}}\left[C_{j_0},-D_{j_0}\right] \rtimes \sigma^-,\pi\right) \leq m\left(\sigma^-,\pi\right).
\]
Then by the induction hypothesis, $\left(\varphi_V+\widetilde{\varphi}_V, \left(\varphi_W^-+\widetilde{\varphi}^-_W\right)^\vee\right)$ is relevant. It follows from Lemma \ref{L:AlterDescriptionRelevance} that
\[
    0\leq \Lambda\left(\rho_1 S_{a_1}, b_1-1;\left(\varphi_W^-+\widetilde{\varphi}_W^-\right)^\vee,\varphi_V+\widetilde{\varphi}_V\right) \leq \Lambda\left(\rho_1 S_{a_1}, b_1-1;\left(\varphi_W+\widetilde{\varphi}_W\right)^\vee,\varphi_V+\widetilde{\varphi}_V\right).
\]
This implies the desired inequality (\ref{E:ReductionToGoodParityOnlyIfInequality-1}).

\vskip 5pt

\emph{Case 2:} $b_1=1$. The proof of this case is similar to (and much easier than) Subcase 1.1. Let $m\in\Z_{>0}$ be as in (\ref{E:DefIndMultiInA-parameter}), and $\psi_V^\ddagger = \psi_V - m\left(\rho_1 S_{a_1}\boxtimes S_{b_1} + \widetilde{\rho}_1 S_{a_1}\boxtimes S_{b_1}\right)$. Then $\psi_V^\ddagger = \varphi_V^\ddagger + \psi_{V,0}+ \widetilde{\varphi}_V^\ddagger$, where
\[
    \varphi_V^\ddagger\coloneqq \begin{cases}
       \varphi_V - m_{\varphi_V}\left(\rho_1 S_{a_1}\right)\rho_1 S_{a_1} - m_{\varphi_V}\left(\widetilde{\rho}_1S_{a_1}\right) \widetilde{\rho}_1S_{a_1}, \quad &\text{if } \rho_1\not\cong \widetilde{\rho}_1;\\
       \varphi_V - m\rho_1 S_{a_1}, \quad &\text{if } \rho_1\cong \widetilde{\rho}_1.
    \end{cases}
\]
It is easy to check that the pair of A-parameters $\left(\psi_V^\ddagger,\psi_W\right)$ still satisfies the conditions in Theorem \ref{T:ReductionToGoodParity}. Let $\pi^\ddagger\coloneqq\pi\Big(\phi_\pi^\ddagger,\eta_\pi\big|_{\calA_{\phi_\pi^\ddagger}}\Big)$, where $\phi_\pi^\ddagger\coloneqq \phi_{\varphi_V^\ddagger} + \phi_{\pi_0} + \phi_{\widetilde{\varphi}_V^\ddagger}$. Then \cite[Prop. 5.3]{MR2504024} implies that $\pi^\ddagger\subset\tau_{\varphi_V^\ddagger}\rtimes\pi_0$. By Theorem \ref{T:LLCProperties}, one knows that $\pi$ is a direct summand of the unitary induction $\Delta_{\rho_1}\left[A_1,-A_1\right]^{\times m}\rtimes \pi^\ddagger$. Using Proposition \ref{P:Frobenius-MackeyInequality} again (together with a standard cuspidal support argument), one can show that
\[
    m\left(\pi,\sigma\right) \leq m\left(\Delta_{\rho_1}\left[A_1,-A_1\right]^{\times m}\rtimes \pi^\ddagger,\sigma\right) \leq m\left(\pi^\ddagger,\sigma\right).
\]
Since $\dim\varphi_V^\ddagger + \dim \varphi_W<\dim\varphi_V + \dim \varphi_W$, the induction hypothesis tells us that $m\left(\pi_0,\sigma_0\right)$ is non-zero and $\left(\varphi_V^\ddagger+\widetilde{\varphi}_V^\ddagger, \left(\varphi_W+\widetilde{\varphi}_W\right)^\vee\right)$ is relevant. Hence, the pair $\left(\varphi_V+\widetilde{\varphi}_V, \left(\varphi_W+\widetilde{\varphi}_W\right)^\vee\right)$ is also relevant. This completes the proof of Theorem \ref{T:ReductionToGoodParity}.


\vskip 10pt


\section{Relevance implies period: Good parity case}\label{sec:relevance_implies_period}

With the reduction step done, we can now focus on the good parity A-parameters. In this section, we prove the ``relevance implies period'' direction. Namely, we want to show that if \(\left(\psi_V,\psi_W^{\vee}\right)\) is a pair of relevant A-parameters, then the unique pair of representations \(\left(\pi,\sigma\right)\in\left(\LP_{\phi_{\psi_V}}\times\LP_{\phi_{\psi_W}}\right)^{\rel}\) indexed by the distinguished character as in Section \ref{sub:local_gan_gross_prasad_conjecture_tempered_case} satisfies \(m\left(\pi,\sigma\right)\neq 0\).

\subsection{Multiplicity one theorem for standard modules}

Firstly, we record a multiplicity one theorem for standard modules, which was proved by M{\oe}glin--Waldspurger \cite[Sect. 1.3]{MR3155346} for (special) orthogonal groups. Here we formulate this theorem for both orthogonal and unitary groups; the proof of M{\oe}glin--Waldspurger extends to the unitary case with only minor modifications.

\begin{Prop}\label{mwprop13}
Let $W\subset V$ be a pair of non-degenerate Hermitian spaces defined over $E$ as in Section \ref{sub:bessel_subgroups_and_bessel_periods}. Let
\[
    \pi=\LS\left(\Delta_{\rho_1}[x_1,y_1],\cdots,\Delta_{\rho_r}[x_r,y_r];\pi_0\right) \quad \text{and} \quad \sigma=\LS(\Delta_{\rho'_1}[x'_1,y'_1],\cdots,\Delta_{\rho'_s}[x'_s,y'_s];\sigma_0)
\]
be irreducible representations of $G^0\left(V\right)$ and $G^0\left(W\right)$ respectively; in the case that one of these groups is a special odd orthogonal group, we extend the representation to the full orthogonal group by letting $\left\{\pm1\right\}$ act trivially. Then we have
\[
    m\left(\calM\left(\pi\right),\calM\left(\sigma\right)\right) = m\left(\pi_0,\sigma_0\right).
\]
\end{Prop}

Combining this result with Theorem \ref{T:LLCProperties}(2) and Theorem \ref{T:TemperedGGP}, we can deduce the following consequence:

\begin{Lem}\label{focuslemma}
Let $W\subset V$ be a pair of non-degenerate Hermitian spaces defined over $E$ as in Section \ref{sub:bessel_subgroups_and_bessel_periods}, and $\phi_V$, $\phi_W$ be L-parameters of $G^0\left(V\right)$, $G^0\left(W\right)$ respectively, not necessarily tempered. Then we have
\[
    \sum_{\left(\pi,\sigma\right)\in \left(\LP_{\phi_{V}}\times \LP_{\phi_{W}}\right)^{\rel}}m\left(\pi,\sigma\right)\leq 1.
\]
If equality holds, then the unique pair of representations $\left(\pi,\sigma\right)$ such that $m\left(\pi,\sigma\right)\neq 0$ is indexed by the distinguished character $\eta_{\phi_V,\phi_W}\boxtimes\eta_{\phi_W,\phi_V}\in\Irr\left(\calA_{\phi_V}\right)\times\Irr\left(\calA_{\phi_W}\right)$ as defined in Section \ref{sub:local_gan_gross_prasad_conjecture_tempered_case}.
\end{Lem}

\vskip 5pt

\begin{proof}
Suppose that we have decompositions
\[
    \phi_V = \left(\phi_1+\cdots+\phi_r\right) + \phi_{V,0} + \left(\widetilde{\phi}_1+\cdots+\widetilde{\phi}_r\right), \quad \phi_W = \left(\phi'_1+\cdots+\phi'_s\right) + \phi_{W,0} + \left(\widetilde{\phi}'_1+\cdots+\widetilde{\phi}'_s\right)
\]
as in Theorem \ref{T:LLCProperties}(2), where $\phi_{V,0}$ and $\phi_{W,0}$ are the tempered parts of $\phi_V$ and $\phi_W$. Then for each pair of representations $\left(\pi,\sigma\right)\in\left(\LP_{\phi_{V}}\times \LP_{\phi_{W}}\right)^{\rel}$, there is a unique pair of representations $\left(\pi_0,\sigma_0\right)\in\left(\LP_{\phi_{V,0}}\times \LP_{\phi_{W,0}}\right)^{\rel}$, such that
\[
    \pi=\LS\left(\tau_1,\cdots,\tau_r;\pi_0\right)\quad \text{and} \quad\sigma=\LS\left(\tau'_1,\cdots,\tau'_s;\sigma_0\right),
\]
where $\tau_i$ and $\tau'_j$ are irreducible essentially discrete series representations corresponding to $\phi_i$ and $\phi'_j$ respectively. It follows from Proposition \ref{mwprop13} that
\[
     \sum_{\left(\pi,\sigma\right)\in \left(\LP_{\phi_{V}}\times \LP_{\phi_{W}}\right)^{\rel}}m\left(\pi,\sigma\right) \leq  \sum_{\left(\pi_0,\sigma_0\right)\in \left(\LP_{\phi_{V,0}}\times \LP_{\phi_{W,0}}\right)^{\rel}}m\left(\pi_0,\sigma_0\right)=1.
\]
Here in the last equality, we make use of Theorem \ref{T:TemperedGGP}. This proves the first statement. The second statement follows from the combination of Theorem \ref{T:LLCProperties}(2) and Theorem \ref{T:TemperedGGP}.

\end{proof}

\vskip 5pt

With this lemma at hand, now to complete the proof of Conjecture \ref{Conj:Non-temperedGGP}(1), it only remains to show:

\begin{Thm}\label{T:RelevanceImpliesPeriod-GoodParity}
Let $W\subset V$ be a pair of non-degenerate Hermitian spaces defined over $E$ as in Section \ref{sub:bessel_subgroups_and_bessel_periods}, and $\psi_V$, $\psi_W$ be good parity A-parameters of $G^0\left(V\right)$, $G^0\left(W\right)$ respectively. If $\left(\psi_V,\psi_W^\vee\right)$ is relevant, then
\[
    \sum_{\left(\pi,\sigma\right)\in \left(\LP_{\phi_{\psi_V}}\times \LP_{\phi_{\psi_W}}\right)^\rel} m\left(\pi,\sigma\right)\neq0.
\]   
\end{Thm}

\vskip 5pt

\subsection{Proof of Theorem \ref{T:RelevanceImpliesPeriod-GoodParity}: I. a special case} 
\label{sub:proof_of_theorem_t:relevanceimpliesperiod-goodparity}

Now we are ready to prove Theorem \ref{T:RelevanceImpliesPeriod-GoodParity}. Suppose that we are given a pair of relevant A-parameters $\left(\psi_V,\psi_W^\vee\right)$, where
\[
    \psi_V = \sum_{i=1}^r \rho_i S_{a_i}\boxtimes S_{b_i} \quad \text{and} \quad \psi_W = \sum_{j=1}^s \rho'_j S_{c_j}\boxtimes S_{d_j},
\]
and these sums are ordered so that $b_1\geq b_2\geq \cdots\geq b_r$ and $d_1\geq d_2\geq \cdots\geq d_s$. Without loss of generality we may assume that $b_1\geq d_1$. For simplicity of notation we set 
\[
    \begin{cases}
       A_i \coloneqq \frac{a_i+b_i}{2}-1,\\
       B_i \coloneqq \frac{a_i-b_i}{2} 
    \end{cases}
    \quad\text{and}\quad
    \begin{cases}
      C_j \coloneqq \frac{c_j+d_j}{2}-1,\\
      D_j \coloneqq \frac{c_j-d_j}{2}.  
    \end{cases}
\]
We shall prove the desired conclusion by induction on the integer
\[
    N\left(\psi_V,\psi_W\right):=\sum^{r}_{i=1}\left(b_i-1\right)+\sum^{s}_{j=1}\left(d_j-1\right).
\]
The base case is when \(N\left(\psi_V,\psi_W\right)=0\), i.e. when both \(\psi_V\) and \(\psi_W\) are tempered, which follows from Theorem \ref{T:TemperedGGP}. Before we proceed to the induction step, we need to deal with a special case first: the case that
\[
    2=b_1=\cdots=b_k>b_{k+1}=\cdots=b_r=1
\]
and $d_1=1$ (in particular $\psi_W$ is tempered). The rest of this subsection is devoted to the proof of this special case. Depending on the multiplicities of non-tempered Jordan blocks, there are three subcases. We treat them separately.

\vskip 5pt

\emph{Subcase 1:} there exists some $i_0\in\left\{1,\cdots,k\right\}$ such that $m_{\psi_V}\left(\rho_{i_0} S_{a_{i_0}}\boxtimes S_2\right) < m_{\psi_W^\vee}\left(\rho_{i_0} S_{a_{i_0}}\right)$. For simplicity of notation, we may assume that $i_0=1$. Let
\[
    \psi_V^- \coloneqq \psi_V - \rho_1 S_{a_1}\boxtimes S_2, \quad \text{and} \quad \psi_{W,0} \coloneqq \psi_W-2\rho_1^\vee S_{a_1}.
\]
Note that our assumption $m_{\psi_V}\left(\rho_{i_0} S_{a_{i_0}}\boxtimes S_2\right) < m_{\psi_W^\vee}\left(\rho_{i_0} S_{a_{i_0}}\right)$ implies that the pair of A-parameters $\left(\psi_V^-,\psi_{W,0}^\vee\right)$ is also relevant; and we have $N\left(\psi_V^-,\psi_W\right)<N\left(\psi_V,\psi_W\right)$. By the induction hypothesis, there exists a unique pair $\left(\pi^-,\sigma_0\right)\in\left(\LP_{\phi_{\psi_V^-}}\times\LP_{\psi_{W,0}}\right)^\rel$ such that $m\left(\pi^-,\sigma_0\right)\neq 0$. Let
\[
    \pi\hookrightarrow\Delta_{\rho_1}\left[B_1,-A_1\right]\rtimes\pi^-
\]
be the unique irreducible subrepresentation. By Theorem \ref{T:LLCProperties}, one knows that $\pi\in\LP_{\phi_{\psi_V}}$, and any irreducible summand of $u_{\rho_1^\vee}\left(a_1,1\right)\rtimes\sigma_0$ belongs to the L-packet $\LP_{\psi_W}$. Then, applying Proposition \ref{P:Frobenius-MackeyInequality}, we have
\[
    m\left(u_{\rho_1^\vee}\left(a_1,1\right)\rtimes\sigma_0,\pi\right)\geq m\left(\sigma_0,\D_{\Delta_{\rho_1}\left[B_1,-A_1\right]}\left(\pi\right)\right) \geq m\left(\sigma_0,\pi^-\right).
\]
Here in the last inequality, we have made use of Frobenius reciprocity, which implies that $\pi^-$ is an irreducible quotient of $\D_{\Delta_{\rho_1}\left[B_1,-A_1\right]}\left(\pi\right)$. It follows that there is an irreducible representation $\sigma\in\LP_{\psi_W}$ such that $m\left(\pi,\sigma\right)\neq 0$.

\vskip 5pt

\emph{Subcase 2:} for any $i\in\left\{1,\cdots,k\right\}$, one has $m_{\psi_V}\left(\rho_{i} S_{a_{i}}\boxtimes S_2\right) = m_{\psi_W^\vee}\left(\rho_{i} S_{a_{i}}\right)$; and there exists some $i_0\in\left\{1,\cdots,k\right\}$ such that $m_{\psi_V}\left(\rho_{i_0} S_{a_{i_0}}\boxtimes S_2\right) = m_{\psi_W^\vee}\left(\rho_{i_0} S_{a_{i_0}}\right) \geq 2$. For simplicity of notation we assume $i_0=1$. Let
\[
    \psi_{V,0} \coloneqq \psi_V - 2\rho_1 S_{a_1}\boxtimes S_2, \quad \text{and} \quad \psi_{W,0} \coloneqq \psi_W-2\rho_1^\vee S_{a_1}.
\]
Note that $\left(\psi_{V,0},\psi_{W,0}^\vee\right)$ is still relevant, and $N\left(\psi_{V,0},\psi_{W,0}\right)<N\left(\psi_V,\psi_W\right)$. The induction hypothesis implies that there exists a pair $\left(\pi_0,\sigma_0\right)\in\left(\LP_{\phi_{\psi_{V,0}}}\times\LP_{\psi_{W,0}}\right)^\rel$ such that $m\left(\pi_0,\sigma_0\right)\neq 0$. Consider the parabolically induced representation $u_{\rho_1}\left(a_1,2\right)\rtimes \pi_0$, and set $\sigma\coloneqq\pi\left(\psi_W,\eta_{\psi_W,\phi_{\psi_V}}\right)$. By Proposition \ref{P:unitaryInd}, $u_{\rho_1}\left(a_1,2\right)\rtimes \pi_0$ contains $\pi\coloneqq \pi\left(\phi_{\psi_V},\eta_{\phi_{\psi_V},\psi_W}\right)$ as an irreducible summand, and by Theorem \ref{T:LLCProperties}(1) $\sigma$ is an irreducible summand of $u_{\rho_1^\vee}\left(a_1,1\right)\rtimes\sigma_0$ and belongs to the L-packet $\LP_{\psi_W}$. Then, applying Proposition \ref{P:Frobenius-MackeyInequality} and Corollary \ref{C:VanishingLemma}, we have
\[
    m\left(u_{\rho_1}\left(a_1,2\right)\rtimes \pi_0,\sigma\right) = m\left(\pi_0, \D_{u_{\rho_1^\vee}\left(a_1,1\right)}\left(\sigma\right)\right)\geq m\left(\pi_0,\sigma_0\right).
\]
Here in the last inequality, once again we have made use of Frobenius reciprocity, which implies that $\sigma_0$ is an irreducible quotient of $\D_{u_{\rho_1^\vee}\left(a_1,1\right)}\left(\sigma\right)$. 

\vskip 5pt

If $m_{\psi_V}\left(\rho_{i_0} S_{a_{i_0}}\boxtimes S_2\right) > 2$, then Proposition \ref{P:unitaryInd} asserts that $\pi= u_{\rho_1}\left(a_1,2\right)\rtimes \pi_0$, and we have $m\left(\pi,\sigma\right)\neq 0$ as desired.
So now it only remains to consider the situation that $m_{\psi_V}\left(\rho_{i_0} S_{a_{i_0}}\boxtimes S_2\right) = 2$. Let
\[
    \psi'_V\coloneqq \psi_{V,0} + \rho_1 S_{a_1}\boxtimes S_2 + \rho_1 S_{a_1+1} + \rho_1 S_{a_1-1}
\]
as in Proposition \ref{P:unitaryInd}. Then irreducible summands of $u_{\rho_1}\left(a_1,2\right)\rtimes \pi_0$ other than $\pi$ all belong to the L-packet $\LP_{\phi_{\psi'_V}}$. Note that $\left(\psi'_V,\psi_W^\vee\right)$ is still relevant, and $N\left(\psi'_V,\psi_W\right)<N\left(\psi_V,\psi_W\right)$. By the induction hypothesis, to show that $m\left(\pi,\sigma\right)\neq 0$ it suffices to show that $\eta_{\psi_W,\phi_{\psi_V}}\neq \eta_{\psi_W,\phi_{\psi'_V}}$. Let \(\alpha\) be the basis element in \(A_{\psi_W}\) corresponding to \(\rho_1^\vee S_{a_1}\). Then by \cite[Lem. A.4]{MR3714507}, we have
\[
    \frac{\eta_{\psi_W,\phi_{\psi'_V}}(\alpha)}{\eta_{\psi_W,\phi_{\psi_V}}\left(\alpha\right)}=\epsilon\left(\rho_1^\vee S_{a_1}\otimes\left(\rho_1 S_{a_1+1}+\rho_1 S_{a_1-1}\right)\right)=-1.
\]
This implies that $\eta_{\psi_W,\phi_{\psi_V}}\neq \eta_{\psi_W,\phi_{\psi'_V}}$, and hence we have $m\left(\pi,\sigma\right)\neq 0$ as desired.

\vskip 5pt

\emph{Subcase 3:} for any $i\in\left\{1,\cdots,k\right\}$, one has $m_{\psi_V}\left(\rho_{i} S_{a_{i}}\boxtimes S_2\right) = m_{\psi_W^\vee}\left(\rho_{i} S_{a_{i}}\right)=1$. In this subcase we appeal to Proposition \ref{P:non-unitaryInd-1}. Let 
\[
    \psi_V^-\coloneqq \psi_V-\rho_1 S_{a_1}\boxtimes S_2 \quad \text{and}\quad \psi'_V \coloneqq \psi_V^- + \rho_1 S_{a_1+1} + \rho_1 S_{a_1-1}.
\]
We also let $\pi\coloneqq\pi\left(\phi_{\psi_V},\eta_{\phi_{\psi_V},\psi_W}\right)$, $\pi^-\coloneqq \pi\Big(\phi_{\psi_V^-},\eta_{\phi_{\psi_V^-},\psi_W}\Big)$ and $\sigma\coloneqq\pi\left(\psi_W,\eta_{\psi_W,\phi_{\psi_V}}\right)$. Then there is a short exact sequence
\[
    0\lra \pi \lra \Delta_{\rho_1}\left[B_1,-A_1\right]\rtimes \pi^- \lra \bigoplus_{\eta'} \pi\left(\phi_{\psi'_V},\eta'\right) \lra 0,
\]
where $\eta'\in\Irr\left(\calA_{\phi_{\psi'_V}}\right)$ runs over all extensions of $\eta_{\phi_{\psi_V^-},\psi_W}$. We claim that for any such $\eta'$, one has $m\left(\pi\left(\phi_{\psi'_V},\eta'\right),\sigma\right)=0$. Note that $\left(\psi'_V,\psi_W^\vee\right)$ is still relevant, and $N\left(\psi'_V,\psi_W\right)<N\left(\psi_V,\psi_W\right)$. Hence by the induction hypothesis, to show our claim, it suffices to show that $\eta_{\psi_W,\phi_{\psi'_V}}\neq \eta_{\psi_W,\phi_{\psi_V}}$. This follows from the same computation as in the previous subcase. With this claim at hand, the usual long exact sequence argument yields
\[
    m\left(\pi,\sigma\right)\geq m\left(\Delta_{\rho_1}\left[B_1,-A_1\right]\rtimes \pi^-,\sigma\right) = m\left(\pi^-,\sigma\right).
\]
Since $\left(\psi_V^-,\psi_W^\vee\right)$ is relevant, and $N\left(\psi_V^-,\psi_W\right)<N\left(\psi_V,\psi_W\right)$, by the induction hypothesis we have $m\left(\pi^-,\sigma\right)\neq 0$. Therefore we also have $m\left(\pi,\sigma\right)\neq 0$.

\vskip 5pt

This completes the proof of the special case.

\vskip 5pt


\subsection{Proof of Theorem \ref{T:RelevanceImpliesPeriod-GoodParity}: II. the induction step} 
\label{sub:proof_of_theorem_t:relevanceimpliesperiod-goodparity_ii_the_induction_step}

Now we may assume that $d_1\geq 2$. We now describe the idea to show the existence of periods between two representations. The following lemma is an analogue of \cite[Lem. 10.1]{chan2022SQIrr} for classical groups.

\begin{Lem}\label{L:3multi1Implies1multi1}
Let $W\subset V$ be a pair of non-degenerate Hermitian spaces defined over $E$ as in Section \ref{sub:bessel_subgroups_and_bessel_periods}, and $\pi$, $\sigma$ be representations of $G\left(V\right)$ and $G\left(W\right)$ respectively. Suppose that there exists a representation $\Pi$ of $G\left(V\right)$ and a representation $\Sigma$ of $G\left(W\right)$, such that:
\begin{itemize}
       \item $\pi$ is a quotient of $\Pi$, and $\sigma$ is a quotient of $\Sigma$;

       \vskip 5pt

       \item $m\left(\Pi,\Sigma\right)=m\left(\Pi,\sigma\right)=m\left(\pi,\Sigma\right)=1$.
   \end{itemize} 
Then we have $m\left(\pi,\sigma\right)=1$.  
\end{Lem}

\vskip 5pt

\begin{proof}
We fix non-zero equivariant maps
\[
    f\in\Hom_{\Bes\left(V,W\right)}\left(\Pi\boxtimes\sigma,\nu\right) \quad \text{and} \quad g\in\Hom_{\Bes\left(V,W\right)}\left(\pi\boxtimes\Sigma,\nu\right).
\]
Since $m\left(\Pi,\Sigma\right)=1$, the pull-back of $f$ along the quotient map $\Sigma\twoheadrightarrow\sigma$ and the pull-back of $g$ along the quotient map $\Pi\twoheadrightarrow\pi$ must be proportional. Therefore $f$ factors through the quotient map $\Pi\twoheadrightarrow\pi$, and it gives a non-zero element in $\Hom_{\Bes\left(V,W\right)}\left(\pi\boxtimes\sigma,\nu\right)$. Hence we have
\[
    0< m\left(\pi,\sigma\right) \leq m\left(\Pi,\sigma\right)=1.
\]
This implies that $m\left(\pi,\sigma\right)=1$ as desired.

\end{proof}

\vskip 5pt

To apply this lemma in our setting, we put
\[
    \left\{\begin{aligned}
        \psi_V^-&\coloneqq \psi_V-\rho_1 S_{a_1}\boxtimes S_{b_1} + \rho_1 S_{a_1}\boxtimes S_{b_1-2},\\
        \pi^-&\coloneqq \pi\left(\phi_{\psi_V^-},\eta_{\phi_{\psi_V^-},\phi_{\psi_W}}\right),\\
        \pi&\coloneqq\pi\left(\phi_{\psi_V},\eta_{\phi_{\psi_V},\phi_{\psi_W}}\right),\\
        \Pi&\coloneqq \Delta_{\rho_1}\left[A_1,-B_1\right]\rtimes\pi^-,
    \end{aligned}\right.
    \quad \text{and}\quad
    \left\{\begin{aligned}
        \psi_W^-&\coloneqq \psi_W-\rho'_1 S_{c_1}\boxtimes S_{d_1} + \rho'_1 S_{c_1}\boxtimes S_{d_1-2},\\
        \sigma^-&\coloneqq \pi\left(\phi_{\psi_W^-},\eta_{\phi_{\psi_W^-},\phi_{\psi_V}}\right),\\
        \sigma&\coloneqq\pi\left(\phi_{\psi_W},\eta_{\phi_{\psi_W},\phi_{\psi_V}}\right),\\
        \Sigma&\coloneqq \Delta_{\rho'_1}\left[C_1,-D_1\right]\rtimes\sigma^-.
    \end{aligned}\right.
\]
Then by Theorem \ref{T:LLCProperties}(2), $\Pi$ and $\Sigma$ are quotients of $\calM\left(\pi\right)$ and $\calM\left(\sigma\right)$ respectively. Therefore by the above lemma and Proposition \ref{mwprop13}, to show $m\left(\pi,\sigma\right)\neq0$, we only need to show that
\[
    m\left(\Pi,\sigma\right)\neq 0 \quad \text{and} \quad m\left(\pi,\Sigma\right)\neq 0.
\]
Applying Proposition \ref{P:Frobenius-MackeyInequality} and Corollary \ref{C:VanishingLemma}, we have
\[
    m\left(\Pi,\sigma\right) = m\left(\pi^-,\sigma\right).
\]
Note that $\left(\psi_V^-,\psi_W^\vee\right)$ is still relevant, and $N\left(\psi_V^-,\psi_W\right)<N\left(\psi_V,\psi_W\right)$. The induction hypothesis tells us that $m\left(\pi^-,\sigma\right)\neq 0$. Therefore it only remains to show $m\left(\pi,\Sigma\right)\neq 0$.

\vskip 5pt

If $d_1=b_1$, then the non-vanishing of $m\left(\pi,\Sigma\right)$ again follows from the combination of Proposition \ref{P:Frobenius-MackeyInequality} and Corollary \ref{C:VanishingLemma}. Thus we can assume that $d_1=b_1-1$. Since $\left(\psi_V,\psi_W^\vee\right)$ is relevant, one knows that $\rho_1^\vee S_{a_1}\boxtimes S_{b_1-1}\subset \psi_W$. Replacing $\psi_W^-$ by $\psi_W-\rho_1^\vee S_{a_1}\boxtimes S_{b_1-1} + \rho_1^\vee S_{a_1}\boxtimes S_{b_1-3}$ (and modifying $\sigma^-$, $\sigma$, $\Sigma$ accordingly) if necessary, we may further assume that $\rho'_1 S_{c_1}\cong\rho_1^\vee S_{a_1}$. Then by Proposition \ref{P:Frobenius-MackeyInequality} we have
\[
    m\left(\Sigma,\pi\right) \geq m\left(\sigma^-,\D_{\Delta_{\rho_1}\left[B_1,-A_1\right]}\left(\pi\right)\right)\geq m\left(\sigma^-,\pi^-\right).
\]
Here in the last inequality, we have made use of the fact that $\pi$ is the unique subrepresentation of $\Delta_{\rho_1}\left[B_1,-A_1\right]\rtimes\pi^-$, which implies that $\pi^-$ is a quotient of $\D_{\Delta_{\rho_1}\left[B_1,-A_1\right]}\left(\pi\right)$ by Frobenius reciprocity. Note that $\left(\psi_V^-,\left(\psi_W^-\right)^\vee\right)$ is still relevant, and $N\left(\psi_V^-,\psi_W^-\right)<N\left(\psi_V,\psi_W\right)$. It follows from the induction hypothesis that $m\left(\pi^-,\sigma^-\right)\neq 0$. Thus $m\left(\pi,\Sigma\right)\neq 0$, and by Lemma \ref{L:3multi1Implies1multi1} we have $m\left(\pi,\sigma\right)\neq 0$.

\vskip 5pt

This completes the proof of Theorem \ref{T:RelevanceImpliesPeriod-GoodParity}.


\vskip 10pt

\section{Period implies relevance: Good parity co-tempered case}\label{sec:period_implies_relevance}

The aim of this section is to show Conjecture \ref{Conj:Non-temperedGGP}(3). Thanks to Theorem \ref{T:ReductionToGoodParity}, it is sufficient to show the following:

\begin{Thm}\label{T:PeriodImpliesRelevanceGPCotempered}
In the setting of Conjecture \ref{Conj:Non-temperedGGP}, suppose that the A-parameters $\psi_V$ and $\psi_W$ are of good parity and co-tempered. If there exists $\left(\pi,\sigma\right)\in\left(\AP_{{\psi_V}}\times \AP_{{\psi_W}}\right)^\rel$ such that $m\left(\pi,\sigma\right)\neq 0$, then $\left(\psi_V,\psi_W^\vee\right)$ must be relevant.
\end{Thm}

\vskip 5pt

The rest of this section is devoted to the proof of this theorem. The idea is in the same spirit as the proof of Case 1 of Theorem \ref{T:ReductionToGoodParity}. We shall argue by induction on $\dim V + \dim W$.
Given a pair of representations $\pi$ and $\sigma$ in the associated A-packets such that $m\left(\pi,\sigma\right)\neq 0$, we will explicitly construct another pair of representations $\pi'$ and $\sigma'$ of smaller groups, such that they are ``more tempered'' and $m\left(\pi',\sigma'\right)\neq 0$. The induction hypothesis then tells us that A-parameters of $\pi'$ and $\sigma'$ are relevant. This already provides a large part of the information. To show the relevance of $\left(\psi_V,\psi_W^\vee\right)$, it only remains to show an inequality concerning multiplicities of certain Jordan blocks. The desired inequality follows from a delicate case-by-case analysis.

\subsection{The base case of the induction} 
\label{sub:the_base_case}

As explicated above, we shall prove Theorem \ref{T:PeriodImpliesRelevanceGPCotempered} by induction on $\dim V + \dim W$ (the outer induction). The base case of the induction is when
\[
    \dim V =2 \quad \text{and} \quad \dim W=1.
\]
We divide the discussion into cases according to the involution $\iota$:
\begin{itemize}
    \item When $\iota={\rm id}$, $G(V)$ and $G(W)$ are orthogonal groups. In this situation co-tempered A-parameters of $G(V)$ or $G(W)$ are tempered. Hence any pair of co-tempered A-parameters $\left(\psi_V,\psi_W^\vee\right)$ is relevant. 

    \vskip 5pt

    \item When $\iota \neq {\rm id}$, $G(V)$ and $G(W)$ are unitary groups. In this situation, any good parity co-tempered A-parameter of $G(V)$ is either tempered, or of the form $\chi\boxtimes S_2$, where $\chi$ is a conjugate orthogonal character of $W_E$; co-tempered A-parameters of $G(W)$ are tempered. If $\psi_V$, $\psi_W$ are tempered, then $\left(\psi_V,\psi_W^\vee\right)$ is automatically relevant. Otherwise, $\psi_V=\chi\boxtimes S_2$ for some $\chi$, and $\pi=\chi\circ\det_V$ is a character of $G(V)$. To have $m\left(\pi,\sigma\right)\neq0$, we must have $\sigma=\chi^\vee$. Hence $\psi_W=\chi^\vee$ and $\left(\psi_V,\psi_W^\vee\right)$ is relevant.
\end{itemize}
\vskip 5pt
This completes the proof of the base case of the induction.

\vskip 5pt

\subsection{The induction step of the proof of Theorem \ref{T:PeriodImpliesRelevanceGPCotempered}} 
\label{sub:the_induction_step_of_the_inner_induction}

Now we do the induction step. Let $N>3$ be a positive integer. Suppose that Theorem \ref{T:PeriodImpliesRelevanceGPCotempered} holds when $\dim V+\dim W<N$, we show that it still holds when $\dim V+\dim W=N$. Suppose that $\left(\pi,\sigma\right)\in\left(\AP_{{\psi_V}}\times \AP_{{\psi_W}}\right)^\rel$ such that $m\left(\pi,\sigma\right)\neq 0$. In the case that both $\pi$ and $\sigma$ are supercuspidal, Gan, Gross and Prasad already proved in \cite[Thm. 7.7]{MR4190046} that $\left(\psi_V,\psi_W^\vee\right)$ is relevant. Therefore we can assume that at least one of $\pi$ and $\sigma$ is not supercuspidal. Set
\[
    \Der\left(\pi,\sigma\right) \coloneqq \left\{\rho\in\calC\left(\GL\right)\,\big|\,D_\rho\left(\pi\right)\neq 0 \text{ or } D_\rho\left(\sigma\right) \neq 0\right\}.
\]
Let $\rho|\cdot|^{-x}$ be a supercuspidal representation in $\Der\left(\pi,\sigma\right)$ with the minimal exponent, where $\rho\in\calC_{\unit}\left(\GL\right)$ and $x\in\R$ (so the exponent $e\left(\rho|\cdot|^{-x}\right)=-x$). Since the multiplicity $m\left(\pi,\sigma\right)$ is symmetric, without loss of generality we may assume that $D_{\rho|\cdot|^{-x}}\left(\pi\right)\neq 0$. By Corollary \ref{C:co-temperedAtobeLem.4.1NStDer}, we have $x\in\frac{1}{2}\Z_{\geq 0}$, and $\rho\boxtimes S_{2x+1} \subset \psi_V$. Moreover, according to Corollary \ref{C:co-temperedShapeofAPara}, $\psi_V$ must be of the form 
\[
    \psi_V = \psi_{V,\rho,x} + m\rho \boxtimes S_{2x+1} + \rho\boxtimes S_{2x+3} + \cdots + \rho\boxtimes S_{2x+2k-1},
\]
where $\psi_{V,\rho,x}$ satisfies the property that for any integer $a\geq 2x+1$, $\rho\boxtimes S_a\not\subset \psi_{V,\rho,x}$; $m\coloneqq m_{\psi_V}\left(\rho\boxtimes S_{2x+1}\right)$, and $k\in\Z_{> 0}$. By the minimality of $-x$, we have $D_{\left(\varrho|\cdot|^{x+1/2}\right)^\vee}\left(\sigma\right)=0$ for any $\varrho\in\calC_{\unit}(\GL)$. To proceed the proof, we consider four cases separately.

\vskip 5pt

\emph{Case 1:} $m_{\psi_V}\left(\rho \boxtimes S_{2x+1}\right)\geq 2$. In this case, by Conjecture \ref{Conj:APacketsNProperties}(2) there is a unique irreducible representation $\pi_0\in\AP_{\psi_{V,0}}\left(G\left(V_0\right)\right)$, where $\psi_{V,0}=\psi_V-2\rho\boxtimes S_{2x+1}$ and $V_0$ is a Hermitian space in the Witt tower containing $V$ of appropriate dimension, such that $\pi$ is an irreducible summand of $Z_\rho\left[x,-x\right]\rtimes \pi_0$.
Then by Proposition \ref{P:Frobenius-MackeyInequality}, we have
\[
    m\left(Z_\rho\left[x,-x\right]\rtimes \pi_0, \sigma\right) \leq m\left(\pi_0,\D_{Z_{\rho^\vee}\left[x-\frac{1}{2},-x-\frac{1}{2}\right]}\left(\sigma\right)\right) + m\left(\pi_0, \D_{Z_{\rho^\vee}\left[x-\frac{1}{2},-x+\frac{1}{2}\right]}\left(\sigma\right)\right).
\]
The minimality of $-x$ implies that $D_{\rho^\vee|\cdot|^{-x-\frac{1}{2}}}\left(\sigma\right) =0$, so we have $\D_{Z_{\rho^\vee}\left[x-\frac{1}{2},-x-\frac{1}{2}\right]}\left(\sigma\right)= 0$, and it follows that
\[
    m\left(\pi_0, \D_{Z_{\rho^\vee}\left[x-\frac{1}{2},-x+\frac{1}{2}\right]}\left(\sigma\right)\right) \neq 0.
\]
Note that by Corollary \ref{C:co-temperedAtobeLem.4.1NStDer}, any irreducible subquotient of $\D_{Z_{\rho^\vee}\left[x-\frac{1}{2},-x+\frac{1}{2}\right]}\left(\sigma\right)$ lies in the A-packet $\AP_{\psi_{W,0}}\left(G\left(W_0\right)\right)$, where $\psi_{W,0}=\psi_W-2\rho^\vee\boxtimes S_{2x}$ and $W_0$ is a Hermitian space in the Witt tower containing $W$ of appropriate dimension.
Obviously $\left(\psi_V,\psi_W^\vee\right)$ is relevant if and only if $\left(\psi_{V,0},\psi_{W,0}^\vee\right)$ is also relevant. Therefore the desired conclusion follows from the induction hypothesis. 

\vskip 5pt

In the remaining three cases, we have $m_{\psi_V}\left(\rho \boxtimes S_{2x+1}\right)=1$. By Corollary \ref{C:co-temperedAtobeLem.4.1NStDer} we must have $x>0$, and $D_{\rho|\cdot|^{-x}}\left(\pi\right) = \pi^-$ is irreducible. Here $\pi^-$ lies in the A-packet associated to the A-parameter
\[
    \psi_V^- = \psi_V -\rho\boxtimes S_{2x+1} + \rho \boxtimes S_{2x-1}.
\]
There is a surjection $\rho|\cdot|^x\rtimes \pi^- \twoheadrightarrow \pi$. Then by the left exactness of the \(\Hom\)-functor and Proposition \ref{P:Frobenius-MackeyInequality}, we have
\[
    m\left(\pi,\sigma\right)\leq m\left(\rho|\cdot|^x\rtimes \pi^-,\sigma\right) \leq m\left(\pi^-,\D_{\rho^\vee|\cdot|^{-x-\frac{1}{2}}}\left(\sigma\right)\right) + m\left(\pi^-,\sigma\right).
\]
Again, $\D_{\rho^\vee|\cdot|^{-x-\frac{1}{2}}}\left(\sigma\right)=0$ due to the minimality of $-x$. Hence $m\left(\pi^-,\sigma\right)\neq 0$. Then by the induction hypothesis, $\left(\psi_V^-,\psi_W^\vee\right)$ is relevant. We then appeal to Lemma \ref{L:AlterDescriptionRelevance} to show the relevance of the pair $\left(\psi_V,\psi_W^\vee\right)$. Note that for any irreducible representation $\varphi$ of $WD_E$ and any $b\in\Z_{>0}$, we have
\[
    \Lambda\left(\varphi,b;\psi_V,\psi_W^\vee\right) \geq \Lambda\left(\varphi,b;\psi_V^-,\psi_W^\vee\right) \geq 0.
\]
On the other hand, for any $\left(\varphi,b\right)\neq\left(\rho, 2x\right)$, the alternating sums
\[
    \Lambda\left(\varphi,b;\psi_W^\vee,\psi_V\right) = \Lambda\left(\varphi,b;\psi_W^\vee,\psi_V^-\right)
\]
are non-negative; for $\left(\varphi,b\right) = \left(\rho, 2x\right)$, by definition
\[
    \Lambda\left(\rho, 2x;\psi_W^\vee,\psi_V\right) = m_{\psi_W^\vee}\left(\rho\boxtimes S_{2x}\right) - m_{\psi_V}\left(\rho\boxtimes S_{2x+1}\right) + \Lambda\left(\rho, 2x+2;\psi_W^\vee,\psi_V\right).
\]
Thus it suffices to show that $m_{\psi_W^\vee}\left(\rho\boxtimes S_{2x}\right)\geq 1$, i.e. $\rho^\vee\boxtimes S_{2x}\subset \psi_W$.

\vskip 5pt

\emph{Case 2:} $m_{\psi_V}\left(\rho\boxtimes S_{2x+1}\right)=1$ and $k>1$. In this case, since $\left(\psi_V^-,\psi_W^\vee\right)$ is relevant, we must have $\rho^\vee\boxtimes S_{2x+2k'} \subset \psi_W$ for some $k'\in\Z_{> 0}$. Once again by the minimality of $-x$ and an argument similar to the proof of Corollary \ref{C:TempShapeofL-parameter}, there is ``no gap'' in the $\rho^\vee$-chain of $\psi_W$ after $x$, in the sense that
\[
    \rho^\vee\boxtimes S_{2x} + \rho^\vee\boxtimes S_{2x+2} + \cdots + \rho^\vee\boxtimes S_{2x+2k'} \subset \psi_W.
\]
In particular $\rho^\vee\boxtimes S_{2x}\subset \psi_W$ as desired.

\vskip 5pt

\emph{Case 3:} $m_{\psi_V}\left(\rho\boxtimes S_{2x+1}\right)=1$, $k=1$, $x>\frac{1}{2}$ and $\rho\boxtimes S_{2x-1}\subset \psi_V$. In this case, since $\left(\psi_V^-,\psi_W^\vee\right)$ is relevant, we must have $\rho^\vee\boxtimes S_{2x+2k'}\not\subset\psi_W$ for any $k'\in\Z_{>0}$. Applying the Aubert--Zelevinsky involution to Lemma \ref{L:non-unitaryInd-0}, we obtain a  surjection $Z_{\rho}\left[x,-x+1\right]\rtimes\pi^\flat\twoheadrightarrow\pi$, where $\pi^\flat$ is an irreducible representation in the A-packet associated to the A-parameter
\[
    \psi_V^\flat = \psi_V - \rho\boxtimes S_{2x+1} - \rho\boxtimes S_{2x-1}.
\]
Then by the left exactness of the \(\Hom\) functor and Proposition \ref{P:Frobenius-MackeyInequality}, we have
\[
    m\left(\pi,\sigma\right)\leq m\left(Z_{\rho}\left[x,-x+1\right]\rtimes\pi^\flat,\sigma\right) \leq m\left(\pi^\flat, \D_{Z_{\rho^\vee}\left[x-\frac{3}{2},-x-\frac{1}{2}\right]}\left(\sigma\right)\right) + m\left(\pi^\flat, \D_{Z_{\rho^\vee}\left[x-\frac{3}{2},-x+\frac{1}{2}\right]}\left(\sigma\right)\right).
\]
It follows from the minimality of $-x$ that $D_{\rho^\vee|\cdot|^{-x-\frac{1}{2}}}\left(\sigma\right)=0$, hence also $\D_{Z_{\rho^\vee}\left[x-\frac{3}{2},-x-\frac{1}{2}\right]}\left(\sigma\right)=0$. This implies that $m\left(\pi^\flat, \D_{Z_{\rho^\vee}\left[x-\frac{3}{2},-x+\frac{1}{2}\right]}\left(\sigma\right)\right)\neq 0$, and in particular $\D_{Z_{\rho^\vee}\left[x-\frac{3}{2},-x+\frac{1}{2}\right]}\left(\sigma\right)\neq0$. Recall that $Z_{\rho^\vee}\left[x-\frac{3}{2},-x+\frac{1}{2}\right]\hookrightarrow\rho^\vee|\cdot|^{-x+\frac{1}{2}}\times\cdots\times\rho^\vee|\cdot|^{x-\frac{3}{2}}$. The non-vanishing of $\D_{Z_{\rho^\vee}\left[x-\frac{3}{2},-x+\frac{1}{2}\right]}\left(\sigma\right)$ implies that $D_{\rho^\vee|\cdot|^{-x+\frac{1}{2}}}\left(\sigma\right)\neq0$. By Corollary \ref{C:co-temperedAtobeLem.4.1NStDer}, we have $\rho^\vee\boxtimes S_{2x}\subset\psi_W$ as desired.

\vskip 5pt

\emph{Case 4:} $m_{\psi_V}\left(\rho\boxtimes S_{2x+1}\right)=1$, $k=1$, and either $x=\frac{1}{2}$ or $\rho\boxtimes S_{2x-1}\not\subset \psi_V$. This is the most delicate case. Consider the set $\Der\left(\pi,\sigma\right)\backslash\left\{\rho|\cdot|^{-x}\right\}$. Depending on this set, there are four subcases, which we treat separately.

\vskip 5pt

$\bullet$\textit{ Subcase 4.1:} $\Der\left(\pi,\sigma\right)\backslash\left\{\rho|\cdot|^{-x}\right\}$ is empty. In particular, $\sigma$ is supercuspidal, and by Corollary \ref{C:co-temperedShapeofAPara} the A-parameter $\psi_V$ is of the form
\[
    \psi_V = \psi_V^{\rho} + \left(\rho\boxtimes S_{\kappa} + \rho\boxtimes S_{\kappa+2} + \cdots + \rho\boxtimes S_{\kappa + 2\ell-2}\right) + \rho\boxtimes S_{2x+1}.
\]
Here $\psi_V^{\rho}$ satisfies the property that for any integer $a\in\Z_{>0}$, $\rho\boxtimes S_a\not\subset \psi_V^\rho$; $\kappa\in\left\{1,2\right\}$ such that $\kappa\equiv 2x+1\pmod 2$, $\ell\in\Z_{\geq 0}$ such that $\kappa+2\ell-2<2x-1$. We claim that in this subcase, $\pi^-$ is supercuspidal. We prove this claim by contradiction. Suppose on the contrary that $\pi^-$ is not supercuspidal. Then by Theorem \ref{T:AtobeJacLLCMainThm4.2} (combined with Aubert--Zelevinsky involution), we must have $x>\frac{1}{2}$, and $D_{\rho|\cdot|^{-x+1}}\left(\pi^-\right)$ is irreducible. There is a surjection
\[
    \rho|\cdot|^{x}\times \rho|\cdot|^{x-1}\rtimes D_{\rho|\cdot|^{-x+1}}\left(\pi^-\right)\twoheadrightarrow \pi.
\]
Since $D_{\rho|\cdot|^{-x+1}}\left(\pi\right)=0$, this surjection factors through $Z_\rho\left[x,x-1\right]\rtimes D_{\rho|\cdot|^{-x+1}}\left(\pi^-\right)\twoheadrightarrow \pi$. Then by Proposition \ref{P:Frobenius-MackeyInequality} we have
\begin{align*}
    m\left(\pi,\sigma\right) &\leq m\left(Z_\rho\left[x,x-1\right]\rtimes D_{\rho|\cdot|^{-x+1}}\left(\pi^-\right),\sigma\right) \\
    &\leq m\left(D_{\rho|\cdot|^{-x+1}}\left(\pi^-\right),\D_{Z_{\rho^\vee}\left[-x+\frac{1}{2},-x-\frac{1}{2}\right]}\left(\sigma\right)\right) + m\left(D_{\rho|\cdot|^{-x+1}}\left(\pi^-\right), \D_{\rho^\vee|\cdot|^{-x+\frac{1}{2}}}\left(\sigma\right)\right).
\end{align*}
However, the cuspidality of $\sigma$ implies that the right-hand side of the last inequality is zero. This contradicts our assumption $m\left(\pi,\sigma\right)\neq0$.

\vskip 5pt

With this claim at hand, now we show that $\rho^\vee\boxtimes S_{2x}\subset \psi_W$. Let $\phi_V\coloneqq \widehat{\psi}_V$ be the L-parameter of $G\left(V\right)$ obtained from $\psi_V$ by exchanging Deligne $\SL_2(\C)$ and Arthur $\SL_2(\C)$, and likewise $\phi_V^-\coloneqq\widehat{\psi}_V^-$ (resp. $\phi_W\coloneqq \widehat{\psi}_W$). Since $\pi^-$ (resp. $\sigma$) is supercuspidal, it also lies in the L-packet $\LP_{\phi_V^-}$ (resp. $\LP_{\phi_W}$). Recall that $m\left(\pi^-,\sigma\right)\neq 0$, it follows from Theorem \ref{T:TemperedGGP} that 
\[
    \pi^- = \pi\left(\phi_V^-,\eta_{\phi_V^-,\phi_W}\right) \quad \text{and} \quad \sigma = \pi\left(\phi_W,\eta_{\phi_W,\phi_V^-}\right).
\]
By \cite[Lem. 5.1]{MR4055180}, there is a short exact sequence of $G\left(V\right)$-representations
\[
    0 \lra \pi_{\rm temp} \lra \rho|\cdot|^x\rtimes \pi^- \lra \pi \lra 0,
\]
Here, we view $\calA_{\phi_V^-}$ as a subgroup of $\calA_{\phi_V}$ through the natural embedding $\calA_{\phi_V^-}\hookrightarrow\calA_{\phi_V}$, and $\pi_{\rm temp}=\pi\left(\phi_V,\eta_{\rm temp}\right)$ for a certain extension $\eta_{\rm temp}\in\Irr\left(\calA_{\phi_V}\right)$ of $\eta_{\phi_V^-,\phi_W}\in\Irr\left(\calA_{\phi_V^-}\right)$. For any basis element $\alpha_i\in\calA_{\phi_V}$ corresponding to the Jordan block $\rho_i S_{b_i}\subset \phi_V$, as long as $\rho_i S_{b_i}\not\cong\rho S_{2x+1}$, by definition we have
\[
    \eta_{\rm temp}\left(\alpha_i\right) = \eta_{\phi_V^-,\phi_W}\left(\alpha_i\right) = \eta_{\phi_V,\phi_W}\left(\alpha_i\right).
\]
Since $\phi_V$ is a discrete L-parameter in the current subcase, this computation implies that $\eta_{\rm temp}=\eta_{\phi_V,\phi_W}$. On the other hand, for any basis element $\beta_j\in\calA_{\phi_W}$ corresponding to the Jordan block $\rho'_j S_{d_j}\subset \phi_W$, we have
\[
    \frac{\eta_{\phi_W,\phi_V}\left(\beta_j\right)}{\eta_{\phi_W,\phi_V^-}\left(\beta_j\right)} = \frac{\epsilon\left(\rho S_{2x+1}\otimes\rho'_jS_{d_j}\right)}{\epsilon\left(\rho S_{2x-1}\otimes\rho'_jS_{d_j}\right)} \cdot
    \begin{cases}
      \det\left(\rho\otimes\rho'_j\right)(-1)^{d_j} \quad & \text{if }\iota = {\rm id};\\
        1 \quad & \text{if }\iota \neq {\rm id}.
    \end{cases}
\] 
According to \cite[Lem. A.4]{MR3714507}, if $\rho'_jS_{d_j}\not\cong\rho^\vee S_{2x}$, then $\eta_{\phi_W,\phi_V}\left(\beta_j\right)=\eta_{\phi_W,\phi_V^-}\left(\beta_j\right)$. Now suppose on the contrary that $\rho^\vee\boxtimes S_{2x}\not\subset \psi_W$, or equivalently, $\rho^\vee S_{2x}\not\subset\phi_W$. Then by the above computations and Theorem \ref{T:TemperedGGP}, we have $m\left(\pi_{\rm temp},\sigma\right)\neq 0$. Note that $\sigma$ is supercuspidal, the short exact sequence above gives
\[
    m\left(\rho|\cdot|^x\rtimes \pi^-,\sigma\right) = m\left(\pi_{\rm temp},\sigma\right) + m\left(\pi,\sigma\right)\geq 2.
\]
However, by the minimality of $-x$ and Proposition \ref{P:Frobenius-MackeyInequality} one knows that
\[
    m\left(\rho|\cdot|^x\rtimes \pi^-,\sigma\right) \leq m\left(\pi^-,\D_{\rho^\vee|\cdot|^{-x-\frac{1}{2}}}\left(\sigma\right)\right) + m\left(\pi^-,\sigma\right)=1.
\]
This leads us to a contradiction. Hence we must have $\rho^\vee\boxtimes S_{2x}\subset \psi_W$ as desired.

\vskip 5pt

In the remaining three subcases, the set $\Der\left(\pi,\sigma\right)\backslash\left\{\rho|\cdot|^{-x}\right\}$ is non-empty. Let $\varrho|\cdot|^{-y}$ be a supercuspidal representation in $\Der\left(\pi,\sigma\right)\backslash\left\{\rho|\cdot|^{-x}\right\}$ with the minimal exponent, where $\varrho\in\calC_{\unit}\left(\GL\right)$ and $y\in\R$ (so the exponent $e\left(\varrho|\cdot|^{-y}\right)=-y$).

\vskip 5pt

$\bullet$\textit{ Subcase 4.2:} $D_{\varrho|\cdot|^{-y}}\left(\pi\right)\neq0$. By Corollary \ref{C:co-temperedAtobeLem.4.1NStDer}, one knows that $y\in\frac{1}{2}\Z_{\geq 0}$ and $\varrho\boxtimes S_{2y+1}\subset\psi_V$. Consider the multiplicity $m_{\psi_V}\left(\varrho\boxtimes S_{2y+1}\right)$ of $\varrho\boxtimes S_{2y+1}$ in $\psi_V$. If $m_{\psi_V}\left(\varrho\boxtimes S_{2y+1}\right)\geq 2$, then we can copy the argument in Case 1 as follows. By Conjecture \ref{Conj:APacketsNProperties}(2) there is an irreducible representation $\pi'_0$ in the A-packet associated to $\psi'_{V,0}\coloneqq \psi_V - 2\varrho\boxtimes S_{2y+1}$, such that $\pi$ is a direct summand of $Z_{\varrho}\left[y,-y\right]\rtimes \pi'_0$. Combining this with Proposition \ref{P:Frobenius-MackeyInequality}, we have
\[
    m\left(\pi,\sigma\right)\leq m\left(Z_{\varrho}\left[y,-y\right]\rtimes \pi'_0,\sigma\right) \leq m\left(\pi'_0, \D_{Z_{\varrho^\vee}\left[y-\frac{1}{2},-y-\frac{1}{2}\right]}\left(\sigma\right)\right) + m\left(\pi'_0, \D_{Z_{\varrho^\vee}\left[y-\frac{1}{2},-y+\frac{1}{2}\right]}\left(\sigma\right)\right).
\]
The minimality of $-y$ implies that $\D_{Z_{\varrho^\vee}\left[y-\frac{1}{2},-y-\frac{1}{2}\right]}\left(\sigma\right)=0$. Thus $m\left(\pi'_0, \D_{Z_{\varrho^\vee}\left[y-\frac{1}{2},-y+\frac{1}{2}\right]}\left(\sigma\right)\right)\neq 0$. By Corollary \ref{C:co-temperedAtobeLem.4.1NStDer}, any irreducible subquotient of $\D_{Z_{\varrho^\vee}\left[y-\frac{1}{2},-y+\frac{1}{2}\right]}\left(\sigma\right)$ lies in the A-packet associated to $\psi'_{W,0}\coloneqq \psi_W - 2\varrho^\vee\boxtimes S_{2y}$. Then the induction hypothesis tells us that $\left(\psi'_{V,0},\left(\psi'_{W,0}\right)^\vee\right)$ is relevant, hence so is $\left(\psi_V,\psi_W^\vee\right)$.

\vskip 5pt

On the other hand if $m_{\psi_V}\left(\varrho\boxtimes S_{2y+1}\right)=1$ (so in particular $y>0$), then we can copy the argument before Case 2. Under this assumption, we know that $D_{\varrho|\cdot|^{-y}}\left(\pi\right)$ is irreducible and lies in the A-packet associated to $\psi_V^\dagger\coloneqq \psi_V - \varrho\boxtimes S_{2y+1} + \varrho\boxtimes S_{2y-1}$. There is a surjection $\varrho|\cdot|^y\rtimes D_{\varrho|\cdot|^{-y}}\left(\pi\right)\twoheadrightarrow \pi$, and this together with Proposition \ref{P:Frobenius-MackeyInequality} implies that
\[
    m\left(\pi,\sigma\right)\leq m\left(\varrho|\cdot|^y\rtimes D_{\varrho|\cdot|^{-y}}\left(\pi\right),\sigma\right)\leq m\left(D_{\varrho|\cdot|^{-y}}\left(\pi\right),\D_{\varrho^\vee|\cdot|^{-y-\frac{1}{2}}}\left(\sigma\right)\right) + m\left(D_{\varrho|\cdot|^{-y}}\left(\pi\right),\sigma\right).
\]
The minimality of $-y$ implies that $\D_{\varrho^\vee|\cdot|^{-y-\frac{1}{2}}}\left(\sigma\right)=0$. Thus $m\left(D_{\varrho|\cdot|^{-y}}\left(\pi\right),\sigma\right)\neq 0$. It follows from the induction hypothesis that $\left(\psi_V^\dagger,\psi_W^\vee\right)$ is relevant, hence
\[
    \Lambda\left(\rho,2x;\psi_W^\vee,\psi_V\right) = \Lambda\left(\rho,2x;\psi_W^\vee,\psi_V^\dagger\right)\geq 0.
\]
This already shows that $\left(\psi_V,\psi_W^\vee\right)$ is relevant, so either $\rho^\vee\boxtimes S_{2x}\subset \psi_W$ or $\rho^\vee \boxtimes S_{2x+2}\subset \psi_W$. In the latter situation, by Theorem \ref{T:AtobeJacLLCMainThm4.2} we must have $\rho^\vee\boxtimes S_{2x}\subset \psi_W$ as well, otherwise $D_{\rho^\vee|\cdot|^{-x-\frac{1}{2}}}\left(\sigma\right)\neq0$, which contradicts the minimality of $-x$.

\vskip 5pt

$\bullet$\textit{ Subcase 4.3:} $D_{\varrho|\cdot|^{-y}}\left(\sigma\right)\neq0$ and $\varrho|\cdot|^{-y}\not\cong\rho^\vee|\cdot|^{-x+\frac{1}{2}}$. Note that $D_{\varrho^\vee|\cdot|^{-y-\frac{1}{2}}}\left(\pi\right)=0$. In this subcase, the argument is almost the same as the previous subcase, except we switch the roles of $\pi$ and $\sigma$. Depending on the integer $m_{\psi_W}\left(\varrho\boxtimes S_{2y+1}\right)$, we can construct a ``more tempered'' pair of representations $\left(\pi''_0,\sigma''_0\right)$ (as in Case 1) or $\left(\pi,\sigma^-\right)$ (as in the paragraph before Case 2) with non-zero period, and then the desired conclusion $\rho^\vee\boxtimes S_{2x}\subset\psi_W$ will follow from the induction hypothesis easily.

\vskip 5pt

$\bullet$\textit{ Subcase 4.4:} $D_{\varrho|\cdot|^{-y}}\left(\sigma\right)\neq0$ and $\varrho|\cdot|^{-y}\cong\rho^\vee|\cdot|^{-x+\frac{1}{2}}$. In this subcase the desired conclusion simply follows from Corollary \ref{C:co-temperedAtobeLem.4.1NStDer}.

\vskip 5pt

The discussions above exhaust all possibilities, and now the proof of Theorem \ref{T:PeriodImpliesRelevanceGPCotempered} has been completed.


\section{Counterexample to ``period implies relevance'' for L-packets of Arthur type}\label{sec:counterexample}

When at least one of the A-parameters in the pair \(\left(\psi_V,\psi_W^{\vee}\right)\) is not co-tempered, the existence of period between the distinguished pair of representations in \(\left(\LP_{\phi_{\psi_V}}\times\LP_{\phi_{\psi_W}}\right)^{\rel}\) may not imply the relevance of \(\left(\psi_V,\psi_W^{\vee}\right)\). In Section \ref{sub:introNTGGPandMainResult} we provide a counterexample for unitary groups. Now we give another counterexample for orthogonal groups.

\vskip 5pt

\begin{Exm}\label{E:so-counterex}
Let $F$ be a $p$-adic field with uniformizer $\varpi$ and the ring of integers $\calO_F$. We choose a unit $u \in \calO_F^\times$ such that $E = F(\sqrt{u})$ is the unique unramified quadratic extension of $F$. Define $V_4 = F^4$ with standard basis $e_1, e_2, e_3, e_4$ and quadratic form
\[
    q\left(x_1, x_2, x_3, x_4\right) = x_1^2 - u x_2^2 - \varpi x_3^2 - x_4^2.
\]
Take $V_3 = \left\langle e_4\right\rangle^\perp$, with restricted form $q\big|_{V_3}$. Then $V_3$ is a non-degenerate anisotropic subspace of $V_4$, and we can identify the orthogonal group $\SO\left(V_3\right)$ as the subgroup of $\SO\left(V_4\right)$ fixing the vector $e_4$. Consider trivial representations $\mathbbm{1}_4$ and $\mathbbm{1}_3$ of $\SO\left(V_4\right)$ and $\SO\left(V_3\right)$ respectively. On the one hand we have
\[
    \Hom_{\SO\left(V_3\right)}\left(\mathbbm{1}_4,\mathbbm{1}_3\right)\neq 0.
\]  
On the other hand, both $\mathbbm{1}_4$ and $\mathbbm{1}_3$ are of Arthur type: 
\begin{itemize}
    \item $\mathbbm{1}_4$ belongs to the L-packet associated to the A-parameter $\chi_{V_4}+\mathbbm{1}\boxtimes S_3$, where $\chi_{V_4}$ is the quadratic character associated to the discriminant of $V_4$;

    \vskip 5pt

    \item since $V_3$ is anisotropic, $\mathbbm{1}_3$ is tempered, hence belongs to a tempered L-packet; one can check that the (tempered) L-parameter of $\mathbbm{1}_3$ is $S_2$.
\end{itemize}
\vskip 5pt 
Obviously the pair of A-parameters $\left(\chi_{V_4}+\mathbbm{1}\boxtimes S_3, S_2\right)$ can not be relevant. Similar to unitary group case, this does not violate Conjecture \ref{Conj:Non-temperedGGP}(4), since $\mathbbm{1}_3$ also lies in the A-packet attached to $\mathbbm{1}\boxtimes S_2$.
\end{Exm}

\vskip 5pt

We would like to propose some ideas to rectify the conjecture. In fact, Conjecture \ref{Conj:Non-temperedGGP}(2) can be regarded as an attempt to make Conjecture \ref{Conj:Non-temperedGGP}(4) more precise. From the point of view of the intersection theory of local A-packets, if $\pi\in\LP_{\phi_{\psi_V}}$ then $\psi_V=\psi^{\rm max}\left(\pi\right)$ is ``\emph{the}'' A-parameter of $\pi$, in the sense of \cite{MR4912257}. Instead of fixing Conjecture \ref{Conj:Non-temperedGGP}(2) only, it would be more interesting if one can answer the following question:

\begin{Que}
In the setting of Conjecture \ref{Conj:Non-temperedGGP}(4), how to specify these A-parameters $\psi_{V,\pi}$ and $\psi_{W,\sigma}$?     
\end{Que}

\vskip 5pt

Recall that by \cite{MR4912257}, there are some special elements in the set of A-parameters $\left\{\psi'_V\,\big|\,\pi\in\AP_{\psi'_V}\right\}$, most notably the maximal A-parameter $\psi^{\rm max}\left(\pi\right)$ and the minimal A-parameter $\psi^{\rm min}\left(\pi\right)$. A naive answer to the question above might be to always take $\psi_{V,\pi}=\psi^{\rm max}\left(\pi\right)$ or $\psi_{V,\pi}=\psi^{\rm min}\left(\pi\right)$, and likewise for $\sigma$. The failure of Conjecture \ref{Conj:Non-temperedGGP}(2) suggests that $\psi^{\rm max}$ is not the correct answer. So, what about $\psi^{\rm min}$? Note that if $\psi_V$ is co-tempered and $\pi\in\AP_{\psi_V}$, then $\psi_V = \psi^{\rm min}\left(\pi\right)$. Therefore, our proof of Conjecture \ref{Conj:Non-temperedGGP}(3) seems to be an evidence that supporting $\psi^{\rm min}$ might be the answer. However, the following example disprove this naive guess. 

\begin{Exm}
Consider the tempered L-parameters
\[
    \phi_V \coloneqq S_4 + S_6 \quad \text{and} \quad \phi_W \coloneqq S_3 + S_5
\]
of $\SO\left(V\right)$ and $\SO\left(W\right)$, where $V$ is an $11$-dimensional quadratic space and $W$ is an $8$-dimensional quadratic space respectively. Let $\pi\boxtimes\sigma$ be the distinguished representation given by the recipe of tempered GGP conjecture (Theorem \ref{T:TemperedGGP}). One can check by hand that
\[
    \psi^{\rm min}\left(\pi\right) = S_5\boxtimes S_2 \quad \text{and} \quad \psi^{\rm min}\left(\sigma\right) = S_4\boxtimes S_2,
\]
and again the pair of A-parameters $\left(\psi^{\rm min}\left(\pi\right),\psi^{\rm min}\left(\sigma\right)^\vee\right)$ is not relevant!   
\end{Exm}

It seems to us that one should choose $\psi_{V,\pi}$ according to $\sigma$ and vice versa. As of writing, we are not aware of any counterexample to Conjecture \ref{Conj:Non-temperedGGP}(4). From the global point of view, the result for split case (i.e. $\GL_n$ case) seems to suggest the following conjecture, which excludes the previous counterexamples above, but its truth is not known.

\begin{Conj}
    If \(\left(\pi,\sigma\right)\) is a pair of representation in \(\left(\Pi_{\phi_{\psi_V}}\times\Pi_{\phi_{\psi_W}}\right)^{\rel}\) such that \(\pi\boxtimes \sigma\) occurs on the split form of \(G\left(V\right)\times G\left(W\right)\) and \(m(\pi,\sigma)=1\), then \(\left(\psi_V,\psi^{\vee}_W\right)\) is relevant.
\end{Conj}

\vskip 5pt

Finally, we consider the notion of \emph{strong Ext relevance} defined by Qadri in his paper \cite[Def. 1.1]{MR4846726}. The definition is recorded below for reference.
\begin{Def}[{\cite[Def. 1.1]{MR4846726}}]
We say the pair of A-parameters \(\left(\psi_V,\psi_W^\vee\right)\) is strong Ext relevant if $\psi_V$ and $\psi_W^\vee$ admit decompositions of the form
    \begin{align*}
        \psi_V&=\sum^r_{i=1}\rho_i S_{c_i}\boxtimes S_{a_i}+ \sum^{r+s}_{i=r+1}\rho_i S_{c_i}\boxtimes S_{b_i-1}+\sum^t_{i=1}\rho'_i S_{f_i}\boxtimes S_{d_i}+ \sum^{t+u}_{i=t+1}\rho'_i S_{e_i-1}\boxtimes S_{f_i}\\
        \psi_W&=\sum^r_{i=1}\rho_i S_{c_i}\boxtimes S_{a_i-1}+ \sum^{r+s}_{i=r+1}\rho_i S_{c_i}\boxtimes S_{b_i}+\sum^t_{i=1}\rho'_i S_{d_i-1}\boxtimes S_{f_i}+ \sum^{t+u}_{i=t+1}\rho'_i S_{f_i}\boxtimes S_{e_i}.
    \end{align*}
\end{Def}

\vskip 5pt

\begin{Conj}
In the setting of Conjecture \ref{Conj:Non-temperedGGP}(2), if \(m\left(\pi,\sigma\right)=1\), then \(\left(\psi_V,\psi^{\vee}_W\right)\) is strong Ext relevant.
\end{Conj}

\vskip 5pt

All counterexamples above satisfy this new modified conjecture. It is also worth noting that, by \cite{MR4661536, MR4912257}, there are several operators that ``act transitively'' on the set of A-parameters $\left\{\psi'_V\,\big|\,\pi\in\AP_{\psi'_V}\right\}$, one of which is the ``\emph{partial Aubert--Zelevinsky dual}''. Applying this operator to the notion of relevance, one obtains this notion of strong Ext relevance. It is interesting to ask that whether the union-intersection operator will also play a role in determining $\psi_{V,\pi}$ and $\psi_{W,\sigma}$? We do not know the answer. 

\vskip 10pt

\appendix

\section{Counting Jordan blocks via Hom spaces} 
\label{App:ProofTechnicalCountingMulti}

The main purpose of this appendix is to prove a technical lemma, which supplies the proof of Case 1, Subcase 1.1 of Theorem \ref{T:ReductionToGoodParity}(2) in Section \ref{sub:proof_of_theorem_t:reductiontogoodparity_the_only_if_part}.

\vskip 5pt

Let $\rho$ be an irreducible unitary supercuspidal representation of $\GL_d\left(E\right)$, and $a,b\in\Z_{>0}$. We set
\[
    S\coloneqq u_\rho\left(a,b\right), \quad \text{and} \quad \begin{cases}
        A\coloneqq \frac{a+b}{2}-1;\\
        B\coloneqq \frac{a-b}{2}.
    \end{cases}
\]
Let $\tau$ be an irreducible Arthur type representation of $\GL_N\left(E\right)$, with a decomposition
\begin{equation}\label{E:DecompTauApp}
    \tau = u_{\rho_1}\left(a_1,b_1\right)\times u_{\rho_2}\left(a_2,b_2\right)\times\cdots\times u_{\rho_r}\left(a_r,b_r\right).
\end{equation}
We define $A_i,B_i\in\frac{1}{2}\Z$ similarly. Recall that for a representation $\pi$ of a general linear group, we set $\widetilde{\pi}\coloneqq \prescript{\iota}{}{\pi}^\vee$. We denote by $n\coloneqq \left|\left\{1\leq i\leq r\,\big|\,\rho_i\in\left\{\rho,\widetilde{\rho}\right\} \text{ and } \left(a_i,b_i\right)=\left(a,b\right)\right\}\right|$. The main result of this appendix is the following.

\begin{Prop}\label{P:TechnicalCountingMultiApp}
Suppose that for all $1\leq i\leq r$ such that $\rho_i\in\left\{\rho,\widetilde{\rho}\right\}$ and $a_i+b_i\equiv a+b \pmod2$, we have $a_i+b_i\leq a+b$ and $b_i\leq b$. Let $m\in\Z_{>0}$ be such that $mdab\leq N$. If
\[
    \Hom_{\GL_N\left(E\right)}\left(\tau, \widetilde{L}^-_1\times\widetilde{L}^-_2\times\cdots\times\widetilde{L}^-_m\times\Xi\times L^+_m\times\cdots\times L_2^+\times L_1^+\right)\neq 0
\]
for some irreducible constituents $L^-_1\boxtimes L^+_1, L^-_2\boxtimes L^+_2, \cdots, L^-_m\boxtimes L^+_m$ of $m^*\left(S\right)$ and some irreducible representation $\Xi$ of $\GL_{N-mdab}\left(E\right)$, then $m\leq n$.    
\end{Prop}

\vskip 5pt

\subsection{Truncated categories and intertwining isomorphisms} 
\label{sub:truncated_category_of_finite_length_representations}

We need some preparations for the proof of Proposition \ref{P:TechnicalCountingMultiApp}. Firstly, we would like to transform the Hom space considered in Proposition \ref{P:TechnicalCountingMultiApp} into another form. To do this, it will be convenient to restrict ourselves to a smaller category. For each $k\in\Z_{\geq 0}$, let $\calR\left(\GL_k\right)$ be the category of \emph{finite length} smooth representations of $\GL_k\left(E\right)$. We consider its full subcategory $\calR\left(\GL_k\right)_{\rho,A}$ consisting of all finite length smooth representations $\Pi$ with the following property: for any irreducible subquotient $\pi$ of $\Pi$ and any $\varrho\in\csupp\left(\pi\right)$, one has
\[
    \varrho\cong\rho|\cdot|^{x}
\]
for some $x\in\left\{A,A-1,\cdots,-A\right\}$. This ``truncated'' category $\calR\left(\GL_k\right)_{\rho,A}$ enjoys an useful property, which we state as below.

\begin{Thm}\label{T:ProdFunctorFullyFaithfulApp}
Let $S'\coloneqq u_\rho\left(a',b'\right)$ for some $a',b'\in\Z_{>0}$ such that $a'+b'=a+b$. Then the parabolic induction functors
\[
    \begin{aligned}
       S'\times-: \calR\left(\GL_k\right)_{\rho,A} &\lra \calR\left(\GL_{k+da'b'}\right)_{\rho,A}\\
       \Pi &\longmapsto S'\times\Pi
    \end{aligned},
    \qquad
    \begin{aligned}
       -\times S': \calR\left(\GL_k\right)_{\rho,A} &\lra \calR\left(\GL_{k+da'b'}\right)_{\rho,A}\\
       \Pi &\longmapsto \Pi\times S'
    \end{aligned}
\] 
are fully-faithful. Moreover, for any irreducible representation $\pi$ in the category $\calR\left(\GL_k\right)_{\rho,A}$, the parabolically induced representations $S'\times\pi$ and $\pi\times S'$ are still irreducible.  
\end{Thm}

\vskip 5pt

\begin{proof}
The assertions for the functor $S'\times-$ follow from \cite[Lem. 8.5, Thm. 9.1]{MR4373242}. 
As for the functor $-\times S'$, the assertions can be deduced from the corresponding conclusions for the functor $S'\times-$ by using the MVW-contragredient trick.
 
\end{proof}

\vskip 5pt

Let $L^-\boxtimes L^+$ be an irreducible constituent of $m^*\left(S\right)$. By the result of Kret--Lapid \cite[Cor. 2.2]{MR2996769}, Jacquet modules of $S$ are semi-simple. Hence $L^-\boxtimes L^+$ is a direct summand of the Jacquet module of $S$ along a maximal parabolic subgroup (say, the one with Levi subgroup $\GL_{D^-}\left(E\right)\times \GL_{D^+}\left(E\right)$ for a pair of non-negative integers $\left(D^-,D^+\right)$ such that $D^-+D^+=dab$). By Frobenius reciprocity and second adjointness, we obtain the canonical inclusion and surjection
\[
    \mathfrak{i}: S\hookrightarrow L^-\times L^+ \quad \text{and} \quad \mathfrak{e}: L^+\times L^-\twoheadrightarrow S.
\]
By using these two maps, we further obtain two intertwining maps:
\[
    \mathfrak{u}\coloneqq \left({\rm id}_{L^-}\times\mathfrak{e}\right)\circ\left(\mathfrak{i}\times {\rm id}_{L^-}\right): S\times L^- \xlongrightarrow{\mathfrak{i}\times {\rm id}_{L^-}} L^-\times L^+\times L^- \xlongrightarrow{{\rm id}_{L^-}\times\mathfrak{e}} L^-\times S,
\]
and
\[
    \mathfrak{v}\coloneqq \left(\mathfrak{e}\times {{\rm id}_{L^+}}\right)\circ\left({\rm id}_{L^+}\times\mathfrak{i}\right): L^+\times S \xlongrightarrow{{\rm id}_{L^+}\times\mathfrak{i}} L^+\times L^-\times L^+ \xlongrightarrow{\mathfrak{e}\times {{\rm id}_{L^+}}} S\times L^+.
\]

\begin{Prop}\label{P:IntMapsIsoApp}
These two intertwining maps $\mathfrak{u}$ and $\mathfrak{v}$ are indeed isomorphisms.    
\end{Prop}

\vskip 5pt

In order to prove this proposition, we need the so-called ``\emph{sandwich lemma}'' as below, which was originally proved by Kang--Kashiwara--Kim--Oh \cite[Lem. 3.1]{MR3314831} for quiver Hecke algebras, and translated into the form we stated via the exact monoidal equivalence of \cite[Thm. 3.10, Prop. 3.11]{MR4549143}.

\vskip 5pt

\begin{Lem}[Sandwich lemma]\label{L:OverlapLemmaApp}
Let $M_i\in\calR\left(\GL_{k_i}\right)_{\rho,A}$ for $i=1,2,3$, and suppose that $X\subset M_1\times M_2$ and $Y\subset M_2\times M_3$ are subrepresentations satisfying
\[
    X\times M_3 \subset M_1\times Y.
\]
Then there exists a subrepresentation $Z\subset M_2$, such that $X\subset M_1\times Z$ and $Z\times M_3\subset Y$.
\end{Lem}

\vskip 5pt

\begin{proof}[Proof of Proposition \ref{P:IntMapsIsoApp}]
We shall only prove that $\mathfrak{u}$ is an isomorphism. The assertion for $\mathfrak{v}$ can be proved similarly.
Thanks to Theorem \ref{T:ProdFunctorFullyFaithfulApp}, we know that both $S\times L^-$ and $L^-\times S$ are irreducible. Therefore to show that $\mathfrak{u}$ is an isomorphism, it suffices to show that $\mathfrak{u}\neq0$. Suppose on the contrary that $\mathfrak{u}\coloneqq \left({\rm id}_{L^-}\times\mathfrak{e}\right)\circ\left(\mathfrak{i}\times {\rm id}_{L^-}\right)=0$. Let $K\subset L^+\times L^-$ be the kernel of the surjection $\mathfrak{e}$. Then we have
\[
    S\times L^-\subset L^-\times K.
\]
By the sandwich lemma (Lemma \ref{L:OverlapLemmaApp}), there exists a subrepresentation $Z\subset L^+$, such that $S\subset L^-\times Z$ and $Z\times L^-\subset K$. However, since $L^+$ is irreducible, the only subrepresentations of $L^+$ are $0$ and itself, and obviously neither meets these requirements. Hence we must have $\mathfrak{u}\neq 0$ as desired.
    
\end{proof}

\vskip 5pt

Now we can define a bijection of Hom spaces by using these four maps $\mathfrak{i}$, $\mathfrak{e}$, $\mathfrak{u}$ and $\mathfrak{v}$. Let $k\geq D^+$ be an integer, $X\in\calR\left(\GL_{k}\right)_{\rho,A}$ and $Y\in\calR\left(\GL_{k-D^+}\right)_{\rho,A}$. Then we set the map
\[
    \Phi: \Hom_{\GL_k\left(E\right)}\left(X,Y\times L^+\right) \lra \Hom_{\GL_{k+D^-}\left(E\right)}\left(X\times L^-,Y\times S\right)
\]
as follows: for $f\in\Hom_{\GL_k\left(E\right)}\left(X,Y\times L^+\right)$, let
\[
    \Phi\left(f\right) \coloneqq \left({\rm id}_Y\times\mathfrak{e}\right)\circ \left(f\times {\rm id}_{L^-}\right): X\times L^- \xlongrightarrow{f\times {\rm id}_{L^-}} Y\times L^+\times L^- \xlongrightarrow{{\rm id}_Y\times\mathfrak{e}} Y\times S.
\]
Conversely, for $g\in\Hom_{\GL_{k+D^-}\left(E\right)}\left(X\times L^-,Y\times S\right)$ we consider the composition 
\[
    X\times S \xlongrightarrow{{\rm id}_X\times\mathfrak{i}} X\times L^-\times L^+ \xlongrightarrow{g\times{\rm id}_{L^+}} Y\times S\times L^+  \xlongrightarrow{{\rm id}_Y\times \mathfrak{v}^{-1}} Y\times L^+ \times S.
\]
This gives an element in the space $\Hom_{\GL_{k+dab}\left(E\right)}\left(X\times S, Y\times L^+\times S\right)$. Note that $X, Y\times L^+\in\calR\left(\GL_k\right)_{\rho,A}$. By the fully-faithfulness of the functor $-\times S$, we have a bijection of the Hom spaces
\[
    \Hom_{\GL_{k+dab}\left(E\right)}\left(X\times S, Y\times L^+\times S\right) \xlongrightarrow{\sim}\Hom_{\GL_k\left(E\right)}\left(X,Y\times L^+\right).
\]
We shall denote by $\Psi\left(g\right)$ the image of the composition $\left({\rm id}_Y\times \mathfrak{v}^{-1}\right)\circ\left(g\times{\rm id}_{L^+}\right)\circ\left({\rm id}_X\times\mathfrak{i}\right)$ in $\Hom_{\GL_k\left(E\right)}\left(X,Y\times L^+\right)$. This assignment $g\mapsto \Psi\left(g\right)$ gives a reverse map
\[
    \Psi: \Hom_{\GL_{k+D^-}\left(E\right)}\left(X\times L^-,Y\times S\right) \lra \Hom_{\GL_k\left(E\right)}\left(X,Y\times L^+\right).
\]

\begin{Prop}\label{P:ConvertHomSpaceApp}
In the context of above discussion, the maps $\Phi$ and $\Psi$ are bijections, and they are mutually inverse.
\end{Prop}

\vskip 5pt

\begin{proof}
We first show that $\Psi\left(\Phi\left(f\right)\right)=f$ for any $f\in\Hom_{\GL_k\left(E\right)}\left(X,Y\times L^+\right)$. By Theorem \ref{T:ProdFunctorFullyFaithfulApp}, this boils down to checking that the diagram below is commutative, which is straightforward; we suppress the details.
\[
\begin{tikzcd}[row sep=2.2em]
  & {X\times L^-\times L^+} 
    \arrow[dr, "{f\times\mathrm{id}_{L^-}\times\mathrm{id}_{L^+}}"{pos=0.38, inner sep=1.5pt}] 
    \arrow[drr, "{\Phi(f)\times\mathrm{id}_{L^+}}", bend left=28] 
  &[-1em] &[3em] \\
{X\times S} 
  \arrow[ur, "{\mathrm{id}_X\times \mathfrak{i}}"{pos=0.55}] 
  \arrow[dr, "{f\times \mathrm{id}_S}"'{pos=0.55}] 
  & & {Y\times L^+\times L^-\times L^+} 
    \arrow[r, "{\mathrm{id}_Y\times\mathfrak{e}\times\mathrm{id}_{L^+}}"] 
  & {Y\times S\times L^+} \\
  & {Y\times L^+\times S} 
    \arrow[ur, "{\mathrm{id}_Y\times\mathrm{id}_{L^+}\times\mathfrak{i}}"'{pos=0.38, inner sep=1.5pt}] 
    \arrow[urr, "{\mathrm{id}_Y\times \mathfrak{v}}"', bend right=28] 
  & &
\end{tikzcd}
\]

\vskip 5pt

Next, we show that $\Psi$ is injective. Note that by the previous paragraph we have $\left(\Psi\circ\Phi\circ\Psi\right)\left(g\right)=\Psi\left(g\right)$ for any $g\in\Hom_{\GL_{k+D^-}\left(E\right)}\left(X\times L^-,Y\times S\right)$. Once the injectivity of $\Psi$ is established, we can deduce that $\Phi$ and $\Psi$ are mutually inverse, and the desired conclusion follows. Let $g\in\Hom_{\GL_{k+D^-}\left(E\right)}\left(X\times L^-,Y\times S\right)$ be such that $\Psi\left(g\right)=0$. Recall that by the definition
\[
    \Psi\left(g\right)\times\id_S = \left({\rm id}_Y\times \mathfrak{v}^{-1}\right)\circ\left(g\times{\rm id}_{L^+}\right)\circ\left({\rm id}_X\times\mathfrak{i}\right).
\]
Since ${\rm id}_Y\times \mathfrak{v}^{-1}$ is an isomorphism, we have $\left(g\times{\rm id}_{L^+}\right)\circ\left({\rm id}_X\times\mathfrak{i}\right)=0$. Consider the commutative diagram
\[
\begin{tikzcd}[row sep=4.0em, column sep=4.5em]
{X\times S\times L^-} 
  \arrow[r, "{\id_X\times\mathfrak{i}\times\id_{L^-}}"] 
  \arrow[dr, "{\id_X\times\mathfrak{u}}"'{pos=0.45}, bend right=20] 
& {X\times L^-\times L^+\times L^-} 
  \arrow[r, "{g\times \id_{L^+}\times\id_{L^-}}"] 
  \arrow[d, "{\id_X\times\id_{L^-}\times\mathfrak{e}}"] 
& {Y\times S\times L^+\times L^-} 
  \arrow[d, "{\id_Y\times\id_S\times\mathfrak{e}}"] \\
& {X\times L^-\times S} 
  \arrow[r, "{g\times\id_S}"'] 
& {Y\times S\times S}
\end{tikzcd}
\]
By assumption, the composition of top horizontal maps is vanishing. This implies that $\left(g\times\id_S\right)\circ\left(\id_X\times \mathfrak{u}\right)=0$. Recall that by Proposition \ref{P:IntMapsIsoApp} $\left(\id_X\times \mathfrak{u}\right)$ is an isomorphism. It follows that $g\times\id_S=0$, and Theorem \ref{T:ProdFunctorFullyFaithfulApp} then implies that $g=0$. Hence $\Psi$ is injective.

\end{proof}

\vskip 5pt


\subsection{Langlands classification versus Zelevinsky classification} 
\label{sub:langlands_classification_versus_zelevinsky_classification}

Secondly, for the purpose of counting Jordan blocks, it seems to be more convenient to use the \emph{Zelevinsky classification} rather than the Langlands classification introduced in Theorem \ref{T:LanglandsClass4GL}. Recall that a \emph{multi-segment} is a finite (multi-)set of segments
\[
    \mathfrak{m} = \left\{\left[x_1,y_1\right]_{\rho_1}, \left[x_2,y_2\right]_{\rho_2}, \cdots, \left[x_r,y_r\right]_{\rho_r}\right\}.
\]
We set $|\mathfrak{m}|\coloneqq \sum_i d_{\rho_i}\left(x_i-y_i+1\right)$, where $d_{\rho_i}\coloneqq\dim\rho_i$. Suppose that the segments in $\mathfrak{m}$ are ordered so that $x_1+y_1\leq x_2+y_2\leq \cdots\leq x_r+y_r$. In Theorem \ref{T:LanglandsClass4GL}, we have introduced the representation
\[
    {\rm L}\left(\mathfrak{m}\right) \coloneqq \LS\left(\Delta_{\rho_1}\left[x_1,y_1\right], \Delta_{\rho_2}\left[x_2,y_2\right], \cdots, \Delta_{\rho_r}\left[x_r,y_r\right]\right).
\]
Zelevinsky also attached another irreducible representation ${\rm Z}\left(\mathfrak{m}\right)$ to this multi-segment $\mathfrak{m}$, defined as the unique irreducible quotient of the parabolically induced representation
\[
    Z_{\rho_1}\left[x_1,y_1\right]\times Z_{\rho_2}\left[x_2,y_2\right]\times \cdots\times Z_{\rho_r}\left[x_r,y_r\right].
\]
The Aubert--Zelevinsky involution induces an involution $\mathfrak{m}\mapsto\mathfrak{m}^t$ on the set of all multi-segment, such that:
\begin{itemize}
    \item the Aubert--Zelevinsky dual of ${\rm L}\left(\mathfrak{m}\right)$ is ${\rm L}\left(\mathfrak{m}^t\right)$; similarly, the Aubert--Zelevinsky dual of ${\rm Z}\left(\mathfrak{m}\right)$ is ${\rm Z}\left(\mathfrak{m}^t\right)$;

    \vskip 5pt

    \item for any multi-segment $\mathfrak{m}$, we have ${\rm L}\left(\mathfrak{m}\right)\cong{\rm Z}\left(\mathfrak{m}^t\right)$.
\end{itemize}

\vskip 5pt

In \cite{MR863522}, M{\oe}glin and Waldspurger gave an algorithm to compute $\mathfrak{m}^t$ in terms of $\mathfrak{m}$, which we shall now describe. For simplicity, from now on we only consider those multi-segment $\mathfrak{m}$ such that ${\rm L}\left(\mathfrak{m}\right), {\rm Z}\left(\mathfrak{m}\right)\in\calR\left(\GL_{|\mathfrak{m}|}\right)_{\rho,A}$. If $\mathfrak{m}$ is such a multi-segment, then for any segment $\left[x,y\right]_{\rho'}\in\mathfrak{m}$, one has $\left[x,y\right]_{\rho'} \subset \left\{\rho|\cdot|^A, \rho|\cdot|^{A-1},\cdots,\rho|\cdot|^{-A}\right\}$. There are three steps in the algorithm:

\vskip 5pt

\begin{itemize}
    \item[-] \emph{Step 1:} Let $y_0$ be the minimal half-integer such that $\left[x,y_0\right]_\rho\in\mathfrak{m}$ for some $x\in\frac{1}{2}\Z$; and let $x_0$ be the minimal half-integer with the property that $\left[x_0,y_0\right]_\rho\in\mathfrak{m}$.

    \vskip 5pt

    \item[-] \emph{Step 2:} Suppose that $\left[x_i,y_i\right]_\rho\in\mathfrak{m}$ has been defined. Let $x_{i+1}$ be the minimal half-integer such that $x_{i+1}>x_i$ and $\left[x_{i+1},y_i+1\right]_\rho\in\mathfrak{m}$. If such $x_{i+1}$ does not exists, we terminate this step; if such $x_{i+1}$ exists, we set $y_{i+1}\coloneqq y_i+1$ and repeat this step.

    \vskip 5pt

    \item[-] \emph{Step 3:} Suppose that Step 2 terminates after $\ell$ rounds, and we obtain a sequence of segments $\left[x_0,y_0\right]_\rho, \left[x_1,y_1\right]_\rho, \cdots, \left[x_\ell,y_\ell\right]_\rho$. Let
    \[
        \mathfrak{m}_-\coloneqq \left(\mathfrak{m}\backslash\left\{\left[x_0,y_0\right]_\rho, \cdots, \left[x_\ell,y_\ell\right]_\rho\right\}\right)\sqcup\left\{\left[x_0,y_0+1\right]_\rho, \cdots, \left[x_\ell,y_\ell+1\right]_\rho\right\}.
    \]
    Then we have
    \[
        \mathfrak{m}^t = \mathfrak{m}_-^t \sqcup \left\{\left[y_0+\ell,y_0\right]_\rho\right\}.
    \]
    We then repeat Steps 1-3 to compute $\mathfrak{m}_-^t$, hence compute $\mathfrak{m}^t$ inductively.
\end{itemize}
\vskip 5pt
Applying this algorithm, we can obtain some combinatorial information on the Zelevinsky data, which will be crucial in later argument.

\begin{Lem}\label{L:CountSegmentInZeleData}
Let $\mathfrak{m}$ be a multi-segment such that ${\rm L}\left(\mathfrak{m}\right)\in\calR\left(\GL_{|\mathfrak{m}|}\right)_{\rho,A}$. Suppose that $\left[x,-A\right]_\rho\notin\mathfrak{m}$ for any $x<B$. Then
\[
    \left|\left\{\left[x,y\right]_\rho\in\mathfrak{m}^t\,\big|\,x\geq -B \text{ and } y=-A\right\}\right| \leq \left|\left\{\left[x,y\right]_\rho\in\mathfrak{m}\,\big|\,x=A\right\}\right|.
\]
Moreover, if $\left[x,y\right]_\rho\in\mathfrak{m}^t$ such that $x\geq -B$ and $y=-A$, then we must have $x=-B$. 
\end{Lem}

\vskip 5pt

\begin{proof}
Consider a round of the above algorithm that produces a segment $\left[x,-A\right]_\rho\in\mathfrak{m}^t$ for some $x\geq -B$. Suppose the sequence of segments involved in Step 2 of this round is 
\[
    \left[x_0,y_0\right]_\rho, \left[x_1,y_1\right]_\rho, \cdots, \left[x_\ell,y_\ell\right]_\rho,
\]
where $y_0=-A$ and $\ell = x+A$. Note that in each round of the algorithm, we replace some segments in $\mathfrak{m}$ (or its substitutions) by some shorter ones, each time we delete the supercuspidal representations with minimal exponents in those segments. Since the representation $\rho|\cdot|^{-A}$ already has the minimal possible exponent, we must have $\left[x_0,y_0\right]_\rho\in\mathfrak{m}$. Then the assumption of the lemma implies that $x_0\geq B$. It follows that
\[
    A \geq x_\ell\geq x_0 + \ell\geq A+B+x.
\]
This inequality forces that $x=-B$, $\ell=b-1$ and $x_i=B+i$ for all $0\leq i\leq b-1$. In particular, the last segment in the sequence is $\left[x_\ell,y_\ell\right]_\rho = \left[A,-B\right]_\rho$. After this round is completed, the segment $\left[x_\ell,y_\ell\right]_\rho$ will be shortened to $\left[A,-B+1\right]_\rho$. As one can see from our argument, this shortened segment will never be involved in any round of the algorithm producing $\left[x',-A\right]_\rho\in\mathfrak{m}^t$ for some $x'\geq -B$. Therefore the desired inequality holds.

\end{proof}

\vskip 5pt


\subsection{Proof of Proposition \ref{P:TechnicalCountingMultiApp}: I. (conjugate) self-dual case} 
\label{sub:proof_of_proposition_p:technicalcountingmultiapp_i_conjugate_self_dual_case}

Now we can start to prove Proposition \ref{P:TechnicalCountingMultiApp}. We first consider the case that $\rho\cong\widetilde{\rho}$ is (conjugate) self-dual. By Bernstein--Zelevinsky's block decomposition, we can write $\tau$ and $\Xi$ as
\[
    \tau = \tau_\rho\times\tau^\rho \quad \text{and} \quad \Xi = \Xi_\rho\times \Xi^\rho,
\]
where $\tau_\rho$ (resp. $\tau^\rho$) is an irreducible representation with cuspidal support contained in (resp. disjoint from) the cuspidal line $\rho|\cdot|^{A+\Z}\coloneqq\left\{\rho|\cdot|^{A+x}\,\big|\,x\in\Z\right\}$, and similarly for $\Xi_\rho$ and $\Xi^\rho$. Then we have
\[
    \Hom_{\GL_N\left(E\right)}\left(\tau, \widetilde{L}^-_1\times\widetilde{L}^-_2\times\cdots\times\widetilde{L}^-_m\times\Xi\times L^+_m\times\cdots\times L_2^+\times L_1^+\right)
\]
\[
    \cong \Hom_{\GL_N\left(E\right)}\left(\tau_\rho\times\tau^\rho, \left(\widetilde{L}^-_1\times\widetilde{L}^-_2\times\cdots\times\widetilde{L}^-_m\times\Xi_\rho\times L^+_m\times\cdots\times L_2^+\times L_1^+\right)\times\Xi^\rho\right).
\]
For the latter Hom space to be non-zero, $\tau_\rho$ and $\widetilde{L}^-_1\times\cdots\times\widetilde{L}^-_m\times\Xi_\rho\times L^+_m\times\cdots\times L_1^+$ must be representations of the same general linear group (say, $\GL_{N'}\left(E\right)$), and 
\[
    \Hom_{\GL_{N'}\left(E\right)}\left(\tau_\rho, \widetilde{L}^-_1\times\widetilde{L}^-_2\times\cdots\times\widetilde{L}^-_m\times\Xi_\rho\times L^+_m\times\cdots\times L_2^+\times L_1^+\right)\neq0.
\]
By the assumption of Proposition \ref{P:TechnicalCountingMultiApp}, $\csupp\left(\tau_\rho\right)\subset \left\{\rho|\cdot|^A, \rho|\cdot|^{A-1},\cdots,\rho|\cdot|^{-A}\right\}$. Therefore we must have $\csupp\left(\Xi_\rho\right)\subset \left\{\rho|\cdot|^A, \rho|\cdot|^{A-1},\cdots,\rho|\cdot|^{-A}\right\}$. Without loss of generality, we can assume that $\tau=\tau_\rho$ and $\Xi=\Xi_\rho$. Under this assumption, both $\tau$ and $\widetilde{L}^-_1\times\cdots\times\widetilde{L}^-_m\times\Xi\times L^+_m\times\cdots\times L_1^+$ belong to the category $\calR\left(\GL_N\right)_{\rho,A}$. Hence in the decomposition (\ref{E:DecompTauApp}) of $\tau$, we have $\rho_i\cong \rho$ for all $1\leq i\leq r$. We may re-order the decomposition (\ref{E:DecompTauApp}) such that
\[
    a+b=a_1+b_1=\cdots=a_k+b_k>a_{k+1}+b_{k+1}\geq\cdots\geq a_r+b_r.
\]
Suppose that for $1\leq j\leq m$, $L^-_j$ is an irreducible ladder representation of $\GL_{\ell^-_j}\left(E\right)$. Let $\ell^-\coloneqq \ell^-_1+\ell^-_2+\cdots+\ell^-_m$, and $\Sigma\coloneqq\widetilde{L}^-_1\times\cdots\times\widetilde{L}^-_m\times\Xi$. Then by Proposition \ref{P:ConvertHomSpaceApp} we have
\begin{align*}
  \Hom_{\GL_N\left(E\right)}\left(\tau, \Sigma\times L^+_m\times\cdots\times L_2^+\times L_1^+\right)&\cong \Hom_{\GL_{N+\ell^-_1}\left(E\right)}\left(\tau\times L^-_1, \Sigma\times L^+_m\times\cdots\times L_2^+\times S\right) \\
  &\cong\Hom_{\GL_{N+\ell^-_1}\left(E\right)}\left(\tau\times L^-_1, \Sigma\times S\times L^+_m\times\cdots\times L_2^+\right)\\
  &\cdots\\
  &\cong\Hom_{\GL_{N+\ell^-}\left(E\right)}\left(\tau\times L^-_1\times\cdots\times L^-_m, \Sigma\times S^{\times m}\right).
\end{align*}
Here in the second isomorphism, we make use of the fact that $L_j^+\times S\cong S\times L_j^+$ is irreducible.

\vskip 5pt

We would like to estimate the multiplicities of the supercuspidal representation $\rho|\cdot|^{-A}$ in the cuspidal supports of certain iterative derivatives of irreducible subquotients of
\[
    \tau\times L^-_1\times\cdots\times L^-_m\cong \left(\bigtimes_{i=1}^k u_\rho\left(a_i,b_i\right)\right) \times \left(\left(\bigtimes_{i=k+1}^r u_\rho\left(a_i,b_i\right)\right)\times L^-_1\times\cdots\times L^-_m\right).
\]
By Theorem \ref{T:ProdFunctorFullyFaithfulApp}, any irreducible subquotient of the above representation must be of the form $\left(\bigtimes_{i\leq k} u_\rho\left(a_i,b_i\right)\right)\times X$ for an irreducible subquotient $X$ of $\left(\bigtimes_{i>k} u_\rho\left(a_i,b_i\right)\right)\times L^-_1\times\cdots\times L^-_m$. Let $\mathfrak{m}_X$ (resp. $\mathfrak{m}_{>k}$, $\mathfrak{m}_{L^-_j}$) be the multi-segment such that ${\rm L}\left(\mathfrak{m}_X\right)\cong X$ (resp. ${\rm L}\left(\mathfrak{m}_{>k}\right)\cong \bigtimes_{i>k} u_\rho\left(a_i,b_i\right)$, ${\rm L}\left(\mathfrak{m}_{L^-_j}\right)\cong L^-_j$). By \cite[Thm. 7.1]{ZeleII}, we know that $\mathfrak{m}_X$ is obtained from $\mathfrak{m}_{>k}\sqcup\mathfrak{m}_{L^-_1}\sqcup\mathfrak{m}_{L^-_2}\sqcup\cdots\sqcup\mathfrak{m}_{L^-_m}$ by applying finitely many times union-intersection of linked segments. Since
\[
    \csupp\left(\bigtimes_{i>k} u_\rho\left(a_i,b_i\right)\right)\subset\left\{\rho|\cdot|^{A-1},\cdots,\rho|\cdot|^{1-A}\right\},
\]
for any $\left[x,-A\right]_\rho\in\mathfrak{m}_{>k}\sqcup\mathfrak{m}_{L^-_1}\sqcup\mathfrak{m}_{L^-_2}\sqcup\cdots\sqcup\mathfrak{m}_{L^-_m}$, there is a unique $j\in\left\{1,2,\cdots,m\right\}$ such that $\left[x,-A\right]_\rho\in\mathfrak{m}_{L^-_j}$. The result of Kret--Lapid (Lemma \ref{L:KretLapidJacFormula}) then implies that $x=B$. Starting from $\mathfrak{m}_{>k}\sqcup\mathfrak{m}_{L^-_1}\sqcup\mathfrak{m}_{L^-_2}\sqcup\cdots\sqcup\mathfrak{m}_{L^-_m}$, each time one performs the union-intersection procedure of two linked segments, a segment $\left[x,-A\right]_\rho$ with $x\geq B$ can only be replaced by $\left[x',-A\right]_\rho$ with $x'\geq x$; the number of such segments remains the same. Hence the irreducible representation $X$ satisfies the assumption in Lemma \ref{L:CountSegmentInZeleData}. It follows that
\begin{align*}
    \left|\left\{\left[x,y\right]_\rho\in\mathfrak{m}_X^t\,\big|\,x\geq -B \text{ and } y=-A\right\}\right| &\leq \left|\left\{\left[x,y\right]_\rho\in\mathfrak{m}_X\,\big|\,x=A\right\}\right|\\
    & = \sum_{j=1}^m\left|\left\{\left[x,y\right]_\rho\in\mathfrak{m}_{L^-_j}\,\big|\,x=A\right\}\right|.
\end{align*}
Consider the iterative highest Bernstein--Zelevinsky derivatives
\[
    D_{\rm max}^{b-1}\left(\bigtimes_{i=1}^k u_\rho\left(a_i,b_i\right)\times X\right)\coloneqq\Big(\cdots\Big(\Big(\Big(\bigtimes_{i=1}^k u_\rho\left(a_i,b_i\right)\times X\Big)^{\left(k_1^{\rm max}\right)}\Big)^{\left(k_2^{\rm max}\right)}\Big)\cdots\Big)^{\left(k_{b-1}^{\rm max}\right)},
\]
where $k_{i+1}^{\rm max}$ is the degree of the highest derivative of $\left(\cdots\left(\left(\bigtimes_{i\leq k}u_\rho\left(a_i,b_i\right)\right)\times X\right)\cdots\right)^{\left(k_i^{\rm max}\right)}$. Let
\[
    \mathfrak{m}_{X,k}^t\coloneqq \left(\bigsqcup_{i=1}^k \left\{\left[A_i,B_i\right]_\rho,\left[A_i-1,B_i-1\right]_\rho,\cdots,\left[-B_i,-A_i\right]_\rho\right\}\right) \sqcup \mathfrak{m}_X^t
\]
be the multi-segment such that ${\rm Z}\left(\mathfrak{m}_{X,k}^t\right)\cong \bigtimes_{i\leq k} u_\rho\left(a_i,b_i\right)\times X$.
Then by \cite[Thm. 8.1]{ZeleII}, the multiplicity of $\rho|\cdot|^{-A}$ in the cuspidal support of $D_{\rm max}^{b-1}\left(\bigtimes_{i\leq k} u_\rho\left(a_i,b_i\right)\times X\right)$ is
\begin{align*}
n_{\rho|\cdot|^{-A}}\left(\bigtimes_{i=1}^k u_\rho\left(a_i,b_i\right)\times X\right) &\coloneqq \left|\left\{\left[x,y\right]_\rho\in\mathfrak{m}_{X,k}^t\,\big|\,x\geq -B \text{ and } y=-A\right\}\right| \\
&\leq  n + \sum_{j=1}^m\left|\left\{\left[x,y\right]_\rho\in\mathfrak{m}_{L^-_j}\,\big|\,x=A\right\}\right|.
\end{align*}

\vskip 5pt

Recall that $\Hom_{\GL_{N+\ell^-}\left(E\right)}\left(\tau\times L^-_1\times\cdots\times L^-_m, \Sigma\times S^{\times m}\right)\neq0$. Hence, there exists a common irreducible subquotient of $\tau\times L^-_1\times\cdots\times L^-_m$ and $\Sigma\times S^{\times m}$. By Theorem \ref{T:ProdFunctorFullyFaithfulApp}, such a common irreducible subquotient can be written as
\[
    \bigtimes_{i=1}^k u_\rho\left(a_i,b_i\right)\times X \cong Y \times S^{\times m},
\]
where $X$ is an irreducible subquotient of $\left(\bigtimes_{i>k} u_\rho\left(a_i,b_i\right)\right)\times L^-_1\times\cdots\times L^-_m$, and $Y$ is an irreducible subquotient of $\Sigma$. Let $\mathfrak{m}^t_{Y}$ (resp. $\mathfrak{m}^t_{\widetilde{L}^-_j}$, $\mathfrak{m}^t_{\Xi}$) be the multi-segment such that ${\rm Z}\left(\mathfrak{m}^t_Y\right)\cong Y$ (resp. ${\rm Z}\left(\mathfrak{m}^t_{\widetilde{L}^-_j}\right)\cong \widetilde{L}^-_j$, ${\rm Z}\left(\mathfrak{m}^t_{\Xi}\right)\cong \Xi$). 
Then
\[
    \mathfrak{m}_{X,k}^t = \left(\bigsqcup_{j=1}^m \left\{\left[A,B\right]_\rho,\left[A-1,B-1\right]_\rho,\cdots,\left[-B,-A\right]_\rho\right\}\right) \sqcup \mathfrak{m}_Y^t.
\]
We then give a lower bound of $n_{\rho|\cdot|^{-A}}\left(\bigtimes_{i\leq k} u_\rho\left(a_i,b_i\right)\times X\right)$ for such a common subquotient. By \cite[Thm. 7.1]{ZeleII}, we know that $\mathfrak{m}_{Y}^t$ is obtained from $\mathfrak{m}^t_{\widetilde{L}^-_1}\sqcup\mathfrak{m}^t_{\widetilde{L}^-_2}\sqcup\cdots\sqcup\mathfrak{m}^t_{\widetilde{L}^-_m}\sqcup\mathfrak{m}^t_{\Xi}$ by applying finitely many times union-intersection of linked segments. Note that for each segment $\left[x,y\right]_\rho\in\mathfrak{m}_{L^-_j}$ such that $x=A$, by the result of Kret--Lapid (Lemma \ref{L:KretLapidJacFormula}) one knows that $\left[-B,-A\right]_\rho\in\mathfrak{m}^t_{\widetilde{L}^-_j}$. Starting from $\mathfrak{m}^t_{\widetilde{L}^-_1}\sqcup\mathfrak{m}^t_{\widetilde{L}^-_2}\sqcup\cdots\sqcup\mathfrak{m}^t_{\widetilde{L}^-_m}\sqcup\mathfrak{m}^t_{\Xi}$, each time one performs the union-intersection procedure of two linked segments, a segment $\left[x,-A\right]_\rho$ with $x\geq -B$ can only be replaced by $\left[x',-A\right]_\rho$ with $x'\geq x$; the number of such segments is non-decreasing. It follows that
\[
    n_{\rho|\cdot|^{-A}}\left(\bigtimes_{i=1}^k u_\rho\left(a_i,b_i\right)\times X\right) \geq  m + \sum_{j=1}^m\left|\left\{\left[x,y\right]_\rho\in\mathfrak{m}_{L^-_j}\,\big|\,x=A\right\}\right|.
\]
Combining this lower bound with the previous upper bound, we deduce that $m\leq n$ as desired.

\vskip 5pt


\subsection{Proof of Proposition \ref{P:TechnicalCountingMultiApp}: II. non (conjugate) self-dual case} 
\label{sub:proof_of_proposition_p:technicalcountingmultiapp_ii_non_conjugate_self_dual_case}

Next we consider the case that $\rho\not\cong\widetilde{\rho}$ is non (conjugate) self-dual. By Bernstein--Zelevinsky's block decomposition, we can write $\tau$ and $\Xi$ as
\[
    \tau = \tau_{\widetilde{\rho}}\times\tau_\rho\times\tau^\rho \quad \text{and} \quad \Xi = \Xi_{\widetilde{\rho}}\times\Xi_\rho\times \Xi^\rho,
\]
where $\tau_\rho$ (resp. $\tau_{\widetilde{\rho}}$, $\tau^\rho$) is an irreducible representation with cuspidal support contained in $\rho|\cdot|^{A+\Z}$ (resp. contained in $\widetilde{\rho}|\cdot|^{A+\Z}$, disjoint from $\rho|\cdot|^{A+\Z}\cup\widetilde{\rho}|\cdot|^{A+\Z}$), and similarly for $\Xi_\rho$, $\Xi_{\widetilde{\rho}}$ and $\Xi^\rho$. Then we have
\[
    \Hom_{\GL_N\left(E\right)}\left(\tau, \widetilde{L}^-_1\times\widetilde{L}^-_2\times\cdots\times\widetilde{L}^-_m\times\Xi\times L^+_m\times\cdots\times L_2^+\times L_1^+\right)
\]
\[
    \cong \Hom_{\GL_N\left(E\right)}\left(\tau_{\widetilde{\rho}}\times\tau_\rho\times\tau^\rho, \left(\widetilde{L}^-_1\times\cdots\times\widetilde{L}^-_m\times\Xi_{\widetilde{\rho}}\right)\times\left(\Xi_\rho\times L^+_m\times\cdots\times L_1^+\right)\times\Xi^\rho\right).
\]
For the later Hom space to be non-zero, we must have:
\begin{itemize}
    \item $\tau_{\widetilde{\rho}}$ and $\widetilde{L}^-_1\times\cdots\times\widetilde{L}^-_m\times\Xi_{\widetilde{\rho}}$ are representations of the same group $\GL_{N''}\left(E\right)$, and 
    \[
        \Hom_{\GL_{N''}\left(E\right)}\left(\tau_{\widetilde{\rho}},\widetilde{L}^-_1\times\cdots\times\widetilde{L}^-_m\times\Xi_{\widetilde{\rho}}\right)\neq 0;
    \]

    \vskip 5pt

    \item $\tau_\rho$ and $\Xi_\rho\times L^+_m\times\cdots\times L_1^+$ are representations of the same group $\GL_{N'}\left(E\right)$, and 
    \[
        \Hom_{\GL_{N'}\left(E\right)}\left(\tau_\rho,\Xi_\rho\times L^+_m\times\cdots\times L_1^+\right)\cong\Hom_{\GL_{N'+\ell^-}\left(E\right)}\left(\tau_\rho\times L^-_1\times\cdots\times L^-_m, \Xi_\rho\times S^{\times m}\right)\neq 0.
    \]
\end{itemize}
\vskip 5pt
Similar to the argument in the previous subsection, by considering the multiplicity of $\widetilde{\rho}|\cdot|^{-A}$ (resp. $\rho|\cdot|^{-A}$) in the cuspidal support of iterative derivatives of $\tau_{\widetilde{\rho}}$ (resp. irreducible subquotients of $\tau_\rho\times L^-_1\times\cdots\times L^-_m$), we obtain that
\[
    \left|\left\{1\leq i\leq r\,\big|\,\rho_i\cong\widetilde{\rho} \text{ and } \left(a_i,b_i\right)=\left(a,b\right)\right\}\right|\geq \sum_{j=1}^m\left|\left\{\left[x,y\right]_\rho\in\mathfrak{m}_{L^-_j}\,\big|\,x=A\right\}\right|
\]
\[
    \left(\text{resp. } \left|\left\{1\leq i\leq r\,\big|\,\rho_i\cong{\rho} \text{ and } \left(a_i,b_i\right)=\left(a,b\right)\right\}\right|+ \sum_{j=1}^m\left|\left\{\left[x,y\right]_\rho\in\mathfrak{m}_{L^-_j}\,\big|\,x=A\right\}\right|\geq m\right).
\]
These two inequalities imply that $m\leq n$. This completes the proof of Proposition \ref{P:TechnicalCountingMultiApp}.



\vskip 15pt

\bibliographystyle{amsalpha}
\bibliography{NTGGP4ClassicalRef}

\end{document}